\documentclass[11pt]{article}
\usepackage{amsfonts,latexsym,graphicx,amssymb,makeidx,amsthm}
\usepackage{amsmath,amsxtra}

\newcounter{f1} 
\newcounter{f2} 
\newcounter{f3}
\newcounter{f4}
\newcounter{f5} 
\newcounter{f6} 
\newcounter{f7} 
\newcounter{f8}
\newcounter{f9}
\newcounter{f10} 
\newcounter{f11} 
\newcounter{f12} 
\begin{document} 
\title{\sc 
immersivity of the contact line bundle of a complex-contact 
manifold and an application to the automorphism group 
}
\author{Osami {\sc Yasukura}} 
\date{}
\maketitle
\ 
{\bf Abstract.} 
A connected Fano complex-contact manifold is isomorphic to 
the k\"ahlerian C-space of Boothby type 
with a natural complex-contact structure 
corresponding to a non-abelian simple complex Lie algebra 
if the contact line bundle is very ample. 
A.~Beauville relaxed the provision 
to two assumptions that the contact line bundle is generically finite 
and that the automorphism group is reductive. 
We relax the provision to another one 
that the contact line bundle is {\it immersive}, that is, 
the manifold admits a holomorphic immersion into some projective space 
associated with some holomorphic sections of the line bundle. 
As an application, we obtain that 
the automorphism group of a connected compact complex-contact manifold 
with immersive contact line bundle is isomorphic to the automorphism group 
of the corresponding simple complex Lie algebra of rank greater than one, 
which is not connected if and only if its type 
is $A_{n}, D_{n+2}$ for $n > 1$ or $E_{6}$. 
\medskip 
\\ \ 
{\bf 1.~Introduction.}
Let $\mathbb{Z}$ be the group of all integers. 
Put $\mathbb{N} := \{ n \in \mathbb{Z} {;~} n > 0 \}$. 
Let $\mathbb{R}$ be the field of all real numbers, 
$\mathbb{C}$ the field of all complex numbers, 
and $\mathbb{C}^{\times}$ 
the group of all non-zero complex numbers. 
Let $M$ be a complex manifold with 
complex dimension ${\rm n}_{M} \in \mathbb{N} \cup \{ 0 \}$ 
in the sense of C.~Chevalley \cite[Chap.III]{Cc1946} 
(cf. S.~Kobayashi \& K.~Nomizu \cite[p.3]{KN1963}). 
Let 
${\rm Aut}(M)$ 
be the group of all biholomorphic transformations on $M$, 
and ${\cal O}(M)$ 
be the $\mathbb{C}$-linear space of all holomorphic functions on $M$. 
Let $P$ be a holomorphic fiber bundle on $M$ 
with the projection $p: P \to M$. 
For any subset $V$ of $M$, put
$P |_{V} := \cup_{\xi \in V} P |_{\xi}$ 
with 
$P |_{\xi} := p^{-1}(\xi)$. 
For any subset $W$ of $P$, let 
$\Gamma^{\rm cont}_{V}(W)$ 
be the set of all continuous maps 
$\sigma: V \to W$ 
from the topological subspace $V$ of $M$ 
to the topological subspace $W$ of $P$ such that 
$p \circ \sigma = {\rm id}_{V}: 
V \to V; \xi \mapsto \xi$. 
Put 
$\Gamma^{\rm cont}_{\rm loc} (W) 
:= \cup_{V: {\rm open}~{\rm in}~M} 
\Gamma^{\rm cont}_{V} (W)$.  
If $W$ is a complex submanifold of $P$, put 
$\Gamma_{\rm loc} (W) 
:= \cup_{V: {\rm open}~{\rm in}~M} 
\Gamma_{V} (W)$ 
with 
$\Gamma_{V}(W) := \{ \sigma \in \Gamma^{\rm cont}_{V}(W); 
\sigma~{\rm is}~{\rm holomorphic} \}$ 
for any open subset $V$ of $M$. 
For any 
$\sigma \in \Gamma_{\rm loc}^{\rm cont} (P)$ 
(resp. $\Gamma_{\rm loc} (P)$), 
let $V_{\sigma}$ 
be the nonempty open subset of $M$ such that 
$\sigma \in \Gamma^{\rm cont}_{V_{\sigma}} (P)$ 
(resp. $\Gamma_{V_{\sigma}} (P)$). 
When $P$ is a holomorphic principal $\mathbb{C}^{\times}$-bundle 
with the right action $R_{\lambda} (\lambda \in G)$ 
and $\sigma \in \Gamma_{\rm loc} (P)$, 
put 
$\phi_{\sigma}: V_{\sigma} \times \mathbb{C}^{\times} 
\to P$; 
$(\xi, \lambda) 
\mapsto \phi_{\sigma}(\xi, \lambda) 
:= R_{\lambda} \sigma (\xi)$ 
as the local trivialization with the holomorphic fiber coordinate 
$\lambda_{\sigma}: 
P |_{V_{\sigma}} 
\to \mathbb{C}^{\times}; 
\eta \mapsto \lambda_{\sigma}(\eta); 
\eta = \phi_{\sigma}(p(\eta), \lambda_{\sigma}(\eta))$. 
Let $L$ be a holomorphic line bundle on $M$ 
with the projection $\pi_L: L \to M$, 
and $L^{*}$ be the dual line bundle 
with the set  
$L^{*} \times_{M} L := 
\cup_{\xi \in M} (L^{*} |_{\xi} \times L |_{\xi})$ 
of dual pairings and the evaluation map: 
$L^{*} \times_{M} L \to \mathbb{C}; 
(\eta, \eta') \mapsto <\!\!\eta,~\eta'\!\!>_{L} 
:= \eta (\eta')$. 
For the image $0$ of the zero section $0_{M}$ 
of $L$ (resp. $L^{*}$), 
put $Q_{L} := L \backslash 0$ (resp. $P_{L} := Q_{L^{*}}$), 
which is defined as a holomorphic principal 
$\mathbb{C}^{\times}$-bundle on $M$ 
with the canonical projection 
$q_{L}: Q_{L} \to M; 
\eta \mapsto q_{L}(\eta) := \pi_{L}(\eta)$ 
(resp. 
$p_{L}: P_{L} \to M; 
\eta \mapsto p_{L}(\eta) := \pi_{L^{*}}(\eta)$) 
where the right action 
$R_{\lambda}$ ($\lambda \in \mathbb{C}^{\times}$) 
is defined as the restriction of $\lambda$-multiple 
on each $\mathbb{C}$-liner fiber $L |_{\pi_{L}(\eta)}$ 
(resp. $L^{*} |_{\pi_{L^{*}}(\eta)}$) 
to the subset $Q_{L} |_{q_{L}(\eta)}$ 
(resp. $P_{L} |_{p_{L}(\eta)}$) 
of all non-zero elements. 
Put $\rho_{L}: Q_{L} \to Q_{L^{*}}$; 
$<\!\! \rho_{L}(\eta), \eta \!\!>_{L} = 1$ 
for $\eta \in Q_{L}$, 
so that 
$q_{L^{*}} \circ \rho_{L} = {\rm id}_{M}$ 
and 
$\rho_{L} \circ R_{\lambda} = R_{\lambda}^{-1} \circ \rho_{L}$ 
($\lambda \in \mathbb{C}^{\times}$). 
Put 
$\iota_L(s): P_{L} \to \mathbb{C}; 
\eta \mapsto \iota_{L}(s) |_{\eta} 
:= <\!\!\eta,~s(p_{L}(\eta))\!\!>_{L}$ 
for a map $s: M \to L$ 
such that $q_{L} \circ s = {\rm id}_{M}$. 
For a $\mathbb{C}$-linear space ${\cal G}$, 
let 
$P_{\cal G} \mathbb{C} := \{ \mathbb{C} X {;~} X \in {\cal G} \}$ 
be the complex projective space with the projection 
$\pi_{\cal G}: 
{\cal G} \backslash \{ 0 \} \to P_{\cal G} \mathbb{C}; 
X \mapsto \mathbb{C} X$. 
For $N \in \mathbb{N}$, 
let $\mathbb{C}_{N+1}$ be the $\mathbb{C}$-linear space 
of all $(N+1)$-row vectors over $\mathbb{C}$, and put 
$P_{N} \mathbb{C} := P_{\mathbb{C}_{N+1}} \mathbb{C}$. 
For a holomorphic line bundle $L$ on a complex manifold $M$ 
of ${\rm n}_{M} > 0$, 
$L$ is called {\it very ample} 
(resp. {\it immersive}; {\it injective}) 
if and only if 
there exist $N \in \mathbb{N}$ and $s_{i} \in \Gamma_{M} (L)$ 
($i \in \{ 0, \ldots, N \}$) such that 
the projectivization 
$\varphi: M \to P_{N} \mathbb{C}; 
p_{L}(\eta) \mapsto \pi_{\mathbb{C}_{N+1}}(\Phi(\eta))$ 
of $\Phi := (\iota_L(s_0), \ldots, \iota_L(s_N)): 
P_{L} \to \mathbb{C}_{N+1}$ 
is well-defined as an injective immersion 
(resp. an immersion; an injection). 
For $m \in \mathbb{N}$, 
let $L^{\otimes m}$ be the $m$-th tensor power of $L$. 
By definition, 
$L$ is called {\it ample} 
(resp. {\it quasi-immersive}) 
if and only if 
there exists $m \in \mathbb{N}$ 
such that $L^{\otimes m}$ is very ample 
(resp. immersive). 
For a complex manifold $M$, let $T' M$ 
be the holomorphic tangent bundle of $(1,0)$-type 
with the projection 
$\Pi_{M}: T' M \to M; x \mapsto \xi$; 
$x \in T_{\xi}' M := T' M |_{\xi}$ 
with the Lie algebra 
$a(M) := \Gamma_{M}(T' M)$ 
of all infinitesimal automorphisms on $M$ and the anti-canonical line bundle
$K_{M}^{*} := \land^{{\rm n}_{M}} T' M$. 
A compact complex manifold $M$ of ${\rm n}_{M} > 0$ 
is defined to be {\it Fano} if and only if  
$K_{M}^{*}$ is positive in the sense of 
K.~Kodaira \cite{Kk1954}, 
so that $K_{M}^{*}$ is ample by Kodaira 
embedding theorem. 
\\ \
C.~LeBrun \cite {Lc1993} and C.~LeBrun \& S.~Salamon \cite{LS1994} 
simplified the notion of the complex-contact distribution 
in S.~Kobayashi \cite{Ks1959}'s notion of the contact line bundle 
for the complex-contact distribution from the complex-contact structures, 
as a ``non-degenerate distribution of codimension one'' 
which improves V.I.~Arnold \cite[App.4]{Avi1978}, as follows: 
Let $E$ 
be any holomorphic vector subbundle of $T' M$ 
having a complex fiber dimension 
${\rm r}_E \in \mathbb{N} \cup \{ 0 \}$ 
and the quotient vector bundle 
$L_{E} := T' M / E$ 
with 
$\varpi_E: T' M \to L_{E}; 
x \mapsto x + E |_{\Pi_{M}(x)}$. 
Put 
$\alpha_{E}: 
E \times_{M} E 
\to L_{E}; 
(x_{1}, x_{2}) \mapsto 
\alpha_{E}(x_{1}, x_{2}) 
:= \varpi_E([X_1, X_2]) |_{\xi}$ 
for 
$X_i \in \Gamma_{\rm loc}(E)$ 
with 
$X_i |_{\xi} = x_{i}~(i = 1, 2)$ 
as O'Neill's {\it involutor} of $E$, 
by which {\it the non-degenerator} of $E$ 
is defined as follows. 
\[
\alpha_{E}^{\perp} 
:= \cup_{\xi \in M} \alpha_{E}^{\perp} |_{\xi}; 
\alpha_{E}^{\perp} |_{\xi} := 
\{ x \in E |_{\xi} {;~}\alpha_{E}(x, y) = 0
~{\rm for}~{\rm all}~y \in E |_{\xi} \}. 
\]
Then 
$\alpha_{E}^{\perp} = E$ 
if and only if $E$ is involutive. 
When $\alpha_{E}^{\perp} = 0$, 
$E$ is said to be {\it non-degenerate}. 
When 
$E$ is non-degenerate and 
that ${\rm r}_{E} = {\rm n}_{M} - 1$, 
$E$ is said to be {\it a c-distribution} on $M$. 
Then ${\rm n}_{M} = 2 {\rm k}_{M} - 1$ 
with some ${\rm k}_{M} \in \mathbb{N}$. 
In this case, 
$(M, E)$ is said to be a {\it c-manifold} 
with {\it the c-line bundle} $L_{E}$. 
By \cite{Ks1959}, 
$L_{E}^{\otimes {\rm k}_{M}} = K_{M}^{*}$ 
up to isomorphisms  (cf. \cite {Lc1993, LS1994}). 
Note that a c-manifold $(M, E)$ of ${\rm n}_{M} = 1$ 
is nothing but a Riemann surface $M$ with $E = 0$ and $L_{E} = T' M = K_{M}^{*}$, 
and that a c-manifold $(M, E)$ of ${\rm n}_{M} > 1$ 
is just a complex-contact manifold with {\it the contact line bundle} $L_{E}$ 
of a {\it complex-contact distribution} $E$ 
which comes from some complex-contact structures 
in the sense of \cite{Ks1959} (cf. Lemma 4.0 (A) (b)). 
Two c-manifolds $(M_{i}, E_{i})$ ($i = 1, 2$) 
are called {\it isomorphic:} 
$(M_{1}, E_{1}) \cong (M_{2}, E_{2})$ 
if and only if 
there exists an {\it isomorphism} between them, 
i.e., a biholomorphism between them 
transferring their c-distributions from one another. 
An isomorphism from a c-manifold to itself is said 
to be an {\it automorphism} of the c-manifold. 
Let ${\rm Aut}(M, E)$ be the group of all automorphisms of $(M, E)$, 
and $a(M, E)$ be the complex Lie algebra of all infinitesimal 
automorphisms of $(M, E)$. 
A c-manifold $(M, E)$ is called {\it homogeneous} 
if and only if $M$ is one orbit of ${\rm Aut}(M, E)$. 
W.M.~Boothby \cite{Bwm1961, Bwm1962} classified a connected, 
simply connected compact homogeneous complex-contact manifold, 
which is called a {\it k\"ahlerian C-space of Boothby type} 
after T.~Nitta \& M.~Takeuchi \cite{NT1987} up to biholomorphisms, 
the list of which is nothing but 
the classification list up to isomorphisms 
since it is known that the complex-contact distributions on a k\"ahlerian C-space 
of Boothby type are unique up to biholomorphisms 
by \cite[Thm.1.7]{NT1987} (cf. \cite[2.2]{LS1994}). 
A k\"ahlerian C-space of Boothby type or $P_{1} \mathbb{C}$ 
is biholomorphic to the {\it adjoint variety} $X_{\mathfrak{g}}$ 
for a simple complex Lie algebra $\mathfrak{g}$ 
of ${\rm n}_{\mathfrak{g}} > 1$ 
in the sense of S.~Mukai \cite{Ms1994}, 
which admits a natural c-distribution 
$E_{\theta_{\mathfrak{g}}}$. 
For a $\mathbb{C}$-linear space ${\cal G}$, 
let $GL_{{\cal G}} \mathbb{C}$ (resp. $PGL_{{\cal G}} \mathbb{C}$) 
be the group of all $\mathbb{C}$-linear transformations on ${\cal G}$ 
(resp. all projective transformations on $P_{{\cal G}} \mathbb{C}$). 
Let 
$\mathfrak{S}_{k}$ 
be the symmetric group of $\{ 1, ..., k \}$, 
so that $\mathfrak{S}_{1} \cong \{ 1 \}$ and 
$\mathfrak{S}_{2} \cong \mathbb{Z}_{2} := \{ \pm 1 \}$ 
in the group $\mathbb{C}^{\times}$. 
Here is the main result. 
\medskip \\ \
\textsc{Theorem 1.1}. 
{\it 
Let $(M, E)$ be a connected compact c-manifold of 
${\rm n}_{M} > 0$ with $\mathfrak{g} := a(M, E)$. 
Assume that $L_{E}$ is immersive. 
Then $\mathfrak{g}$ is a simple complex Lie algebra 
of ${\rm n}_{\mathfrak{g}} > 1$, 
and that 
$(M, E) \cong (X_{\mathfrak{g}}, E_{\theta_{\mathfrak{g}}})$. 
Moreover, 
both $L_{E}$ and $K_{M}^{*}$ are very ample, 
so that $M$ is Fano. 
}
\medskip \\ \
\textsc{Theorem 1.2}. 
{\it 
Let 
$\mathfrak{g}$ be a simple complex Lie algebra of 
${\rm n}_{\mathfrak{g}} > 1$. 
}
\smallskip \\ \
(1) 
{\it 
For $\beta \in GL_{\mathfrak{g}} \mathbb{C}$, 
let $[\beta] \in PGL_{\mathfrak{g}} \mathbb{C}$ 
be defined such that 
$[\beta] (\pi_{\mathfrak{g}} X) := \pi_{\mathfrak{g}} (\beta X)$ 
for all $X \in \mathfrak{g} \backslash \{ 0 \}$. 
Then 
${\rm Aut}(\mathfrak{g}) \to {\rm Aut}(X_{\mathfrak{g}}, E_{\theta_{\mathfrak{g}}}); 
\beta \mapsto [\beta] |_{X_{\mathfrak{g}}}$ 
is a complex Lie group isomorphism. 
}
\smallskip \\ \
(2)
{\it 
Put $G := {\rm Aut}(X_{\mathfrak{g}}, E_{\theta_{\mathfrak{g}}})$. 
Let $G^{\circ}$ be the identity connected component of $G$. 
Then the quotient group $G/G^{\circ}$ 
is isomorphic to one of the following groups 
for $n > 1$: 
}
\[
G/G^{\circ} \cong \left\{ 
\begin{array}{cc} 
\{ 1 \} & ({\rm if}~\mathfrak{g}~{\rm is~of~type}
~A_{1}, C_{n}, B_{n+1}, G_{2}, F_{4}, E_{7}~{\rm or}~E_{8}), \\
\mathbb{Z}_{2} & 
({\rm if}~\mathfrak{g}~{\rm is~of~type}~A_{n}, D_{n+3}~{\rm or}~E_{6}), \\
\mathfrak{S}_{3}
& ({\rm if}~\mathfrak{g}~{\rm is~of~type}~D_{4}).~~~~~~~~~~~~~~~ \\ 
\end{array} 
\right.
\]
\smallskip \\ \
\textsc{Corollary 1.3.}
{\it 
Let $M$ be a positive quaternionic-K\"ahler manifold 
with the Lie algebra $i(M)$ of all infinitesimal isometries on $M$ 
and the complex-contact twistor space $(Z, E)$. 
Put $\mathfrak{g} := i(M) \otimes_{\mathbb{R}} \mathbb{C}$. 
Assume that $L_{E}$ is immersive. Then 
$\mathfrak{g}$ 
is a simple complex Lie algebra of ${\rm rank}(\mathfrak{g}) > 1$ 
and that 
$(X_{\mathfrak{g}}, E_{\theta_{\mathfrak{g}}}) \cong (Z, E)$, 
so that 
${\rm Aut}(\mathfrak{g})$ is isomorphic to ${\rm Aut}(Z, E)$. 
}
\smallskip \\ \ 
A connected 
smooth Riemannian real $4n$-manifold $M$ 
is called {\it quaternion-K\"ahler} 
when there exists a fixed quaternion-K\"ahler structure on $M$ 
\cite[p.146, \S 2]{NT1987} (cf. H.~Tasaki \cite{Th1986}), 
i.e., 
the holonomy group of $M$ is contained in $Sp(n) Sp(1)$ 
(cf. \cite[p.396, p.403]{Bal1987}), 
so that $M$ is orientable.  
After \cite{Lc1993, LS1994}, 
a {\it positive quaternionic-K\"ahler manifold} 
is defined to be a geodesically complete quaternion-K\"ahler real 
$4n$-manifold $M$ of positive scalar curvature 
such that either $n > 1$ 
or $M$ is an Einstein half-conformally flat 4-manifold 
(cf. \cite{Ss1982}, \cite{NT1987}), 
which is compact (cf. S.~Salamon \cite[p.103]{Ss1989}). 
J.A.~Wolf \cite{Wja1965} 
constructed a positive quaternionic-K\"ahler manifold 
$W_{\mathfrak{u}}$ for a simple compact Lie algebra $\mathfrak{u}$ of 
${\rm rank}(\mathfrak{u}) > 1$ as a Riemannian symmetric space of certain type 
such that the Lie algebra $i(W_{\mathfrak{u}})$ 
of all infinitesimal isometries on 
$W_{\mathfrak{u}}$ is isomorphic to $\mathfrak{u}$, 
which is called {\it the Wolf space} for $\mathfrak{u}$ 
after D.V.~Alekseevskii \cite{Adv1968} (cf. \cite[14.53, 14.54]{Bal1987}), 
which is also characterized by \cite{Th1986} 
as {\it a semisimple quaternionic symmetric space of compact type}. 
S.~Salamon \cite{Ss1982} constructed 
{\it the twistor space} $(Z, E)$ of 
any positive quaternionic-K\"ahler $4n$-manifold $M$
as a Fano complex-contact manifold 
(cf. \cite{Wja1965}, 
N.J.~Hitchin \cite{Hnj1981}, 
L.~B\'erard Bergery \cite[14.80]{Bal1987}). 
After \cite{Hnj1981} (cf. S.~Mukai \cite{Ms1989}), 
Th.~Friedrich \& H.~Kurke \cite{FK1982}, 
\cite{Ss1982}, \cite{NT1987}, 
Y.S.~Poon \& S.~Salamon \cite{PS1991} 
and \cite{Lc1993, LS1994}, 
it is a problem in the literature whether 
{\it 
any positive quaternionic-K\"ahler manifold 
is homothetic to some Wolf space or not}  
(cf. \cite[13.30 (1)]{Bal1987}, \cite[7.8]{Ss1989}). 
This paper comes from a study about 
complex-contact geometry on this problem. 
\begin{center}
{\bf Table of Contents}
\end{center}
\
1. Introduction. 
\\ \ 
2. 
Geometric topology and Euler operator. 
\\ \
3. 
Quasi-immersivity of holomorphic line bundles. 
\\ \
4. 
Local symplectifications and infinitesimal automorphisms. 
\\ \
5. 
The adjoint varieties. 
\\ \
6. 
Proof of Theorems 1.1 and 1.2. 
\\ \
7. 
Proof of Corollary 1.3. 
\\ \
In \S 2, 
{\it geometric topology} on a complex manifold $M$ 
is studied by virtue of O.~Forster \cite[pp.28-29]{Fott1981} 
and K.~Fritzsche \& H.~Grauert \cite[pp.37-38]{FG2002}. 
In general, 
{\it Euler operator} is defined and studied 
for a holomorphic principal $\mathbb{C}^{\times}$-bundle $P$ on $M$. 
\\ \
In \S 3, 
immersivity and quasi-immersivity of a holomorphic line bundle $L$ on $M$ 
are studied by homogeneous functions of degree $m \in \mathbb{N}$ 
on $P_{L}$. 
\\ \
In \S 4, 
the notion of a local symplectification $(W, \theta)$ of degree 
$\epsilon_{\theta} \in \mathbb{C}^{\times}$ 
is defined in a holomorphic principal $\mathbb{C}^{\times}$-bundle $P$ 
on a complex manifold $M$ 
to unify symplectifications given in 
\cite[p.145]{Bwm1961}, \cite{Avi1978}, T.~Oshima \& H.~Komatsu \cite{OK1977} 
and \cite{NT1987}. 
Moreover, the prolongation of automorphisms of a c-manifold is described. 
\\ \
In \S 5, 
N.~Tanaka \cite{Tn1976}'s theory of 
$\mathbb{Z}$-graded simple Lie algebras 
is studied as Lemma 5.0, 
by which the adjoint variety $X_{\mathfrak{g}}$ 
for a simple complex Lie algebra $\mathfrak{g}$ 
of ${\rm n}_{\mathfrak{g}} > 1$ 
is realized as Lemma 5.1 (A) 
(cf. \cite{Yo1998}, H.~Kaji, M.~Ohno \& O.Y. \cite{KOY1999}, 
H.~Kaji \& O.Y. \cite{KY2004}). 
As Lemma 5.1 (B), 
$X_{\mathfrak{g}}$ 
is characterized by the use of 
W-L.~Chow \cite{Cwl1949}'s theorem applied for 
compact complex submanifolds of $P_{n} \mathbb{C}~(n \in \mathbb{N})$ 
and A.~Borel \cite{Ba1991}'s fixed point theorem 
applied for complex projective algebraic varieties 
by virtue of T.~Ohta \& T.~Nishiyama \cite[p.28]{ON2015}. 
In Proposition 5.2 (D) (cf. \cite{Yo1998}), 
it is obtained that 
{\it 
there exists a natural c-distribution on $X_{\mathfrak{g}}$ 
such that the c-line bundle is very ample, so that 
$X_{\mathfrak{g}}$ is Fano and simply connected}. 
\\ \
In \S 6, 
we prove Lemma 6.0 describing 
{\it 
the equivalence of 
the immersivity of the c-line bundle and the infinitesimal transitivity 
of the canonical symplectification} 
after 
\cite{Yo1993, Yo1998} 
obtained independently from Wisniewskii \cite{Ba1998}. 
By virtue of A.~Swann \cite{Sa2001}'s moment map, 
we obtain Lemmas 6.1 and 6.2. 
By the use of A.J.~Sommes \cite{Saj1973}'s 
fixed point theorem applied for 
a simply conneced compact k\"ahlerian manifold, 
we obtain Theorems 1.1 and 1.2 as Propositon 6.3. 
Especially, \cite[(2.5)]{Wja1965}'s 
observation on the infinitesimal automorphisms 
is improved as Corollary 6.4. 
In \S 7, 
we prove Corollary 1.3 from Proposition 6.3 
by the use of \cite{NT1987}'s theorems of prolongation and Matsushima type 
by virtue of Lemma 7.0 (A) generalizing \cite[p.248, Thm.VI.4.6 (1)]{KN1963}. 
\smallskip \\ \
{\bf 2. 
Geometric topology and Euler operator.} 
\smallskip \\ \
\textsc{Lemma 2.0.} 
(A) 
{\it 
Let $A$ be a topological space (resp. a complex manifold) 
with open subsets $A_{j} (j \in J)$ of $A$ 
such that $A = \cup_{j \in J} A_{j}$. 
Assume that $C$ is a topological space 
(resp. a complex manifold) 
with continuous (resp. holomorphic) maps 
$\beta_{j}: A_{j} \to C$ $(j \in J)$ 
such that 
$\beta_{j} = \beta_{j'}$ on $A_{j} \cap A_{j'}$ 
for any two $j, j' \in J$. 
Then 
$\beta := \cup_{j \in J} \beta_{j}: A \to C; 
\beta(\eta) := \beta_{j}(\eta)$ $(\eta \in A_{j})$ 
is well-defined as a continuous (resp. holomorphic) map. 
}
\\ \
(B) (\cite[p.252, Lem.1]{KN1963}) 
{\it 
Let $W$ be a connected complex manifold with $f_{1}, f_{2} \in {\cal O}(W)$ 
admitting a non-empty open subset $U$ of $W$ 
such that $f_{1} |_{U} = f_{2} |_{U}$. 
Then $f_{1} = f_{2}$. 
}
\\ \
(C) (cf. O.~Forster \cite[p.29, 4.22.Thm.]{Fott1981}) 
{\it 
Let $A$ be a compact Hausdorff space, 
and $B$ be a Hausdorff space. 
Assume that 
$\nu: A \to B$ is a locally homeomorphic surjection such that 
$\nu^{-1}(b)$ is a discrete subset of $A$ for every $b \in B$. 
Then $\nu$ is a covering map. 
}
\\ \
(D) 
{\it 
Let $A$ be a compact smooth manifold, 
and $B$ be a smooth manifold. 
Assume that 
$\nu: A \to B$ is a surjective smooth immersion. 
Then $\nu$ is a covering map. 
}
\smallskip \\ \
{\it Proof}. 
(A) 
$\beta$ is well-defined as a map. 
Let $U$ be any open subset of the topological space $C$. 
$\beta^{-1}(U) = \beta^{-1}(U) \cap  (\cup_{j \in J} A_{j}) 
= \cup_{j \in J} \beta^{-1} (U) \cap A_{j} 
= \cup_{j \in J} \beta_{j}^{-1}(U)$. 
For each $j \in J$, 
$\beta_{j}^{-1}(U)$ is open in the open subset $A_{j}$ of $A$, 
so that $\beta_{j}^{-1}(U)$ is open in $A$, 
so is $\cup_{j \in J} \beta_{j}^{-1}(U) = \beta^{-1}(U)$. 
Hence, $\beta$ is continuous. 
If $\beta_{j}$ is holomorphic, 
then $\beta$ is holomorphic by Cauchy-Riemann equation on $A_{j}$. 
\\ \ 
(B) 
This is a special result of a more generalized 
lemma \cite[p.252, Lem.1]{KN1963} 
given for $\mathbb{R}$-valued real analytic functions $f_{1}, f_{2}$ 
defined on a real analytic manifold. 
\\ \
(C) 
(\cite[p.28, 4.21; p.29, 4.22]{Fott1981}) 
Take any $b \in B$. 
Since $B$ is a Hausdorff space, $\{ b \}$ is closed in $B$. 
Hence, 
$\nu^{-1}(b)$ is compact as a closed subset of the compact space $A$. 
Since $\nu$ is surjective, 
$\nu^{-1}(b)$ is a nonempty discrete subset. 
Then there exist $n \in \mathbb{N}$ and a nonempty finite subset 
$\{ a_{1}, \ldots, a_{n} \}$ of $A$ 
such that $a_{i} \ne a_{j}$ for $i \ne j$ 
and that 
$\nu^{-1}(b) = \{ a_{1}, \ldots, a_{n} \}$. 
For every $j \in \{ 1, \ldots, n \}$ 
there exists an open neighbourhood $W_{j}$ of $a_{j}$ 
in $A$ and an open neighbourhood $U_{j}$ of $b$ in $B$ 
such that the restriction 
$\nu_{j} := \nu |_{W_{j} \to U_{j}}$ 
is a homeomorphism, where $W_{j}$'s can be taken 
to be pairwise disjoint since $A$ is Hausdorff 
and $\{ a_{1}, \ldots, a_{n} \}$ is a finite set. 
Put $V := \cup_{j = 1}^{n} W_{j}$, 
which is an open subset of $A$ such that 
$\nu^{-1}(b) \subseteq V$, 
so that 
$\nu^{-1}(b) \cap V^{c} = \emptyset$ 
for $V^{c} := A \backslash V$. 
Since $V^{c}$ is a closed subset of the compact space $A$, 
$V^{c}$ is compact, so that 
$\nu(V^{c})$ 
is compact 
in the Hausdorff space $B$. 
Then 
$\nu(V^{c})$ is closed in $B$. 
Put 
$U_{0} := B \backslash \nu(V^{c})$, which is open in $\nu(V)$ 
such that $\nu^{-1}(U_{0}) \subseteq V$. 
Put $U := \cap_{i = 0}^{n} U_{i}$. 
Then 
$b \in U \subseteq \cup_{j = 1}^{n} U_{j}$ 
such that 
$\nu^{-1}(U) \subseteq \cup_{j = 1}^{n} W_{j}$. 
Put 
$V_{j} := W_{j} \cap \nu^{-1}(U)$, 
so that $V_{j}$'s are disjoint and that 
$\nu^{-1}(U) = \cup_{j = 1}^{n} V_{j}$. 
Moreover, 
$\nu(V_{j}) = \nu(W_{j}) \cap U = U_{j} \cap U = U$, 
so that 
$\nu |_{V_{j} \to U} = \nu_{j} |_{V_{j} \to U}$ 
is a homeomorphism.
\\ \
(D) 
Because of $\nu(A) = B$, 
the topological dimension of $A$ 
is non less than the one of $B$ 
by Harish-Chandra 
(e.g. S.~Helgason \cite[p.345, Lem.VII.12.3]{Hs1978}), 
so that ${\rm n}_{A} \geq {\rm n}_{B}$ 
by W.~Hurewicz \& H.~Wallman 
(e.g. \cite[p.344, Ch.VII, \S 12]{Hs1978}). 
Since $\nu$ is an immersion, 
${\rm n}_{A} \leq {\rm n}_{B}$, 
so that the real differential map 
$\nu_{*a}: T A |_{a} \to T B |_{\nu(a)}$ 
of $\nu$ is an $\mathbb{R}$-linear isomorphism 
at every $a \in A$. 
Because of the constant rank map theorem 
(e.g. \cite[p.86, Thm.I.15.5]{Hs1978}), 
any $(b, a) \in B \times \nu^{-1}(a)$ 
admits an open neighbourhood 
$U_{a}$ of $a$ in $A$ and an open neighbourhood 
$V_{b}$ of $b$ in $B$ 
such that the restriction 
$\nu_{U_{a}} := \nu |_{U_{a} \to V_{b}}$ 
of $\nu$ is a smooth diffeomorphism. 
Especially, 
$U_{a} \cap \nu^{-1}(b) = \nu_{U_{a}}^{-1}(b) = \{ a \}$. 
Hence, 
$\nu$ is locally homeomorphic 
and that $\nu^{-1}(b)$ is discrete in $A$ 
for every $b \in B$. 
Note that $A$ is a compact Hausdorff space, and that 
$B$ is a Hausdorff space, 
so that 
$\nu$ is a covering map because of the claim (C).
~\qed 
\smallskip \\ \
For any set $W$ 
and any function $f: W \to \mathbb{C}$, put 
$W_{f} := \{ \eta \in W {;~} f(\eta) = 0 \}$ 
and $W^{f} := W \backslash W_{f}$. 
As preliminaries, 
we relax Riemann extension theorem \cite[I.8.2]{FG2002} 
and a proof of the connectedness theorem of 
\cite[IV.1.6]{FG2002} 
on a domain in $\mathbb{C}_{N}$ 
to one's on a complex manifold as follows. 
\smallskip \\ \
\textsc{Lemma 2.1.} 
{\it 
Let $W$ be a 
complex manifold 
such that $W \ne \emptyset$ 
with some $f \in {\cal O}(W)$. 
Let ${\rm cl}_{W}(W^{f})$ be the closure of $W^{f}$ in $W$. 
Then one has the following results.
} 
\\ \
(A) 
{\it 
Assume that $W$ is connected. 
Then 
$W = {\rm cl}_{W}(W^{f})$ if and only if $W^{f} \ne \emptyset$. 
}
\\ \
(B) 
{\it 
Assume that $W^{f} \ne \emptyset$, 
$g \in {\cal O}(W^{f})$ 
and that 
there exists an open covering ${\cal U}$ of $W$ 
such that the restriction $g |_{U^{f}}$ to $U^{f}$ 
is bounded for any $U \in {\cal U}$ such as $U^{f} \ne \emptyset$. 
Then there exists unique $h \in {\cal O}(W)$ such that $h |_{W^{f}} = g$. 
}
\\ \
(C) 
{\it 
Assume that $W$ is connected. 
Then $W^{f}$ is connected. 
}
\smallskip \\ \
{\it Proof}. 
(A) 
Assume that $W = {\rm cl}_{W}(W^{f})$. 
Suppose that $W^{f} = \emptyset$. 
Then $W = \emptyset$, 
which contradicts with $W \ne \emptyset$. 
Thus $W^{f} \ne \emptyset$. 
Conversely, 
assume that $W \ne {\rm cl}_{W}(W^{f})$. 
Put $U := W \backslash {\rm cl}_{W}(W^{f}) \ne \emptyset$. 
Then $U$ is an open subset of $W$ 
and that 
$W_{f} := W \backslash W^{f} 
\supseteq W \backslash {\rm cl}_{W}(W^{f}) = U$, 
so that $f |_{U} = 0$. 
By Lemma 2.0 (B), 
$W = W_{f}$ since $W$ is connected. 
Then 
$W^{f} = W \backslash W_{f} = \emptyset$, 
as required. 
\\ \ 
(B1) Let $W = \cup_{\lambda \in \Lambda} W_{\lambda}$ 
be the disjoint union of connected components. 
By assumption, there exists 
$\lambda \in \Lambda$ such that 
$W_{\lambda}^{f} \ne \emptyset$. 
Take any $\lambda \in \Lambda$ 
such that 
$W_{\lambda}^{f} \ne \emptyset$. 
{\it 
Claim. 
$U^{f} \ne \emptyset$ for any non-empty open subset 
$U$ of $W_{\lambda}$:} 
{\it In fact}, 
suppose that there exists a non-empty open subset $U$ of $W_{\lambda}$ 
such that $U^{f} = \emptyset$. 
Then 
$(W_{\lambda})_{f} \supseteq U_{f} := U \backslash U^{f} = U \ne \emptyset$, 
so that $f |_{U} = 0$. 
Since $W_{\lambda}$ is connected, 
$f |_{W_{\lambda}} = 0$ by Lemma 2.0 (B), 
so that $W_{\lambda}^{f} = \emptyset$, 
{\it a contradiction}. 
\\ \
(B2) Take any $\lambda \in \Lambda$. 
When $W_{\lambda}^{f} = \emptyset$: 
Put $h_{\lambda}: W_{\lambda} \to \mathbb{C}; \eta \mapsto 0$. 
When $W_{\lambda}^{f} \ne \emptyset$: 
Put 
${\cal U}_{\lambda} := \{ U' := U \cap W_{\lambda}; U \in {\cal U}, 
U \cap W_{\lambda} \ne \emptyset \}$. 
By $W = \cup_{U \in {\cal U}} U$, 
$W_{\lambda} = \cup_{U' \in {\cal U}_{\lambda}} U'$. 
Take an open refinement ${\cal U}_{\lambda}'$ of ${\cal U}_{\lambda}$ 
such that 
$W_{\lambda} = \cup_{U \in {\cal U}_{\lambda}'} U$ 
and that each $U \in {\cal U}_{\lambda}'$ 
is biholomorphic to a non-empty connected open subset 
$\hat{U}$ of $\mathbb{C}_{{\rm n}_{W}}$. 
Take any $U \in {\cal U}_{\lambda}'$. 
By the step (B1), $U^{f} \ne \emptyset$. 
The restriction $g |_{U^{f}}$ is bounded, 
since there exists $\widetilde{U} \in {\cal U}$ 
such that $U \subseteq \widetilde{U}$ 
on which the restriction $g |_{\widetilde{U}^{f}}$ is bounded by assumption. 
By Riemann extension theorem on $\hat{U}$ \cite[p.38, I.8.2]{FG2002}, 
there exists $h_{U} \in {\cal O}(U)$ 
such that $h_{U} |_{U^{f}} = g |_{U^{f}}$. 
For any $U_{1}, U_{2} \in {\cal U}_{\lambda}'$ 
such as $U_{1} \cap U_{2} \ne \emptyset$, 
let $U_{1} \cap U_{2} 
= \cup_{\lambda'' \in \Lambda''} U_{\lambda''}$ 
be the disjoint union of non-empty connected components of $U_{1} \cap U_{2}$. 
For each $\lambda'' \in \Lambda''$, 
$U_{\lambda''}^{f} \ne \emptyset$ by the step (B1). 
Moreover, 
$h_{U_{1}} |_{U_{\lambda''}^{f}} 
= g |_{U_{\lambda''}^{f}} 
= h_{U_{2}} |_{U_{\lambda''}^{f}}$. 
Because of 
$U_{\lambda''} = {\rm cl}_{U_{\lambda''}}(U_{\lambda''}^{f})$ 
by the claim (A), one has that 
$h_{U_{1}} |_{U_{\lambda''}} = h_{U_{2}} |_{U_{\lambda''}}$ 
since they are continuous on $U_{\lambda''}$ in $U_{1} \cap U_{2}$. 
Hence, 
$h_{U_{1}} |_{U_{1} \cap U_{2}} = h_{U_{2}} |_{U_{1} \cap U_{2}}$. 
By Lemma 2.0 (A), put 
$h_{\lambda} := \cup_{U \in {\cal U}_{\lambda}'}~h_{U} \in {\cal O}(W_{\lambda})$. 
Then 
$h_{\lambda} |_{W_{\lambda}^{f}} 
= \cup_{U \in {\cal U}_{\lambda}'} h_{U} |_{U^{f}} 
= \cup_{U \in {\cal U}_{\lambda}'} g |_{U^{f}} 
= g |_{W_{\lambda}^{f}}$. 
By Lemma 2.0 (A), 
put $h := \cup_{\lambda} h_{\lambda} \in {\cal O}(W)$. 
Then 
$h |_{W^{f}} 
= \cup_{\lambda \in \Lambda}~h_{\lambda} |_{W_{\lambda}^{f}} 
= \cup_{\lambda \in \Lambda}~g |_{W_{\lambda}^{f}} 
= g |_{W^{f}}$. 
\\ \
(C) 
Suppose that $W^{f}$ is not connected, 
so that $W^{f} \ne \emptyset$ 
and that there exist two non-empty open subsets 
$U_{1}, U_{2}$ of $W^{f}$ 
such that $W^{f} = U_{1} \cup U_{2}$ 
and $U_{1} \cap U_{2} = \emptyset$. 
Put $g (\eta) := i$ if $\eta \in U_{i}$ ($i = 1, 2$). 
Then $g \in {\cal O}(W^{f})$, 
which is bounded. 
By the claim (B),  
there exists $h \in {\cal O}(W)$ 
such that $h |_{W^{f}} = g$. 
By the claim (A), $W = {\rm cl}_{W}(W^{f})$, 
so that any $\eta \in W$ 
admits a complex coordinates neighbourhood $U_{\eta}$ of $\eta$ 
in $W$ biholomorphic to an open subset $\hat{U}_{\eta}$ 
of $\mathbb{C}_{{\rm n}_{W}}$ 
and a sequence $\{ \eta_{n} \}$ in $W^{f} \cap U_{\eta}$ 
such that $\lim_{n \to \infty} \eta_{n} = \eta$ 
with respect to a metric on $U_{\eta}$ induced from 
a natural metric on $\hat{U}_{\eta}$ in $\mathbb{C}_{{\rm n}_{W}}$. 
Since $h$ is continuous, 
$\mathbb{C} \ni h(\eta) = \lim_{n \to \infty} h(\eta_{n})
= \lim_{n \to \infty} g(\eta_{n}) \in \{ 1, 2 \}$. 
Then $W = h^{-1}(1) \cup h^{-1}(2)$ 
and $h^{-1}(1) \cap h^{-1}(2) = \emptyset$. 
Since $h^{-1}(1)$ and $h^{-1}(2)$ are closed in $W$, 
$h^{-1}(1) = W \backslash h^{-1}(2)$ 
and $h^{-1}(2) = W \backslash h^{-1}(1)$ 
are open in $W$. 
Since $h^{-1}(i) \supseteq g^{-1}(i) = U_{i} \ne \emptyset$ 
($i = 1, 2$), $W$ is not connected, 
which contradicts with the assumption that $W$ is connected.~\qed 
\smallskip \\ \
Let $P$ 
be a holomorphic principal $\mathbb{C}^{\times}$-bundle on 
a complex manifold $M$ 
with the right action 
$R_{\lambda}$ ($\lambda \in \mathbb{C}^{\times}$) 
and the canonical projection 
$p: P \to M$. 
For $\eta \in P$ and any $S \subseteq P$, 
put 
$\breve{\eta}: 
\mathbb{C}^{\times} \to P |_{p(\eta)}; 
\lambda \mapsto R_{\lambda} \eta$, 
$\breve{\mathbb{C}}^{\times}_{\eta, S} 
:= \breve{\eta}^{-1}(S |_{p(\eta)})$; 
$\hat{\eta}: 
\mathbb{C}^{\times} \to P |_{p(\eta)}; 
\lambda \mapsto R_{\lambda}^{-1} \eta$ 
and 
$\hat{\mathbb{C}}^{\times}_{\eta, S} 
:= \hat{\eta}^{-1}(S |_{p(\eta)})$. 
Let ${\cal U}(P)$ be the set of all open subsets of $P$, 
and put 
${\cal O}_{P} := \cup_{W \in {\cal U}(P)}~{\cal O} (W)$, 
on which we define the {\it Euler operator} ${\cal E}_{P}$ of $P$ 
as follows. 
\[ 
{\cal E}_{P}: {\cal O}_{P} \to {\cal O}_{P}; 
f \mapsto ({\cal E}_{P} f)(\eta) 
:= \left. \frac{d}{d \lambda} 
(f \circ \breve{\eta}) \right|_{\lambda = 1} 
= \left. \frac{d}{d z} (f (R_{e^{z \cdot 1}} (\eta))) 
\right|_{z = 0,} 
\]
which is identified with (1,0)-part of 
{\it the fundamental vector field of $P$ 
with respect to $1 \in \mathbb{C} = {\it Lie}~\mathbb{C}^{\times}$} 
(cf. \cite[pages 4, 42 and 51]{KN1963}), 
so that 
${\cal E}_{P} \in \Gamma_{M} (T' P)$, 
${\cal E}_{P} |_{\eta} \ne 0 (\eta \in P)$ 
and 
$p_{*'}~{\cal E}_{P} = 0$, 
where 
$p_{*'}: T' P \to T' M$ 
means the restriction of the complexified total differential 
$p_{*}: T P \otimes_{\mathbb{R}} \mathbb{C} \to T M \otimes_{\mathbb{R}} \mathbb{C}$ 
of $p$ to the $(1, 0)$-tangent vector bundles $T' P$ and $T' M$. 
For any open subset $W$ of $P$, 
the $\mathbb{C}$-linear space of all 
{\it 
homogeneous functions of degree} $\delta \in \mathbb{C}$ 
on $W$ is defined as 
${\cal O}^{\delta}(W) 
:= \{ f \in {\cal O}(W){;~} 
{\cal E}_{P} f = \delta \cdot f \}$. 
Put 
$\mathbb{Z}_{m} := \{ \omega \in \mathbb{C}{;~} \omega^{m} = 1 \}$. 
\smallskip \\ \
\textsc{Lemma 2.2.} 
{\it 
For 
$k \in \mathbb{Z}$ 
and 
$m \in \mathbb{N}$, 
one has the following results.
}
\\ \
(A) 
{\it 
Let $U$ be a connected open subset of $\mathbb{C}^{\times}$ 
with $g \in {\cal O}(U)$ such that 
$\lambda \frac{d}{d \lambda} (g(\lambda)^{m}) = k \cdot g(\lambda)^{m}$ 
at all $\lambda \in U$. 
For each $\lambda_{0} \in U$, 
put $c_{0} := g(\lambda_{0})^{m} \cdot \lambda_{0}^{-k}$. 
Then $g(\lambda)^{m} = c_{0} \cdot \lambda^{k}$ 
for all $\lambda \in U$. 
Especially, 
$c_{0} = g(1)^{m}$ if $U \ni 1$. 
}
\\ \
(B) 
{\it 
If $S |_{p(\eta)}$ is an open subset of $P |_{p(\eta)}$ 
for every $\eta \in S$, 
then $\breve{\mathbb{C}}^{\times}_{\eta, S}$ 
and $\hat{\mathbb{C}}^{\times}_{\eta, S}$ 
are open neighbourhoods of $1$ in $\mathbb{C}^{\times}$ 
for every $\eta \in S$. 
If $S |_{p(\eta)}$ is a connected subset of $P$ for every $\eta \in S$, 
then $\breve{\mathbb{C}}^{\times}_{\eta, S}$ 
and $\hat{\mathbb{C}}^{\times}_{\eta, S}$ 
are connected for $\eta \in S$. 
}
\\ \
(C) 
{\it 
Let $W$ be an open subset of $P$. 
Then one has the following results.
} 
\\ \
(a) 
{\it 
${\cal O}^{k/m}(W) = 
\{ f \in {\cal O}(W){;~} f^{m} \in {\cal O}^{k}(W) \}$. 
If moreover 
$W |_{p(\eta)}$ is connected for every $\eta \in W$, 
then 
\[
{\cal O}^{k/m}(W) = \{ f \in {\cal O}(W){;~} 
f(R_{\lambda} \eta)^{m} = \lambda^{k} \cdot f(\eta)^{m} 
~(\eta \in W, \lambda \in \breve{\mathbb{C}}^{\times}_{\eta, W}) \}. 
\]
}
\
(b) 
{\it 
Assume that $W |_{p(\eta)}$ is connected for every $\eta \in W$. 
For any subset $V$ of $M$ and $f \in {\cal O}^{k/m}(W)$, put 
$V^{f^{\sim}} := p((P |_{V})^{f})$. 
Then 
$P^{f} |_{\xi} 
= P |_{\xi}$ 
for any $\xi \in V^{f^{\sim}}$, 
so that 
$(P |_{V})^{f} = P |_{V^{f^{\sim}}}$. 
If $V$ is an open subset (resp. a connected subset) of $M$, 
then $V^{f^{\sim}}$ is an open subset (resp. a connected subset) of $M$. 
} 
\\ \
{\it Proof}. 
(A) 
Put $h(\lambda) := g(\lambda)^{m} \cdot \lambda^{-k}$ 
($\lambda \in U$). 
Since $U$ is open and nonempty, 
$\frac{d}{d \lambda} h(\lambda)$ 
is well-defined. 
At $\lambda \in U$, 
$\lambda \frac{d}{d \lambda} h(\lambda) = 
\lambda (\frac{d}{d \lambda} (g(\lambda)^{m})) \cdot \lambda^{-k} 
- k \cdot g(\lambda)^{m} \cdot \lambda^{-k} = 0$ 
by assumption. 
Hence, 
$\frac{d}{d \lambda} h(\lambda) = 0 = \frac{d}{d \lambda} h(\lambda_{0})$. 
By Lemma 2.0 (B), 
$h(\lambda) = h(\lambda_{0})$ ($\lambda \in U$) 
since $U$ is connected. 
Hence, 
$g(\lambda)^{m} = c_{0} \cdot \lambda^{k}$ 
for $c_{0} = h(\lambda_{0})$. 
Especially, $g(1)^{m} = c_{0}$ if $1 \in U$. 
\\ \ 
(B) 
Take any $\eta \in S$. 
Then 
$1 \in \breve{\mathbb{C}}^{\times}_{\eta, S} 
\cap \hat{\mathbb{C}}^{\times}_{\eta, S}$. 
Since $P$ is a holomorphic principal $\mathbb{C}^{\times}$-bundle, 
$\breve{\eta}$ is biholomorphic. 
If $S |_{p(\eta)}$ is open (resp. connected) in $P |_{p(\eta)}$, 
then 
$\breve{\mathbb{C}}^{\times}_{\eta, S} 
= \breve{\eta}^{-1}(S |_{p(\eta)})$ 
and 
$\hat{\mathbb{C}}^{\times}_{\eta, S} 
= \{ \lambda^{-1}; \lambda \in \breve{\mathbb{C}}^{\times}_{\eta, S} \}$ 
are open (resp. connected) in $\mathbb{C}^{\times}$. 
\\ \ 
(C) (a) 
For $f \in {\cal O}(W)$ and $\eta \in W$, 
${\cal E}_{P} |_{\eta } (f^{m}) 
= \left. \frac{d}{d \lambda} \right|_{\lambda = 1} f(R_{\lambda} \eta)^{m} 
= m f(\eta)^{m-1} \cdot 
\left. \frac{d}{d \lambda} \right|_{\lambda = 1} f(R_{\lambda} \eta) 
= m f(\eta)^{m-1} \cdot {\cal E}_{P} |_{\eta} f$. 
Let 
$f \in {\cal O}(W)$ 
and 
$f^{m} \in {\cal O}^{k}(W)$. 
When $f = 0 \in {\cal O}(W)$: 
${\cal E}_{P} f = 0 = \frac{k}{m} \cdot f$, 
so that 
$f \in {\cal O}^{k/m}(W)$. 
When $f \ne 0 \in {\cal O}(W)$: 
By 
$m f(\eta)^{m-1} \cdot {\cal E}_{P} |_{\eta} f 
= k \cdot f(\eta)^{m}$, 
${\cal E}_{P} |_{\eta} f = \frac{k}{m} \cdot f(\eta)$ 
at $\eta \in W^{f}$. 
Put 
${\cal O}(W) \ni g: W \to \mathbb{C}; \eta \mapsto g(\eta) 
:= {\cal E}_{P} |_{\eta} f - \frac{k}{m} \cdot f(\eta)$, 
so that $g |_{W^{f}} = 0$. 
Let 
$W = \cup_{\lambda \in \Lambda} W_{\lambda}$ 
be the disjoint union of connected open subsets 
$W_{\lambda}$ of $W$. 
By $f \ne 0$, 
there exists $\lambda \in \Lambda$ 
such as $W_{\lambda}^{f} \ne \emptyset$. 
For any $\lambda \in \Lambda$ 
such as $W_{\lambda}^{f} \ne \emptyset$, 
$W_{\lambda} = {\rm cl}({W_{\lambda}}^{f})$ 
by Lemma 2.1 (A), so that $g |_{W_{\lambda}} = 0$ 
because of $g |_{W_{\lambda}^{f}} = 0$. 
By Lemma 2.0 (A), 
${\cal O}(W) \ni g = \cup_{\lambda \in \Lambda} g |_{W_{\lambda}} = 0$, 
i.e., $f \in {\cal O}^{k/m}(W)$. 
Conversely, take any $f \in {\cal O}^{k/m}(W)$. 
For any $\eta \in W$, 
${\cal E}_{P} |_{\eta } (f^{m}) 
= m f(\eta)^{m-1} \cdot {\cal E}_{P} |_{\eta} f 
= m f(\eta)^{m-1} \cdot \frac{k}{m} f(\eta) 
= k \cdot f(\eta)^{m}$, 
so that $f^{m} \in {\cal O}^{k}(W)$. 
Hence, 
${\cal O}^{k/m}(W) = 
\{ f \in {\cal O}(W){;~} f^{m} \in {\cal O}^{k}(W) \}$, 
as required. 
\\ \
For $f \in {\cal O}^{k/m}(W)$, put 
$g := f \circ \breve{\eta} 
\in {\cal O}(\breve{\mathbb{C}}^{\times}_{\eta, W})$. 
For any 
$\lambda \in \breve{\mathbb{C}}^{\times}_{\eta, W}$, 
by 
a new variable $\lambda' := \lambda \lambda_{1}$, 
$k \cdot g(\lambda)^{m} 
= k \cdot f(R_{\lambda} \eta)^{m} 
= {\cal E}_{P} |_{R_{\lambda} \eta} (f^{m}) 
= \frac{d}{d \lambda_{1}} f(R_{\lambda_{1}} R_{\lambda} \eta)^{m} 
|_{\lambda_{1} = 1} 
= \frac{d \lambda'}{d \lambda_{1}} |_{\lambda_{1} = 1} \cdot 
\frac{d}{d \lambda'} (f(R_{\lambda'} \eta)^{m}) 
|_{\lambda' = \lambda} 
= \lambda \frac{d}{d \lambda} (f(R_{\lambda} \eta)^{m}) 
= \lambda \frac{d}{d \lambda} (g(\lambda)^{m})$. 
Assume that 
$W |_{p(\eta)}$ is connected for every $\eta \in W$. 
Then 
$1 \in \breve{\mathbb{C}}^{\times}_{\eta, W}$, 
which is a connected open subset of $\mathbb{C}^{\times}$ 
for every $\eta \in W$ by the claim (B) 
since $W |_{p(\eta)} = W \cap P |_{p(\eta)}$ is open in $P |_{p(\eta)}$. 
By the claim (A),  
$g(\lambda)^{m} = g(1)^{m} \cdot \lambda^{k}$, 
i.e., 
$f(R_{\lambda} \eta)^{m} = f(\eta)^{m} \cdot \lambda^{k}$. 
Conversely, let $f \in {\cal O}(W)$ satisfy that 
$f(R_{\lambda} \eta)^{m} = \lambda^{k} \cdot f(\eta)^{m}$ 
$(\eta \in W, \lambda \in \breve{\mathbb{C}}^{\times}_{\eta, W})$. 
Then 
${\cal E}_{P} |_{\eta } (f^{m}) 
= \left. \frac{d}{d \lambda} \right|_{\lambda = 1} f(R_{\lambda} \eta)^{m} 
= \left. \frac{d}{d \lambda} \right|_{\lambda = 1} 
\lambda^{k} \cdot f(\eta)^{m} = k \cdot f(\eta)^{m}$, 
so that 
$f^{m} \in {\cal O}^{k}(W)$. 
Hence, 
$f \in {\cal O}^{k/m}(W)$. 
\\ \
(b) 
Take any $\xi \in V^{f^{\sim}} = p((P |_{V})^{f})$. 
Then 
$P^{f} |_{\xi} = P^{f} \cap p^{-1}(\xi) \subseteq p^{-1}(\xi) = 
P |_{\xi}$. 
Conversely, take any $\eta' \in P |_{\xi}$. 
There exists $\eta \in P^{f}$ such that $p(\eta) = \xi$. 
By $p(\eta') = \xi = p(\eta)$, 
there exists $\lambda \in \mathbb{C}^{\times}$ 
such that $\eta' = R_{\lambda} \eta$. 
Then 
$\lambda \in \breve{\mathbb{C}}^{\times}_{\hat{\eta}, P}$. 
By the claim (a), 
$f(\eta')^{m} =f(R_{\lambda} \eta)^{m} 
= \lambda^{k} \cdot f(\eta)^{m} \ne 0$. 
Hence, $\eta' \in P^{f} |_{\xi}$. 
Then 
$P |_{\xi} \subseteq P^{f} |_{\xi}$, 
so that 
$P |_{\xi} = P^{f} |_{\xi} = p^{-1}(\xi) \cap P^{f}$ 
and 
$P |_{V^{f^{\sim}}} = p^{-1}(V) \cap P^{f} = (P |_{V})^{f}$. 
Assume that $V$ is open in $M$. 
Then 
$(P |_{V})^{f} = p^{-1}(V)^{f}$ 
is open in $P$. 
Since the canonical projection $p$ is an open map, 
$V^{f^{\sim}} = p((P |_{V})^{f})$ 
is open in $M$. 
Assume that $V$ is connected. 
Since $V$ and each fiber of $P$ is path-connected, 
$P |_{V}$ is a connected open complex submanifold of $P$. 
By Lemma 2.1 (C), 
$(P |_{V})^{f}$ is connected in $P$. 
Hence, $p((P |_{V})^{f}) = V^{f^{\sim}}$ is connected. 
\qed
\smallskip \\ \
\textsc{Lemma 2.3.} 
{\it 
Let 
$P_{1}$ and $P_{2}$ 
be holomorphic principal $\mathbb{C}^{\times}$-bundles 
on a complex manifold $M$ with the projections 
$p_{1}: P_{1} \to M$ and $p_{2}: P_{2} \to M$ 
and the right actions 
$R_{\lambda}$ 
of 
$\lambda \in \mathbb{C}^{\times}$ 
on $P_{1}$ and $P_{2}$, 
respectively. 
Let 
$\kappa: P_{1} \to P_{2}$ 
be a holomorphic immersion such that 
$p_{2} \circ \kappa = p_{1}$. 
Assume that there exists 
$\ell \in \mathbb{C}^{\times} \cap \mathbb{Z}$ 
such that 
$\kappa \circ R_{\lambda} = \lambda^{\ell} \cdot \kappa$ 
$(\lambda \in \mathbb{C}^{\times})$. 
Then one has the following results. 
}
\\ \
(A) 
{\it 
$\kappa$ is a covering map of $|\ell|$-sheets onto $P_{2}$, 
and that 
$\kappa_{*\eta} ({\cal E}_{P_{1}}) 
= \ell \cdot {\cal E}_{P_{2}} |_{\kappa(\eta)}$ 
for all $\eta \in P_{1}$. 
}
\\ \
(B) 
{\it 
$\kappa^{*}: {\cal O}^{1}(P_{2}) \to {\cal O}^{\ell}(P_{1}); 
g \mapsto \kappa^{*} g$ 
is a $\mathbb{C}$-linear isomorphism. 
}
\smallskip \\ \
{\it Proof}. 
(A) 
For any $\sigma \in \Gamma_{\rm loc}(P_{1})$, 
put 
$\kappa_{\sigma}: P_{1} |_{V_{\sigma}} \to P_{2} |_{V_{\sigma}}; 
\phi_{\sigma}(\xi, \lambda) = R_{\lambda} \sigma(\xi) \mapsto 
R_{\lambda^{\ell}} \kappa(\sigma(\xi)) 
= \phi_{\kappa \circ \sigma}(\xi, \lambda^{\ell})$, 
which is a covering map of $|\ell|$-sheets since 
$\hat{\mu}_{\ell}: \mathbb{C}^{\times} \to \mathbb{C}^{\times}; 
\lambda \to \lambda^{\ell}$ 
is a covering map of $|\ell|$-sheets. 
For another 
$\sigma' \in \Gamma_{\rm loc}(P_{1})$, 
there exists 
$\lambda_{\sigma, \sigma'}: 
V_{\sigma} \cap V_{\sigma'} \to \mathbb{C}^{\times}; 
\xi \mapsto \lambda_{\sigma, \sigma'}(\xi)$ 
such that 
$\sigma (\xi) = R_{\lambda_{\sigma, \sigma'}(\xi)} \sigma'(\xi)$, 
$\phi_{\sigma}(\xi, \lambda) 
= \phi_{\sigma'}(\xi, \lambda_{\sigma, \sigma'}(\xi) \lambda)$ 
and 
$\kappa \circ \sigma 
= \lambda_{\sigma, \sigma'}^{\ell} \cdot \kappa \circ \sigma'$. 
Then 
$\kappa_{\sigma'} 
(\phi_{\sigma}(\xi, \lambda)) 
= \kappa_{\sigma'} 
(\phi_{\sigma'}(\xi, \lambda_{\sigma, \sigma'}(\xi) \lambda)) 
= \phi_{\kappa \circ \sigma'} 
(\xi, \lambda_{\sigma, \sigma'}(\xi)^{\ell} \cdot \lambda^{\ell}) 
= \phi_{\lambda_{\sigma, \sigma'}(\xi)^{\ell} \cdot \kappa \circ \sigma'}
(\xi, \lambda^{\ell}) 
= \phi_{\kappa \circ \sigma} (\xi, \lambda^{\ell}) 
= \kappa_{\sigma} (\phi_{\sigma}(\xi, \lambda))$. 
Hence, 
$\kappa = \cup_{\sigma \in \Gamma_{\rm loc}(P_{1})} \kappa_{\sigma}$ 
in the sense of Lemma 2.0 (A). 
Any $\eta \in P_{2}$ admits 
$\sigma \in \Gamma_{\rm loc}(P_{2})$ 
such that $p_{2}(\eta) \in V_{\sigma}$, 
so that $\eta \in P_{2} |_{V_{\sigma}}$. 
Since $\kappa_{\sigma}$ is a covering map of $|\ell|$-sheets, 
there exists an open neighbourhood $W_{\eta}$ of $\eta$ 
in $P_{2} |_{V_{\sigma}}$ and disjoint open subsets $W_{j}$ 
of $P_{1} |_{V_{\sigma}}$ ($j \in \{ 1, \ldots, \ell \}$) 
such that 
$\kappa_{\sigma}^{-1}(W_{\eta}) 
= \cup_{j = 1}^{\ell} W_{j}$ 
and that the restriction 
$\kappa_{\sigma} |_{W_{j} \to W_{\eta}}$ 
is a homeomorphism. 
Especially, $\kappa$ is a surjection. 
Recall that 
$p_{1} = p_{2} \circ \kappa$, 
so that 
$p_{1}(\kappa^{-1}(W_{\eta})) 
= p_{2}(\kappa (\kappa^{-1} W_{\eta})) 
\subseteq p_{2} (W_{\eta}) \subseteq V_{\sigma}$. 
Then 
$\kappa^{-1}(W_{\eta}) 
= \kappa^{-1}(W_{\eta}) \cap P |_{V_{\sigma}} 
= \kappa_{\sigma}^{-1}(W_{\eta}) = \cup_{j = 1}^{\ell} W_{j}$ 
and that $\kappa |_{W_{j} \to W_{\eta}} 
= \kappa_{\sigma} |_{W_{j} \to W_{\eta}}$ 
as a homeomorphism. 
For any 
$g \in {\cal O}(P_{2})$ and $\eta \in P_{1}$, 
by means of a new variable 
$\lambda' := \lambda^{\ell}$, 
one has that 
$\kappa_{*} ({\cal E}_{P_{1}} |_{\eta}) g   
= \left. \frac{d}{d \lambda} 
(g \circ \kappa) (R_{\lambda} \eta) 
\right|_{\lambda = 1} 
= \left. \frac{d}{d \lambda} 
g ( R_{\lambda^{\ell}} \kappa(\eta) ) \right|_{\lambda = 1} 
= \left. 
\frac{d \lambda'}{d \lambda} \right|_{\lambda = 1} 
\cdot \left. 
\frac{d}{d \lambda'} g(R_{\lambda'} \kappa(\eta))  
\right|_{\lambda' = 1} 
= (\ell \cdot {\cal E}_{P_{2}} |_{\kappa(\eta)}) g$. 
Hence,  
$\kappa_{*} ({\cal E}_{P_{1}} |_{\eta}) 
= \ell \cdot {\cal E}_{P_{2}} |_{\kappa(\eta)}$, 
as required. 
\\ \
(B) 
For any $g \in {\cal O}^{1}(P_{2})$, 
${\cal E}_{P_{1}} (\kappa^{*} g) 
= g_{*} (\kappa_{*} {\cal E}_{P_{1}}) \circ \kappa
= g_{*} (\ell {\cal E}_{P_{2}} |_{\kappa}) 
= \ell \cdot g \circ \kappa 
= \ell \cdot \kappa^{*} g$. 
Hence, 
the $\mathbb{C}$-linear map 
$\kappa^{*}: {\cal O}^{1}(P_{2}) \to {\cal O}^{\ell}(P_{1})$ 
is well-defined. 
Take any $f \in {\cal O}^{\ell}(P_{1})$. 
Put 
$\breve{f}: P_{2} \to \mathbb{C}; 
\kappa (\eta) \mapsto 
\breve{f}(\kappa(\eta)) := f(\eta)$, 
which is well-defined as a map because of 
$f(\kappa^{-1}(\kappa(\eta))) 
= f(\{ R_{\omega} \eta{;~} \omega \in \mathbb{Z}_{|\ell|} \}) 
= \{ f(\eta) \}$ for any $\eta \in P$. 
Moreover, 
$f = \breve{f} \circ \kappa$. 
Since $\kappa$ is locally biholomorphic by (A),  
$\breve{f} \in {\cal O}(P_{2})$ by $f \in {\cal O}(P_{1})$. 
Moreover, 
$({\cal E}_{P_{2}} \breve{f}) \circ \kappa 
= \breve{f}_{*} ({\cal E}_{P_{2}}) \circ \kappa 
= \breve{f}_{*} (\kappa_{*} {\cal E}_{P_{1}}/\ell) \circ \kappa 
= f_{*} ({\cal E}_{P_{1}} /\ell) 
= \ell f /\ell 
= \breve{f} \circ \kappa$. 
Hence, $\breve{f} \in {\cal O}^{1}(P_{2})$ 
such that $\kappa^{*} \breve{f} = f$. 
Then $\kappa^{*}$ is surjective. 
If $\kappa^{*} g = 0$, then 
$\{ g(\eta_{2}); \eta_{2} \in P_{2} \} 
= \{ g(\kappa(\eta_{1}); \eta_{1} \in P_{1} \} 
= \{ (\kappa^{*} g)(\eta_{1}); \eta_{1} \in P_{1} \} 
= \{ 0 \}$ since $\kappa$ is surjective, 
so that $g = 0 \in {\cal O}(P_{1})$. 
Then $\kappa^{*}$ is injective. 
Hence, $\kappa^{*}$ is a $\mathbb{C}$-linear isomorphism.~\qed 
\smallskip \\ \
{\bf 3. Immersivity of holomorphic line bundles}. 
Let $L$ be a holomorphic line bundle on a complex manifold $M$. 
Fix $m \in \mathbb{N}$. 
Put 
\begin{eqnarray*}
& &\mu'_{L}: 
Q_{L} \to Q_{L^{\otimes m}}; \eta' \mapsto (\eta')^{\otimes m}; 
\\
& &\mu_{L}: P_{L} \to 
P_{L^{\otimes m}}; 
\eta \mapsto \eta^{\otimes m}; 
<\!\! \eta^{\otimes m}, \sum_{i} \mu'_{L}(\eta_{i}') \!\!>_{L^{\otimes m}} 
:= \sum_{i} <\!\! \eta, \eta_{i}' \!\!>_{L}^{m}. 
\end{eqnarray*}
\
\textsc{Lemma 3.0.} 
(A) 
{\it 
Put 
$\mu := \mu'_{L}$ (resp. $\mu_{L}$), 
$P := Q_{L}$ (resp. $P_{L}$) 
and $\breve{P} := Q_{L^{\otimes m}}$ (resp. $P_{L^{\otimes m}}$) 
with the canonical projections $p := q_{L}$ (resp. $p_{L}$) 
and 
$\breve{p} := q_{L^{\otimes m}}$ (resp. $p_{L^{\otimes m}}$). 
Then 
$\mu$ is a covering map with $m$-sheets such that 
$\breve{p} \circ \mu = p$, 
$\mu \circ R_{\lambda} = \lambda^{m} \cdot \mu$ 
$(\lambda \in \mathbb{C}^{\times})$ 
and that 
$\mu_{*\eta} ({\cal E}_{P}) = m \cdot {\cal E}_{\breve{P}} |_{\mu(\eta)}$ 
for all $\eta \in P$.
For $\eta_{1}, \eta_{2} \in P$, 
$\mu(\eta_{1}) = \mu(\eta_{2})$ 
if and only if there exists $\omega \in \mathbb{Z}_{m}$ 
such that $\eta_{1} = R_{\omega} \eta_{2}$. 
Moreover, 
$\widetilde{\mu}^{*}: {\cal O}^{1}(\breve{P}) \to {\cal O}^{m}(P); 
g \mapsto \widetilde{\mu}^{*}(g) := \mu^{*} g$ 
is a $\mathbb{C}$-linear isomorphism. 
} 
\\ \ 
(B) 
{\it 
For 
$f \in {\cal O}^{m}(P_{L})$, 
put $\breve{f} := (\widetilde{\mu}_{L}^{*})^{-1} (f) 
\in {\cal O}^{1}(P_{L^{\otimes m}})$ 
so that $f = \widetilde{\mu}_{L}^* \breve{f}$. 
If 
$M^{f^{\sim}} \ne \emptyset$, 
then there exists unique 
$1_{\breve{f}} \in \Gamma_{M^{f^{\sim}}}^{} 
(P_{L^{\otimes m}}^{\breve{f}} )$ 
such that  
$(\breve{f} \circ 1_{\breve{f}})(M^{f^{\sim}}) 
= \{ 1 \}$ in $\mathbb{C}^{\times}$. 
} 
\smallskip \\ \
{\it Proof}. 
(A) 
Note that $\mu$ is well-defined as a map. 
For the covering map 
$\hat{\mu}_{m}: \mathbb{C}^{\times} \to \mathbb{C}^{\times}; 
\lambda \mapsto \lambda^{m}$ 
and any $\sigma \in \Gamma^{\rm cont}_{\rm loc} (P)$, 
one has that 
$\mu (\phi_{\sigma}(\xi, \lambda)) 
= \phi_{\sigma}(\xi, \lambda)^{\otimes m} 
= \phi_{\sigma^{\otimes m}}(\xi, \hat{\mu}_{m}(\lambda))$ 
($(\xi, \lambda) \in V_{\sigma} \times \mathbb{C}^{\times}$), 
so that 
$\mu$ is a holomorphic immersion such that 
$\mu \circ R_{\lambda} = \lambda^{m} \cdot \mu$ 
($\lambda \in \mathbb{C}^{\times}$). 
Then the result follows from Lemma 2.3 (A, B) 
for $\kappa := \mu$. 
\\ \
(B) 
By Lemma 2.2 (C) (b), 
$M^{f^{\sim}} = p (P^{f}) 
= \breve{p} (\breve{P}^{\breve{f}})$. 
For 
$s \in \Gamma_{\rm loc}^{} (\breve{P}^{\breve{f}})$, 
put 
$1_{s}: V_{s} \to \breve{P}^{\breve{f}}; 
\xi' \mapsto 
R_{\breve{f}( s(\xi') )^{-1}} s(\xi')$. 
Then 
$1_{s} \in \Gamma_{V_{s}}(\breve{P}^{\breve{f}})$ 
such that 
$\breve{f} \circ 1_{s} \equiv 1$ on $V_{s}$. 
Put 
$1_{\breve{f}} 
:= \cup_{s \in \Gamma_{\rm loc}^{} (\breve{P}^{\breve{f}})} 1_{s}$. 
By Lemma 2.0 (A), 
$1_{\breve{f}} \in \Gamma_{M^{f^{\sim}}}(\breve{P})$ 
such that $\breve{f} \circ 1_{\breve{f}} \equiv 1$ 
on $M^{f^{\sim}}$. 
For another $t \in \Gamma_{M^{f^{\sim}}}(\breve{P})$ 
such as $\breve{f} \circ t = 1$ on $M^{f^{\sim}}$, 
there exist $\lambda(\xi) \in \mathbb{C}^{\times}$ 
such that $t(\xi) = R_{\lambda(\xi)} 1_{\breve{f}} (\xi)$ 
$(\xi \in M^{f^{\sim}})$. 
Because of 
$\breve{f} \in {\cal O}^{1}(\breve{P})$, 
$1 = \breve{f} (t(\xi)) 
= \breve{f}(R_{\lambda(\xi)} 1_{\breve{f}} (\xi)) 
= \lambda(\xi)~(\xi \in M^{f^{\sim}})$. 
Then $t = 1_{\breve{f}}$, as required.~\qed 
\smallskip \\ \
\textsc{Lemma 3.1.} 
(A) 
{\it 
Put 
$\widetilde{\iota}_{L^{\otimes m}}: 
\Gamma_{M}(L^{\otimes m}) \to {\cal O}^{m}(P_{L}); 
\breve{s} \mapsto \iota_{L^{\otimes m}}(\breve{s}) \circ \mu_{L}$. 
Then 
$\widetilde{\iota}_{L^{\otimes m}}$ 
is a $\mathbb{C}$-linear isomorphism. 
For any $f \in {\cal O}^{m}(P_{L})$, 
put 
$\breve{s}_{f} := \widetilde{\iota}_{L^{\otimes m}}^{-1}(f) 
\in \Gamma_{M}(L^{\otimes m})$. 
Then 
$<\! 1_{\breve{f}}, \breve{s}_{f}\!>_{L^{\otimes m}} = 1$ on $M^{f^{\sim}}$. 
} 
\\ \ 
(B) 
{\it 
$L^{\otimes m}$ is immersive 
if and only if there exist $N \in \mathbb{N}$ 
and $\mathbb{C}$-linearly independent 
$f_{j} \in {\cal O}^{m}(P_{L})$ 
$(j \in \{ 0, \ldots, N \})$ 
such that 
$(f_{0}, \ldots, f_{N}): P_{L} \to \mathbb{C}_{N+1}$ 
is an immersion.  
In this case, the map 
$(f_{0}, \ldots, f_{N})$ 
is non-vanishing. 
}
\\ \
(C) 
{\it 
$L$ is very ample 
if and only if there exist $N \in \mathbb{N}$ 
and $\mathbb{C}$-linearly independent 
$f_{j} \in {\cal O}^{1}(P_{L})$ 
$(j \in \{ 0, \ldots, N \})$ 
such that 
$(f_{0}, \ldots, f_{N}): P_{L} \to \mathbb{C}_{N+1}$ 
is an injective immersion. 
}
\smallskip \\ \
{\it Proof}. 
(A) Because of 
$\widetilde{\iota}_{L^{\otimes m}} \circ R_{\lambda} 
= \lambda^{m} \cdot \widetilde{\iota}_{L^{\otimes m}}$ 
($\lambda \in \mathbb{C}^{\times}$), 
one has that 
$\widetilde{\iota}_{L^{\otimes m}}$ 
is well-defined as a $\mathbb{C}$-linear map. 
Take any $f \in {\cal O}^{m}(P_{L})$. 
Put 
$\breve{f} := (\widetilde{\mu}_{L}^{*})^{-1}(f) 
\in {\cal O}^{1}(P_{L^{\otimes m}})$ 
by Lemma 3.0 (A), 
so that $\breve{f} \circ \mu_{L} = f$. 
For 
$\sigma' \in \Gamma_{\rm loc}^{} (\breve{P})$, 
let 
$s_{\sigma'} \in 
\Gamma_{V_{\sigma'}}^{} (L^{\otimes m})$ 
be taken such as 
$<\!\! \sigma'(\xi), s_{\sigma'}(\xi) \!\!>_{L^{\otimes m}} = 1$ 
($\xi \in V_{\sigma'}$). 
Then 
$\Gamma_{M} (L^{\otimes m}) \ni \breve{s}_{f} 
:= \cup_{\sigma' \in \Gamma_{\rm loc}(\breve{P})} 
(\breve{f} \circ \sigma') \cdot s_{\sigma'}$ 
by Lemma 2.0 (A). 
For 
$(\xi, \lambda) \in V_{\sigma'} \times \mathbb{C}^{\times}$, 
$\breve{f}( \lambda \cdot \sigma'(\xi)) 
= \lambda \cdot \breve{f}(\sigma'(\xi)) 
= <\! \lambda \cdot \sigma'(\xi), 
\breve{f}(\sigma'(\xi)) \cdot 
s_{\sigma'}(\xi) \!>_{L^{\otimes m}} 
= \iota_{L^{\otimes m}} (\breve{s}_{f}) 
|_{\lambda \cdot \sigma'(\xi),}$ 
so that 
$\breve{f}= \iota_{L^{\otimes m}}(\breve{s}_{f})$ 
and 
$f = \breve{f} \circ \mu_{L} 
= \iota_{L^{\otimes m}}(\breve{s}_{f}) \circ \mu_{L}$. 
Hence, 
$\widetilde{\iota}_{L^{\otimes m}}$ 
is surjective. 
Assume that 
$\widetilde{\iota}_{L^{\otimes m}}(\breve{s}) = 0$. 
For $\eta \in P_{L}$,
$0 = <\!\! \mu_{L}(\eta), 
\breve{s}(p_{L^{\otimes m}}(\mu_{L}(\eta))) \!\!>_{L^{\otimes m}}$. 
By 
$M = p_{L^{\otimes m}}(\mu_{L}(P_{L}))$, 
$\breve{s} = 0$. 
Hence, 
$\widetilde{\iota}_{L^{\otimes m}}$ 
is injective. 
Moreover, 
$< 1_{\breve{f}}, \breve{s}_{f} >_{L^{\otimes m}} 
= \iota_{L^{\otimes m}} (\breve{s}_{f}) 
\circ 1_{\breve{f}} 
= \breve{f} \circ 1_{\breve{f}} 
= 1$ on $M^{f^{\sim}}$. 
\\ \ 
(B0) {\it Claim. 
Let $F: P_{L} \to \mathbb{C}_{N+1}$ be an immersion 
such that there exists $g \in {\cal O}(\mathbb{C}^{\times})$ 
such as $F \circ R_{\lambda} = g(\lambda) \cdot F$ 
$(\lambda \in \mathbb{C}^{\times})$. 
Then $F$ is non-vanishing:} 
In fact, {\it suppose that $F(\eta) = 0$ for some $\eta \in P_{L}$}. 
Then 
$F_{*\eta} {\cal E}_{P_{L}} 
= \left. \frac{d}{d \lambda} \right|_{\lambda = 1} F(R_{\lambda} \eta) 
= \left. \frac{d}{d \lambda} \right|_{\lambda = 1} g(\lambda) \cdot F(\eta) 
= \left. \frac{d}{d \lambda} \right|_{\lambda = 1} 0 = 0$, 
so that $F$ is not an immersion, {\it a contradiction}. 
\\ \ 
(B1) Assume that $L^{\otimes m}$ 
is immersive. 
Then there exist 
$s_{j} \in \Gamma (L^{\otimes m})~(j \in \{ 0, \ldots, N \})$ 
such that the projectivization 
$\varphi: M \to P_{N} \mathbb{C}; 
p_{L^{\otimes m}}(\eta) \to \pi_{N}(\Phi(\eta))$ 
of $\Phi := (\iota_{L}(s_{0}), \ldots, \iota_{L}(s_{N})): 
P_{L^{\otimes m}} \to \mathbb{C}_{N+1}$ 
is an immersion such that 
$\varphi \circ p_{L^{\otimes m}} = \pi_{N} \circ \Phi$. 
\\ \
{\it Claim. $\Phi$ is an immersion:} 
In fact, for $X \in T' P_{L^{\otimes m}}$, 
assume that $\Phi_{*'} X = 0$. 
Then 
$\varphi_{*'}(p_{L^{\otimes m}*'} X) 
= \pi_{N*'} (\Phi_{*'} X) = 0$. 
Since $\varphi$ is an immersion, 
$p_{L^{\otimes m}*'} X = 0$. 
Hence, 
$X = \lambda \cdot {\cal E}_{P_{L}}$ 
for some $\lambda \in \mathbb{C}$. 
Then  
$0 = \Phi_{*'} X 
= \lambda \cdot {\cal E}_{\mathbb{C}_{N+1} \backslash \{ 0 \}}$, 
so that $\lambda = 0$ and $X = 0 \cdot {\cal E}_{P_{L}} = 0$. 
Hence, $\Phi$ is an immersion, as required. 
\\ \
By the step (B0), 
one has that 
$\Phi (P_{L}) \subseteq \mathbb{C}_{N+1} \backslash \{ 0 \}$. 
Note that 
$\Phi \circ R_{\lambda} = \lambda \cdot \Phi$ 
($\lambda \in \mathbb{C}^{\times}$), 
and that 
$\mathbb{C}_{N+1} \backslash \{ 0 \}$ 
is a holomorphic principal $\mathbb{C}^{\times}$-bundle 
on $P_{N} \mathbb{C}$. 
By Lemma 2.3 (A), 
$\Phi_{*'\eta} {\cal E}_{P_{L}} 
= {\cal E}_{\mathbb{C}_{N+1}} |_{\Phi(\eta)}$ 
for $\eta \in P_{L}$. 
Put 
$(f_{0}, \ldots, f_{N}) 
:= \Phi \circ \mu_{L}: P_{L} \to \mathbb{C}_{N+1}\backslash \{ 0 \}$, 
which is an immersion such that 
$\Phi \circ \mu_{L} \circ R_{\lambda} 
= \lambda^{m} \cdot \Phi \circ \mu_{L}$ 
($\lambda \in \mathbb{C}^{\times}$) 
because of Lemma 3.0 (A), so that 
$f_{i} \in {\cal O}^{m}(P_{L})$ 
($i \in \{ 0, \ldots, N \}$), 
as required. 
\\ \
(B2) 
Conversely, 
assume that there exist 
$N \in \mathbb{N}$ 
and 
$f_{j} \in {\cal O}^{m}(P_{L}) (j \in \{ 0, \ldots, N \})$ 
such that 
$F := (f_{0}, \ldots, f_{N}): P_{L} \to \mathbb{C}_{N+1}$ 
is an immersion. 
Then  
$F \circ R_{\lambda} = \lambda^{m} \cdot F$ 
($\lambda \in \mathbb{C}^{\times}$). 
By Lemma 5.0 (A), 
$F_{*} {\cal E}_{P_{L}} 
= m \cdot {\cal E}_{\mathbb{C}_{N+1} \backslash \{ 0 \}}$. 
By the step (B0), 
$F (P_{L}) \subseteq \mathbb{C}_{N+1} \backslash \{ 0 \}$. 
Then the projectivization of $F$ 
is well-defined as a holomorphic map 
$\varphi: M \to P_{N} \mathbb{C}; 
p_{L}(\eta) \to \pi_{N}(F(\eta)))$ 
such that 
$\varphi \circ p_{L} = \pi_{N} \circ F$. 
\\ \
{\it Claim. $\varphi$ is an immersion:} 
In fact, suppose that 
$\varphi_{*'} (p_{L*'} X) = 0$ and 
$p_{L*'} X \ne 0$ for some $(X, \eta) \in T' P_{L} |_{\eta} \times P_{L}$. 
Then 
$X \ne 0$ and 
$\pi_{N*'} (F_{*'} X) = \varphi_{*'} (p_{L*'} X) = 0$. 
There exists 
$\lambda \in \mathbb{C}$ 
such that 
$F_{*'} X 
= \lambda \cdot {\cal E}_{\mathbb{C}_{N+1} \backslash \{ 0 \}} |_{F(\eta)} 
= F_{*'} (\lambda \cdot {\cal E}_{P_{L}}/m)$, 
so that 
$X = \lambda \cdot {\cal E}_{P_{L}}/m$, 
since $F$ is an immersion. 
Hence, 
$p_{L*'} X = \lambda \cdot p_{L*'}({\cal E}_{P_{L}} |_{\eta}) = 0$, 
a contradiction. 
\\ \
(BC3) 
In the step (B2), there exists a subset 
$\{j_{0}, \ldots, j_{N'} \}$ of $\{ 0, \ldots, N \}$ 
such that 
$\{ f_{j_{0}}, \ldots, f_{j_{N'}} \}$ 
is a $\mathbb{C}$-linear base of 
${\cal F} := \sum_{j = 0}^{N} \mathbb{C} f_{j}$. 
Put 
$\Phi' := (f_{j_{0}}, \ldots, f_{j_{N'}}): 
P_{L} \to \mathbb{C}_{N' + 1}$. 
Then $\Phi'$ is non-vanishing or an immersion 
(resp. an injection) 
if and only if $\Phi \circ \mu_{L} = (f_{0}, \ldots, f_{N})$ 
is non-vanishing or an immersion (resp. an injection), 
respectively. 
\\ \
(C1) 
Assume that $L$ is very ample. 
Then there exist 
$s_{j} \in \Gamma (L)~(j \in \{ 0, \ldots, N \})$ 
such that the projectivization 
$\varphi: M \to P_{N} \mathbb{C}; 
p_{L}(\eta) \to \pi_{N}(\Phi(\eta))$ 
of $\Phi := (\iota_{L}(s_{0}), \ldots, \iota_{L}(s_{N})): 
P_{L} \to \mathbb{C}_{N+1}$ 
is an injective immersion. 
By the step (B1) for $m = 1$, 
$\Phi$ is an immersion. 
\\ \
{\it Claim. $\Phi$ is injective:} 
For $\eta_{1}, \eta_{2} \in P_{L}$, 
assume that $\Phi(\eta_{1}) = \Phi(\eta_{2})$. 
Then 
$\varphi (p_{L} (\eta_{1})) 
= \pi_{N} (\Phi (\eta_{1})) 
= \pi_{N} (\Phi (\eta_{2})) 
= \varphi (p_{L} (\eta_{2}))$. 
Since $\varphi$ is injective, 
$p_{L}(\eta_{1}) = p_{L}(\eta_{2})$. 
Hence, 
$\eta_{2} = R_{\lambda} \eta_{1}$ 
for some $\lambda \in \mathbb{C}^{\times}$. 
By $\Phi \circ R_{\lambda} = \lambda \cdot \Phi$ 
($\lambda \in \mathbb{C}^{\times}$), 
$\Phi (\eta_{2}) = \lambda \cdot \Phi (\eta_{1}) 
= \lambda \cdot \Phi (\eta_{2})$. 
By the step (B1), $\Phi(\eta_{2}) \ne 0$. 
Hence, $\lambda = 1$, 
so that $\eta_{2} = \eta_{1}$. 
\\ \
(C2) 
Conversely, assume that there exist 
$N \in \mathbb{N}$ 
and 
$f_{j} \in {\cal O}^{1}(P_{L}) (j \in \{ 0, \ldots, N \})$ 
such that 
$F := (f_{0}, \ldots, f_{N}): P_{L} \to \mathbb{C}_{N+1}$ 
is an injective immersion. 
By the step (B2), the projectivization $\varphi$ of $F$ 
is well-defined as an immersion such that 
$\Pi_{N} \circ F = \varphi \circ p_{L}$. 
\\ \
{\it Claim. $\varphi$ is injective:} 
In fact, 
assume that $\varphi(p_{L}(\eta_{1})) = \varphi(p_{L}(\eta_{2}))$. 
Then $\Pi_{N}(F(\eta_{1})) = \Pi_{N}(F(\eta_{2}))$. 
There exists $\lambda \in \mathbb{C}^{\times}$ 
such that $F(\eta_{1}) = \lambda \cdot F(\eta_{2}) 
= F(R_{\lambda} \eta_{2})$, 
so that $\eta_{1} = R_{\lambda} \eta_{2}$ 
since $F$ is injective. 
Hence, $p_{L}(\eta_{1}) = p_{L}(R_{\lambda} \eta_{2}) 
= p_{L}(\eta_{2})$, as required.~\qed 
\smallskip \\ \
\textsc{Lemma 3.2.} 
{\it 
Let $L$ be a holomorphic line bundle on $M$. 
Assume that $L$ is immersive (resp. very ample). 
Then $L^{\otimes m}$ is immersive (resp. very ample) 
for any $m \in \mathbb{N}$. 
} 
\smallskip \\ \
{\it Proof.} 
(0) 
Assume that $L$ is immersive. 
By Lemma 3.1 (B), 
there exist 
$N \in \mathbb{N}$ 
and $f_{i} \in {\cal O}^{1}(P_{L})$ 
($i \in \{ 0, \ldots, N \}$) 
such that 
$F := (f_{0}, \ldots, f_{N}): P_{L} \to \mathbb{C}_{N_{1} +1}$ 
is an immersion, which is non-vanishing by the step (B0) of 
the proof of Lemma 3.1 (B). 
Let 
$F^{\otimes m}:= (f_{J})_{J \in \{ 0, \ldots, N \}^{\times m}}: 
P_{L} \to \mathbb{C}_{(N_{1} + 1)^{m}}$ 
be defined as 
$f_{J} := f_{j_{1}} \cdots f_{j_{m}} \in {\cal O}^{m}(P_{L})$ 
for any $J := (j_{1}, \ldots, j_{m}) \in \{ 0, \ldots, N \}^{\times m}$. 
Take any $(\eta, X) \in P_{L} \times T_{\eta}' P_{L}$ 
such as $F^{\otimes m}_{*\eta} X = 0$. 
Since $F$ is non-vanishing, there exists 
$j_{\eta} \in \{ 0, \ldots, N \}$ 
such that $f_{j_{\eta}}(\eta) \ne 0$. 
Then 
$0 = (f_{j_{\eta}})^{m}_{*} X 
= m f_{j_{\eta}}(\eta)^{m-1} \cdot f_{j_{\eta}*} X$, 
so that 
$f_{j_{\eta}*} X = 0$. 
When $m = 1$: $X = 0$ by assumption. 
When $m \geq 2$: 
For any $k \in \{ 0, \ldots, N \}$, 
$0 = (f_{j_{\eta}}^{m-1} f_{k})_{*\eta} X 
= f_{j_{\eta}}(\eta)^{m-1} \cdot f_{k*\eta} X 
+ (m-1) f_{j_{\eta}}(\eta)^{m-2} f_{k}(\eta) \cdot f_{j_{\eta}*\eta} X 
= f_{j_{\eta}}(\eta)^{m-1} \cdot f_{k*\eta} X$, 
so that $f_{k*\eta} X = 0$. 
Then 
$F_{*\eta} X = (f_{0}, \ldots, f_{N})_{*\eta} X = 0$. 
Since $F$ is an immersion, 
$X = 0$. 
Hence, 
$F^{\otimes m}$ is an immersion. 
By the converse part of Lemma 3.1 (B), 
$L^{\otimes m}$ is immersive. 
\\ \
(1) 
Assume that there exist 
$N \in \mathbb{N}$ and $f_{i} \in {\cal O}^{1}(P_{L})$ 
($i \in \{ 0, \ldots, N \}$) 
such that 
$F := (f_{0}, \ldots, f_{N}): P_{L} \to \mathbb{C}_{N_{1} +1}$ 
is an injection. 
\\ \
{\it Claim 1. 
$F$ is non-vanishing:} 
In fact, if $F(\eta') = 0$ for some $\eta' \in P_{L}$, 
then $F(R_{\lambda} \eta') = \lambda \cdot \Phi(\eta') = 0$ 
for all $\lambda \in \mathbb{C}^{\times}$, 
which contradicts with that $F$ is an injection, 
as required. 
\\ \
{\it Claim 2. 
For each $J \in \{0, \ldots, N \}^{m}$, 
there exists $\breve{f}_{J} := (\mu_{L^{*}})^{-1} (f_{J}) 
\in {\cal O}^{1}(P_{L^{\otimes m}})$ 
such that $f_{J} = \breve{f}_{J} \circ \mu_{L}$:}
It follows from Lemma 3.0 
because of $f_{J} = f_{j_{1}} \cdots f_{j_{m}} \in {\cal O}^{m}(P_{L})$ 
for each $J = (j_{1}, \ldots, j_{m}) \in \{0, \ldots, N \}^{m}$, 
as required. 
\\ \
Put 
$\breve{F}^{\otimes m} := (\breve{f}_{J})_{J \in \{0, \ldots, N \}^{m}}: 
P_{L^{\otimes m}} \to \mathbb{C}_{(N+1)^{m}}$. 
Then $F^{\otimes m} = \breve{F}^{\otimes m} \circ \mu_{L}$. 
Let 
$\eta, \eta' \in P_{L}$ 
be such that 
$\breve{F}^{\otimes m}(\mu_{L}(\eta)) 
= \breve{F}^{\otimes m}(\mu_{L}(\eta'))$. 
Then 
$F^{\otimes m}(\eta) = F^{\otimes m}(\eta')$. 
Because of $F(\eta') \ne 0$, 
there exists $j \in \{ 0, \ldots, N \}$ 
such that $f_{j}(\eta') \ne 0$. 
By  
$f_{j}(\eta)^{m} = f_{j}(\eta')^{m} \ne 0$, 
there exists 
$\omega \in \mathbb{Z}_{m}$ 
such that 
$f_{j}(\eta) = \omega \cdot f_{j}(\eta')$. 
For any $i \in \{ 0, \ldots, N \}$, 
$f_{i}(\eta) f_{j}(\eta)^{m-1} 
= f_{i}(\eta') f_{j}(\eta')^{m-1}$, 
so that 
$f_{i}(R_{\omega^{m-1}} \eta) 
= f_{i}(\eta) \omega^{m-1} = f_{i}(\eta')$, 
i.e., 
$F(R_{\omega^{m-1}} \eta) = F(\eta')$, 
so that 
$R_{\omega^{m-1}} \eta = \eta'$ 
since $F$ is an injection. 
Then 
$\mu_{L}(\eta') = \mu_{L} (R_{\omega^{m-1}} \eta) 
= \mu_{L} (\eta)$. 
Hence, 
$\breve{F}^{\otimes m}$ 
is injective. 
\\ \
(2) By the steps (0) and (1), 
$F^{\otimes m}$ 
is an injective immersion if $F$ is an injective immersion. 
Hence, $L^{\otimes m}$ is very ample if $L$ is very ample.
~\qed 
\smallskip \\ \
Let $L$ be a holomorphic line bundle on a complex manifold $M$, 
and ${\rm Her}(L\otimes_{\mathbb{R}}L, \mathbb{C})$ 
be the smooth fiber bundle of every positive definite hermitian form on 
complex line $L|_{\xi} := q_{L}^{-1}(\xi) (\xi \in M)$, 
a smooth global section of which is said to be 
{\it a hermitian fiber metric on $L$}. 
By definition, 
$L$ is said to be {\it strictly positive} 
if and only if there exists a hermitian fiber metric 
$h_{L}$ on $L$ 
such that each $\xi \in M$ 
admits some 
$e \in \Gamma_{\rm loc}(Q_{L})$ 
such that 
$h^{e}_{L} := - \log \circ h_{L} \circ e^{\times 2}: 
V_{e} \to \mathbb{R}; \xi' \mapsto := - \log h_{L}(e(\xi'), e(\xi'))$ 
is {\it strictly plurisubharmonic}, 
that is, 
for each $\xi' \in V_{e}$ and any open neighbourhood $U$ of $\xi'$ in $V_{e}$ 
with holomorphic complex coordinates 
$(z_{1}, \ldots, z_{{\rm n}_{M}}): 
U \to \mathbb{C}_{{\rm n}_{M}}$, 
{\it the Levi form} of $h^{e}_{L} |_{U}$ 
defined as 
\begin{eqnarray*}
& &{\rm Lev}(h^{e}_{L}): 
T' U |_{\xi''} \to \mathbb{R}; 
\\
& &
X := \sum_{i = 1}^{{\rm n}_{M}} 
\left. X_{i} 
\frac{\partial}{\partial z_{i}} 
\right|_{\xi''} \mapsto {\rm Lev}(h^{e}_{L}) X 
:= \sum_{i = 1}^{\rm n_{M}} 
\left. X_{i} \bar{X}_{j} 
\frac{\partial^{2} h^{e}_{L}}
{\partial z_{i} \partial \bar{z}_{j}} 
\right|_{\xi''} 
\end{eqnarray*} 
satisfies that 
${\rm Lev}(h^{e}_{L}) (T' U |_{\xi''} \backslash \{ 0 \}) 
\subseteq \mathbb{R}_{+} := \{ r \in \mathbb{R}; r > 0 \}$ 
($\xi'' \in U$) 
(cf. S.~Nakano \cite[(17.3)]{Ns1981}, 
S.~Kobayashi \cite[p.200]{Ks2005}; 
\cite[II.2.9, VI.6.1, VI.6.2]{FG2002}). 
\\ \
\textsc{Proposition 3.3}. 
(A) 
{\it 
For $\ell \in \{ -1, 1 \}$, 
let 
$Q_{L} \times_{\rho_{\ell}} \mathbb{C}$ 
be the fiber bundle associated with the principal 
$\mathbb{C}^{\times}$-bundle $Q_{L}$ 
with respect to the representation 
$\rho_{\ell}: \mathbb{C}^{\times} \to \mathbb{C}^{\times}; 
\lambda \mapsto \lambda^{\ell}$ 
of $\mathbb{C}^{\times} \cong GL_{1} \mathbb{C}$, 
i.e., 
$Q_{L} \times_{\rho_{\ell}} \mathbb{C}$ 
is defined as the quotient space of 
the product space $Q_{L} \times \mathbb{C}$ 
with the equivalent relation $\sim_{\ell}$ 
for any $(\eta, w), (\eta', w') \in Q_{L} \times \mathbb{C}$ 
such that 
$(\eta, w) \sim_{\ell} (\eta', w')$ 
if and only if 
$(\eta, w) = (R_{\lambda}^{-1} \eta, \rho_{\ell}(\lambda) w) 
:= (\lambda^{-1} \cdot \eta, \lambda^{\ell} \cdot w)$ 
for some $\lambda \in \mathbb{C}^{\times}$ 
with the equivalence class $[\eta, w]_{\ell}$ 
and the canonical projection 
$\pi_{L}: Q_{L} \times_{\rho_{\ell}} \mathbb{C} \to M; 
[\eta, w]_{\ell} \mapsto q_{L}(\eta)$. 
For $e \in \Gamma_{\rm loc}(Q_{L})$, put 
$\phi_{e, \ell}: \pi^{-1}(V_{e}) \to V_{e} \times \mathbb{C}; 
[\eta, w]_{\ell} \mapsto (\pi(\eta), \rho_{\ell}(\lambda) w)$ 
with unique 
$\lambda \in \mathbb{C}^{\times}$ 
such as $R_{\lambda} e(\pi(\eta)) = \eta$. 
Then $Q_{L} \times_{\rho_{\ell}} \mathbb{C}$ 
is a holomorphic line bundle on $M$ such that 
$\phi_{e, \ell}$ is a local trivialization for evey 
$e \in \Gamma_{\rm loc}(Q_{L})$, 
so that 
$Q_{L} \times_{\rho_{1}} \mathbb{C} \cong L$ 
and 
$Q_{L} \times_{\rho_{-1}} \mathbb{C} \cong L^{*}$ 
as holomorphic line bundles on $M$, 
and that 
$Q_{L} \times_{\rho_{1}} \mathbb{C}^{\times} \cong Q_{L}$ 
and 
$Q_{L} \times_{\rho_{-1}} \mathbb{C}^{\times} \cong Q_{L^{*}}$ 
as holomorphic principal $\mathbb{C}^{\times}$-bundles. 
}
\\ \
(B) 
{\it 
For $j \in \{ 1, 2 \}$, 
let $L_{j}$ be a holomorphic line bundle on a complex manifold $M_{j}$ 
with the canonical projection $q_{L_{j}}: L_{j} \to M_{j}$. 
Let $\varphi: M_{1} \to M_{2}$ and 
$\Phi: Q_{L_{1}} \to Q_{L_{2}}$ be holomorphic immersions such that 
$q_{L_{2}} \circ \Phi = \varphi \circ q_{L_{1}} |_{Q_{L_{1}}}$ 
and that $\Phi (Q_{L_{1}} |_{\xi}) = Q_{L_{2}} |_{\varphi(\xi)}$ 
for every $\xi \in M_{1}$. 
Then each $\xi_{1} \in M_{1}$ admits 
an open neighbourhood $U'$ of $\xi_{1}$ in $M_{1}$ 
such that 
$\varphi(U')$ is a complex submanifold of $M_{2}$ 
and $\varphi |_{U' \to \varphi(U')}$ 
is biholomorphic and that 
$\Phi |_{Q_{L_{1}} |_{U'} \to Q_{L_{2}} |_{\varphi(U')}}$ 
is surjective and locally biholomorphic. 
Assume that 
$\Phi \circ R_{\lambda} = R_{\lambda} \circ \Phi$ 
$(\lambda \in \mathbb{C}^{\times})$ 
and that 
$L_{2}$ is strictly positive. 
Then 
$L_{1}$ is strictly positive. 
If moreover $M$ is compact, 
then $L_{1}$ is positive in the sense of K.~Kodaira. 
}
\\ \
(C) 
{\it 
Let $L$ be an immersive holomorphic line bundle on a complex manifold $M$. 
Then $L$ is strictly positive. 
If moreover $M$ is compact, 
then $L$ is positive in the sense of K.~Kodaira. 
}
\\ \
(D) 
{\it 
Let $M$ be a compact complex manifold such that $K_{M}^{*}$ is 
immersive. 
Then $M$ is Fano, 
so that each connected component of $M$ is simply connected. 
}
\smallskip \\ \
{\it Proof}. 
(A) 
Put 
$f_{1}: Q_{L} \times_{\rho_{1}} \mathbb{C} \to L; 
[\eta, w]_{1} \mapsto w \cdot \eta$, 
which is well-defined as an isomorphism as well as \cite[Prop.5.4]{KN1963}, 
so that $Q_{L} \times_{\rho_{1}} \mathbb{C} \cong L$. 
Note that 
$f_{1}(Q_{L} \times_{\rho_{1}} \mathbb{C}^{\times}) 
= L \backslash 0 = Q_{L}$, 
so that 
$Q_{L} \times_{\rho_{1}} \mathbb{C}^{\times} \cong Q_{L}$. 
Put 
$f_{-1}: Q_{L} \times_{\rho_{-1}} \mathbb{C} \to L^{*}$; 
$<\!\! f_{-1}([\eta, w]_{-1}), [\eta, w']_{1} \!\!>_{L} 
:= w \cdot w'$, 
which is well-defined as an isomorphism, 
so that $Q_{L} \times_{\rho_{-1}} \mathbb{C} \cong L^{*}$. 
Note that 
$f_{-1}(Q_{L} \times_{\rho_{-1}} \mathbb{C}^{\times}) 
= L^{*} \backslash 0 = Q_{L^{*}}$, 
so that 
$Q_{L} \times_{\rho_{-1}} \mathbb{C}^{\times} 
\cong Q_{L^{*}}$. 
\\ \
(B) 
Take any $\xi_{1} \in M_{1}$. 
Since $\varphi$ is an immersion, 
there exists an open neighbourhood $U$ of $\xi_{1}$ in $M_{1}$ 
such that the restriction 
$\varphi |_{U \to \varphi(U)}$ 
is a biholomorphism (cf. \cite[Thm.I.15.5]{Hs1978}). 
Similarly, $\Phi$ is locally biholomorphic. 
There exists $\sigma \in \Gamma_{\rm loc}(Q_{L_{1}})$ 
such that $V_{\sigma} \ni \xi_{1}$. 
Put $U' := U \cap V_{\sigma} \ni \xi_{1}$. 
By $p_{L_{2}} \circ \Phi = \varphi \circ p_{L_{1}}$, 
$\Phi (\phi_{\sigma}(\xi, \lambda)) 
= \phi_{\Phi \circ \sigma}(\xi, g(\xi, \lambda))$ 
for some $g: U_{1} \times \mathbb{C}^{\times} \to \mathbb{C}^{\times}$. 
Moreover, 
$\Phi(Q_{L_{1}} |_{U'}) 
= \cup_{\xi \in U'} \Phi(Q_{L_{1}} |_{\xi}) 
= \cup_{\xi \in U'} Q_{L_{2}} |_{\varphi(\xi)} 
= Q_{L_{2}} |_{\varphi(U')}$, 
as required. 
If 
$(\eta', w') = (R_{\lambda}^{-1} \eta, \rho_{\ell}(\lambda) w) 
:= (\lambda^{-1} \cdot \eta, \lambda^{\ell} \cdot w) 
\in Q_{L_{1}} \times \mathbb{C}$, 
then 
$[\Phi(\eta'), w'] = [\Phi(\eta), w]$, 
so that 
$\breve{\Phi}: 
Q_{L_{1}} \times_{\rho_{1}} \mathbb{C} \to 
Q_{L_{2}} \times_{\rho_{1}} \mathbb{C}; 
[\eta, w] \mapsto [\Phi(\eta), w]$ 
is well-defined as a holomorphic map such that 
$(\pi_{L_{2}} \circ \breve{\Phi})[\eta, w]_{1} 
= q_{L_{2}} (\Phi(\eta)) 
= q_{L_{1}} (\eta) 
= \pi_{L_{1}}[\eta, w]_{1}$ 
and that the restriction to the fiber 
$\breve{\Phi} |_{\pi^{-1}(\xi)}: 
Q_{L_{1}} \times_{\rho_{1}} \mathbb{C} |_{\xi} 
\to Q_{L_{2}} \times_{\rho_{1}} \mathbb{C} |_{\xi}$ 
is a $\mathbb{C}$-linear isomorphism at every $\xi \in M$, 
so that 
$\breve{\Phi} |_{\pi^{-1}(\xi)}(
Q_{L_{1}} \times_{\rho_{1}} \mathbb{C}^{\times}) 
= Q_{L_{2}} \times_{\rho_{1}} \mathbb{C}^{\times}$. 
There exists a hermitian fiber metric $h_{L_{2}}$ on $L_{2}$ 
such that each $\xi \in M_{2}$ admits 
$e \in \Gamma_{\rm loc}(Q_{L_{2}})$ 
such that $V_{e} \ni \xi$ 
and that $h_{L_{2}}^{e}$ 
is strictly plurisubharmonic. 
Put 
$h_{L_{1}} := h_{L_{2}} 
\circ (f_{1} \circ \breve{\Phi} \circ f_{1}^{-1})^{\times 2}$. 
{\it Claim. 
$h_{L_{1}}$ 
is a hermitian fiber metric on $L_{1}$:} 
{\it In fact}, 
since $h_{L_{2}}$ is hermitian, so is $h_{L_{1}}$. 
Moreover, for $\eta \in L_{1} \backslash 0$, 
$\eta' := (f_{1} \circ \breve{\Phi} \circ f_{1}^{-1}) \eta 
\in L_{2} \backslash 0$ 
and 
$h_{L_{1}}(\eta, \eta) = h_{L_{2}}(\eta', \eta') > 0$, 
i.e., 
$h_{L_{1}}$ is positive definite, {\it as required}. 
Take any $\xi_{1} \in M_{1}$. 
Put 
$\xi_{2} := \varphi(\xi_{1}) \in M_{2}$. 
Take 
$e_{2} \in \Gamma_{\rm loc}(Q_{L_{2}})$ 
such that 
$\xi_{2} \in V_{e_{2}}$ 
and that 
$h_{L_{2}}^{e_{2}}: V_{e_{2}} \to \mathbb{R}$ 
is strictly plurisubharmonic. 
Since $\varphi$ is continuous, 
there exists an open neighbourhood 
$U_{1}$ of $\xi_{1}$ in $U'$ 
such that $\varphi(U_{1}) \subseteq V_{e_{2}}$ 
and that 
$U_{1}$ has holomorphic complex coordinates 
$\phi_{1}: U_{1} \to \mathbb{C}_{{\rm n}_{M_{1}}}$. 
Moreover, 
there exist an open subset $W_{1}$ of $Q_{L_{1}} |_{U_{1}}$ 
and an open neighbourhood $W_{2}$ of $e_{2}(\xi_{2})$ 
in $Q_{L_{2}} |_{\varphi(U_{1})}$ 
such that $\Phi |_{W_{1} \to W_{2}}$ is biholomorphic. 
Put 
$U_{11} := q_{L_{1}}(W_{1})$ 
and 
$U_{22} := q_{L_{2}}(W_{2})$. 
Then 
$(q_{L_{2}} \circ \Phi |_{W_{1} \to W_{2}}) |_{W_{1}} 
= (\varphi |_{U_{11} \to U_{22}} \circ q_{L_{1}}) |_{W_{1}}$. 
Put 
$e_{1} := 
(\Phi |_{W_{1} \to W_{2}})^{-1} \circ e_{2} 
\circ \varphi |_{U_{11}.}$ 
Then 
$\varphi |_{U_{11} \to U_{22}} \circ q_{L_{1}} \circ e_{1} 
= q_{L_{2}} \circ e_{2} \circ \varphi |_{U_{11}} 
= \varphi |_{U_{11}}$ 
so that $q_{L_{1}} \circ e_{1} = {\rm id}_{U_{11}.}$ 
Hence, 
$e_{1} \in \Gamma_{U_{11}}(Q_{L_{1}} |_{U_{11}})$. 
{\it 
Claim. 
$h_{L_{1}}^{e_{1}}: U_{1} \to \mathbb{R}$ 
is strictly prulisubharmonic:} 
{\it In fact}, 
$h_{L_{1}}^{e_{1}} |_{U_{11}} 
= - \log \circ h_{L_{1}} \circ e_{1}^{\times 2}
= - \log \circ h_{L_{2}} \circ (f_{1} \circ \breve{\Phi} \circ f_{1}^{-1} 
\circ (\Phi |_{W_{1} \to W_{2}})^{-1} \circ e_{2} 
\circ \varphi |_{U_{11}} )^{\times 2} 
= - \log \circ h_{L_{2}} \circ (e_{2} \circ \varphi |_{U_{11}})^{\times 2} 
= h_{L_{2}}^{e_{2}} \circ \varphi |_{U_{11}}: U_{11} \to \mathbb{R}$. 
For any 
$X \in T' U_{11} \backslash 0$, 
$\varphi_{*} X \in T' V \backslash 0$ 
since $\varphi$ is an immersion. 
As \cite[p.57, Exercise 4]{FG2002}, 
${\rm Lev}(h_{L_{1}}^{e_{1}}) X
= 
{\rm Lev}(h_{L_{2}}^{e_{2}} |_{\varphi(U_{1})} \circ \varphi |_{U_{11}}) X 
= {\rm Lev}(h_{L_{2}}^{e_{2}} |_{\varphi(U_{1})}) (\varphi_{*} X) 
> 0$, 
{\it as required}. 
Hence, $L_{1}$ is strictly positive. 
Assume that $M_{1}$ is compact. 
By virtue of S.~Kobayashi \cite[p.4, line 4]{Ks2005}, 
$L_{1}$ is positive in the sense of K.~Kodaira \cite{Kk1954}. 
\\ \
(C) 
Since $L$ is immersive, 
there exist 
$N \in \mathbb{N}$ and 
$s_{i} \in \Gamma(L) (i \in \{ 0, \ldots, N \})$ 
such that the projectivization 
$\varphi: M \to P_{N} \mathbb{C}; 
p_{L}(\eta) \mapsto \pi_{\mathbb{C}_{N+1}}(\Phi(\eta))$ 
of $\Phi := (\iota_{L}(s_{0}, \ldots, \iota_{L}(s_{N})): 
P_{L} \to \mathbb{C}_{N}$ is an immersion, 
so that $\Phi$ is an immersion such that 
$\pi_{\mathbb{C}_{N+1}} \circ \Phi = \varphi \circ p_{L}$, 
$\Phi \circ R_{\lambda} = \lambda \cdot \Phi$ 
($\lambda \in \mathbb{C}^{\times}$) 
and $\Phi(P_{L}) \subseteq \mathbb{C}_{N+1} \backslash \{ 0 \}$ 
(cf. Lemma 3.1 (B)). 
Put 
\[
F := \{ (\xi, \eta) \in P_{N} \mathbb{C} \times \mathbb{C}_{N+1}; 
\pi_{\mathbb{C}_{N+1}}^{-1}(\xi) \cup \{ 0 \} \ni \eta \} 
\]
with the canonical projection 
$\pi_{F}: F \to P_{N} \mathbb{C}; 
(\xi, \eta) \mapsto \xi$ 
as {\it the tautological line bundle on $P_{N} \mathbb{C}$} 
after Y.~Matsushima \cite[(5.5.3)]{My1974} (cf. \cite[p.231]{FG2002}). 
Then 
$Q_{F} := F \backslash 0 
= \{ \widetilde{\eta} := (\pi_{\mathbb{C}_{N+1}}(\eta), \eta);
~\eta \in \mathbb{C}_{N+1} \backslash \{ 0 \} \}$. 
Put 
$\Pi_{2}: Q_{F} \to \mathbb{C}_{N+1} \backslash \{ 0 \}; 
\widetilde{\eta} \to \eta$ 
as a biholomorphism, so that 
$\Phi_{1} := \Pi_{2}^{-1} \circ \Phi: Q_{L^{*}} \to Q_{F}$ 
is an immersion such that 
$q_{L^{*}} \circ \Phi_{1} = \varphi \circ q_{F}$ 
and $\Phi_{1} \circ R_{\lambda} = R_{\lambda} \circ \Phi_{1}$ 
($\lambda \in \mathbb{C}^{\times}$). 
Then $\Phi_{2} 
:= \rho_{F} \circ \Phi_{1} \circ \rho_{L^{\otimes m}}: 
Q_{L} \to Q_{F^{*}}$ 
is an immersion such that 
$q_{F^{*}} \circ \Phi_{2} = \varphi \circ q_{L}$ 
and $\Phi_{2} \circ R_{\lambda} = R_{\lambda} \circ \Phi_{2}$ 
($\lambda \in \mathbb{C}^{\times}$). 
Hence, to show that $L$ is strictly positive, 
it is enough to show that $F^{*}$ is strictly positive
by the claim (B) (Then $L$ is posituve in the sense of K.~Kodaira 
if moreover $M$ is compact). 
It is proved as follows. 
Put 
$H := Q_{F} \times_{\rho_{-1}} \mathbb{C} 
= \{ [\widetilde{\eta}, w]_{-1}; \widetilde{\eta} \in Q_{F}, w \in \mathbb{C} \}$. 
By the claim (A), $H$ is isomorphic to $F^{*}$. 
{\it Claim. $H$ is strictly positive:} 
In fact, 
for 
$[\widetilde{\eta}_{j}, w_{j}]_{-1} \in H$ 
such as 
$\eta_{j} := (z_{0}(\eta_{j}), \ldots, z_{N}(\eta_{j})) 
\in \mathbb{C}_{N+1}$ and $w_{j} \in \mathbb{C}$ ($j \in \{ 1, 2 \}$) 
with $|| \eta_{j} || := 
\sqrt{\sum_{i = 0}^{N} z_{i}(\eta_{j}) \bar{z}_{i}(\eta_{j})}$, 
let $h$ be a hermitian fiber metric on $H$ 
defined as follows 
(cf. \cite[(17.5)]{Ns1981}): 
\[
h([\widetilde{\eta}_{1}, w_{1}]_{-1}, [\widetilde{\eta}_{2}, w_{2}]_{-1}) 
:= w_{1} \bar{w}_{2} / \sum_{i = 0}^{N} z_{i}(\eta_{1}) \bar{z}_{i}(\eta_{2}). 
\]
For $k \in \{ 0, \ldots, N \}$, 
put 
$U_{k} := \{ (z_{0}, \ldots, z_{N}) 
\in \mathbb{C}_{N+1} \backslash \{ 0 \}; z_{k} \ne 0 \}$, 
so that $P_{N} \mathbb{C} = \cup_{k = 0}^{N} \pi_{\mathbb{C}_{N+1}}(U_{k})$ 
and $\pi_{\mathbb{C}_{N+1}}^{-1}(\pi_{\mathbb{C}_{N+1}}(U_{k})) = U_{k}$. 
Put 
$\phi_{k}: \pi_{\mathbb{C}_{N+1}}(U_{k}) \to \mathbb{C}_{N}; 
\xi = [z_{0}, \ldots, z_{N}] 
\mapsto 
(z_{0}, \ldots, \hat{z}_{k}, \ldots, z_{N})/z_{k} 
=: (\zeta_{1}, \ldots, \zeta_{N})$ 
as the complex coordinates of $\pi_{\mathbb{C}_{N+1}}(U_{k})$. 
Put 
$\widetilde{e}_{k}: \widetilde{\pi}(U_{k}) \to Q_{H}; 
\xi = [z_{0}, \ldots, z_{N}] \mapsto 
[\widetilde{\eta}_{k}, 1]_{-1} = 
[(\pi_{\mathbb{C}_{N+1}}(\eta_{k}), \eta_{k}), 1]$ 
for 
$\eta_{k} = (\eta_{k 0}, \ldots, \eta_{k N}) = 
(z_{0}, \ldots, z_{N})/z_{k} \in \mathbb{C}_{N+1} \backslash \{ 0 \}$, 
so that 
$\widetilde{e}_{k} \in \Gamma_{\pi_{\mathbb{C}_{N+1}}(U_{k})}(Q_{H})$. 
Then 
$h^{\widetilde{e}_{k}} |_{\xi}
= - \log \circ h \circ \widetilde{e}_{k}^{\times 2} |_{\xi} 
= \log || \eta_{k} ||^{2} 
= \log (1 + \sum_{\ell = 1}^{N} \zeta_{\ell} \bar{\zeta}_{\ell})$ 
and 
\begin{eqnarray*}
& &
{\rm Lev}(h^{\widetilde{e}_{k}}) ( 
\sum_{\ell = 1}^{N} 
\left. 
X_{\ell} \frac{\partial}{\partial \zeta_{\ell}} 
\right|_{\xi}) 
= 
\left. 
\sum_{\alpha = 1}^{N} \sum_{\beta = 1}^{N} 
X_{\alpha} \bar{X}_{\beta} 
\frac{\partial^{2} 
\log(1 + \sum_{\ell = 1}^{N} \zeta_{\ell} \bar{\zeta}_{\ell})}
{\partial \zeta_{\alpha} \partial \zeta_{\beta}} \right|_{\xi}
\\
&=& \left. 
\sum_{\alpha = 1}^{N} \sum_{\beta = 1}^{N}
X_{\alpha} \bar{X}_{\beta} 
\left( \frac{\delta_{\alpha, \beta}}
{1 + \sum_{\ell = 1}^{N} \zeta_{\ell} \bar{\zeta}_{\ell}} 
- \frac{\bar{\zeta}_{\alpha} \zeta_{\beta}}
{(1 + \sum_{\ell = 1}^{N} \zeta_{\ell} \bar{\zeta}_{\ell})^{2}} \right) 
\right|_{\xi.} 
\end{eqnarray*}
Put 
$U(N) := 
\{ A \in GL_{N} \mathbb{C}; {}^{t} \bar{A} = A^{-1} \}$ 
and 
$R_{A}: U_{k} \to U_{k}; \xi \mapsto \phi_{k}^{-1}(\phi_{k}(\xi) A)$ 
for $A \in U(N)$, 
so that 
$h^{\widetilde{e}_{k}} \circ R_{A} = h^{\widetilde{e}_{k}}$. 
Take any 
$X = \sum_{\ell = 1}^{N} 
X_{\ell} \frac{\partial}{\partial \zeta_{\ell}} |_{\xi} 
\in (T' U_{k} \backslash \{ 0 \}) |_{\xi}$. 
There exists $A \in U(N)$ such that 
$R_{A*} X 
= c \cdot \left. \frac{\partial}{\partial \zeta_{1}} \right|_{R_{A} \xi}$ 
for $c := \sqrt{\sum_{\ell = 1}^{N} X_{\ell} \bar{X}_{\ell}} \ne 0$. 
Then 
\begin{eqnarray*}
{\rm Lev}(h^{\widetilde{e}_{k}}) ( X ) 
&=& {\rm Lev}(h^{\widetilde{e}_{k}} \circ R_{A}) ( X ) 
= {\rm Lev}(h^{\widetilde{e}_{k}}) (R_{A*} X) 
\\
&=& c^{2} \cdot \left. 
\frac{1 + \sum_{\ell = 2}^{N} \zeta_{\ell} \bar{\zeta}_{\ell}}
{1 + \sum_{\ell = 1}^{N} \zeta_{\ell} \bar{\zeta}_{\ell}} 
\right|_{R_{A} \xi} > 0 
\end{eqnarray*} 
Hence, $H$ is strictly positive, {\it as required}. 
\\ \
(D) 
By the claim (C), $K_{M}^{*}$ is positive in the sense of 
K.~Kodaira \cite{Kk1954}, i.e., $M$ is Fano. 
By virtue of S.~Kobayashi \cite{Ks1961} 
and S.T.~Yau \cite{Yst1977, Yst1978} 
(cf. H.~Nakajima \cite[p.61, Thm.2.3]{Nh1999}), 
it is well-known after A.L.~Besse \cite[11.26, 11.15]{Bal1987} 
or O.~Fujino \cite[Thm.6.1]{Fosa2014} 
that each connected component of 
a Fano manifold is simply connected. 
\qed 
\smallskip \\ \
{\bf 4. 
Local symplectifications and infinitesimal automorphisms.} 
Let $P$ be a holomorphic principal $\mathbb{C}^{\times}$-bundle 
on a complex manifold $M$ with the right action 
$R_{\lambda}$ ($\lambda \in \mathbb{C}^{\times}$) 
and the canonical projection $p: P \to M$. 
Let $W$ be an open subset of $P$, 
and $\theta$ be a holomorphic (1,0)-form on $W$ 
with $D_{\theta} := \{ X \in {T'} P{;~} \theta (X) = 0 \}$ 
such that 
$p_{*' \eta}^{-1}(0) \subseteq D_{\theta} |_{\eta} \ne {T'} P |_{\eta}$ 
($\eta \in P$). 
In the sense of \cite[pp.27-31]{OK1977}, 
let ${\cal L}_{{\cal E}_{P}}$ 
be the Lie derivative by the (1,0)-vector field ${\cal E}_{P}$. 
Assume that there exists 
$\epsilon_{\theta} \in \mathbb{C}^{\times}$ 
such that 
\[
{\cal L}_{{\cal E}_{P}} \theta 
:= \left. \frac{d}{d \lambda} R_{\lambda}^{*} \theta \right|_{\lambda = 1} 
= \epsilon_{\theta} \cdot \theta. 
\]
Then 
$(W, \theta)$ is said to be 
a {\it local pre-symplectification of degree $\epsilon_{\theta}$ 
in $P$ on $M$}, 
and put 
$E_{\theta} := p_{*'} D_{\theta} \subseteq T' M$ 
so that $p_{*'}^{-1} E_{\theta} = D_{\theta}$. 
Moreover, put 
\[
\flat 
:= \flat_{d \theta}: {T'} W \to {T'}^{*} W; 
v \mapsto \flat(v) := - \iota_{v} d \theta 
\] 
and 
$(d \theta)^{\perp} := \cup_{\eta \in W} 
\{ X \in {T'} W |_{\eta}{;~}d \theta(X_{1}, X) = 0
~(X_{1} \in {T'} W |_{\eta}) \}$. 
A local pre-symplectification 
$(W, \theta)$ of degree $\epsilon_{\theta}$ 
in $P$ on $M$ is said to be 
{\it a local symplectification of degree $\epsilon_{\theta}$ in $P$ on $M$} 
if and only if $(d \theta)^{\perp} = 0$. 
A local (pre-) symplectification $(P, \theta)$ 
of degree $\epsilon_{\theta}$ in $P$ on $M$ 
is simply said to be 
{\it a (pre-) symplectification of degree $\varepsilon_{\theta}$ on $M$}. 
For a c-manifold $(M, E)$, 
a local (pre-) symplectification $(W, \theta)$ 
of degree $\epsilon_{\theta}$ in $P$ on $M$ 
is said to be {\it related to $E$} or {\it $E$-related} 
if and only if $E_{\theta} = E$. 
\smallskip \\ \
\textsc{Lemma 4.0.} 
(A) 
{\it 
Let $(W, \theta)$ 
be a local pre-symplectification of degree 
$\epsilon_{\theta} \in \mathbb{C}^{\times}$ 
in $P$ on a complex manifold $M$. 
For $\sigma \in \Gamma_{\rm loc}(W)$, 
we define 
$\gamma_{\sigma}^{\theta} \in \Gamma_{V_{\sigma}} (T' W)$ 
by 
$\gamma_{\sigma}^{\theta} |_{\xi} 
:= \sigma^{*} (\theta |_{\phi_{\sigma}(\xi, 1)})$ 
$(\xi \in V_{\sigma})$, so that 
$(p^{*} \gamma_{\sigma}) X = \theta (X)$ 
for $X \in T' \sigma(V_{\sigma})$. 
Then one has the following results. 
}
\\ \
(a) 
{\it 
$E_{\theta}$ is a holomorphic distribution 
of complex codimension one on $p(W)$ 
such that 
$E_{\theta} |_{V_{\sigma}} 
= \{ Y \in T' V_{\sigma};~\gamma_{\sigma}^{\theta}(Y) = 0 \}$, 
$p_{*'\eta}^{-1} E_{\theta} = D_{\theta} |_{\eta}$, 
$p_{*\eta} D_{\theta} = E_{\theta} |_{p(\eta)}$, 
$R_{\lambda*} (D_{\theta} |_{R_{\lambda}^{-1} \eta}) = D_{\theta} |_{\eta}$, 
$R_{\lambda^{-1}}^{*} (\alpha_{D_{\theta}} |_{\eta}) 
= \alpha_{D_{\theta}} |_{R_{\lambda^{-1}} \eta}$, 
$R_{\lambda*} (\alpha_{D_{\theta}}^{\perp} |_{R_{\lambda^{-1}} \eta}) 
= \alpha_{D_{\theta}}^{\perp} |_{\eta}$ 
$(\eta \in W, \lambda \in \hat{\mathbb{C}}^{\times}_{\eta, W})$; 
$[{\cal E}_{P}, X] \in \Gamma_{\rm loc}(D_{\theta})~
(X \in \Gamma_{\rm loc}(D_{\theta}))$ 
and ${\cal E}_{P} |_{\eta} \in \alpha_{D_{\theta}}^{\perp}$ $(\eta \in W)$. 
Moreover, 
$\flat ({\cal E}_{P}) = - \epsilon_{\theta} \cdot \theta$, 
so that $d \theta ({\cal E}_{P}, X) = \epsilon_{\theta} \cdot \theta (X)$ 
for all $X \in T' W$. 
}
\\ \
(b) 
{\it 
For $\sigma \in \Gamma_{\rm loc}(W)$, 
$\alpha_{E_{\theta}}^{\perp} |_{\xi} 
= \{ Y \in E_{\theta} |_{\xi};~d \gamma_{\sigma}^{\theta} (Y_{1}, Y) = 0
~(Y_{1} \in E_{\theta} |_{\xi}) \}$ 
at $\xi \in V_{\sigma}$. 
At $\eta \in W$, 
$\alpha_{D_{\theta}}^{\perp} |_{\eta} 
= \{ X \in D_{\theta} |_{\eta};~d \theta (X_{1}, X) 
= 0~(X_{1} \in D_{\theta} |_{\eta}) \}$ 
and 
$p_{*'} (\alpha_{D_{\theta}}^{\perp} |_{\eta}) 
= \alpha_{E_{\theta}}^{\perp} |_{p(\eta)}$. 
} 
\\ \ 
(c) 
{\it 
At $\eta \in W$, 
$\alpha_{D_{\theta}}^{\perp} |_{\eta} 
= (d \theta)^{\perp} |_{\eta} \oplus \mathbb{C} {\cal E}_{P} |_{\eta}$. 
}
\\ \ 
(d) 
{\it 
For any $\eta \in W$, 
$\breve{p}_{\eta}: 
(d \theta)^{\perp} |_{\eta} \to \alpha_{E_{\theta}}^{\perp} |_{p(\eta)}; 
X \to p_{*} X$ 
is 
well-defibed as 
a $\mathbb{C}$-linear isomorphism. 
Especially, 
$(d \theta)^{\perp} |_{\eta} = 0$ 
if and only if $\alpha_{E_{\theta}}^{\perp} |_{p(\eta)} = 0$. 
Equivalently, 
$(W, \theta)$ 
is a local symplectification of degree $\epsilon_{\theta}$ 
in $P$ on $M$ if and only if 
$(p(W), E_{\theta})$ is a c-manifold. 
}
\\ \
(e) 
{\it 
If $\epsilon_{\theta} = \ell$ 
for some $\ell \in \mathbb{C}^{\times} \cap \mathbb{Z}$, 
and that 
$W |_{p(\eta)}$ is connected for all $\eta \in W$, 
then 
$R_{\lambda}^{*} (\theta |_{R_{\lambda} \eta}) 
= \lambda^{\ell} \cdot \theta |_{\eta}$ 
for all $(\eta, \lambda) \in W \times \breve{\mathbb{C}}^{\times}_{\eta, W}$. 
}
\\ \
(B) 
(cf. \cite[p.145]{Bwm1961}, \cite{Avi1978}, \cite{OK1977}) 
{\it 
Let $(M, E)$ be a c-manifold $(M, E)$. 
On 
$(P_{E}, p_{E}) := (P_{L_{E}}, p_{L_{E}})$, 
put 
\[
\theta_{E}: {T'} P_{E} \to \mathbb{C}; 
X \mapsto \theta_{E}(X) 
:= <\!\! \Pi_{P_{E}}(X), \varpi_E(p_{E*} X) \!\!>_{L_{E}}. 
\] 
Then 
$(P_{E}, \theta_{E})$ 
is a symplectification of degree $1$ on $M$ 
related to $E$, 
which is said to be the canonical symplectification of $(M, E)$. 
}
\smallskip \\ \
{\it Proof}. 
(A) 
(a) 
Take any $\eta \in W$. 
Then 
$p_{*'\eta}^{-1} E_{\theta} 
= p_{*'\eta}^{-1} (p_{*'\eta} D_{\theta}) 
= p_{*'\eta}^{-1}(0) + D_{\theta} |_{\eta} 
= D_{\theta} |_{\eta}$. 
Since $p_{*'\eta}$ is surjective, 
$E_{\theta} |_{p(\eta)} = p_{*' \eta}(p_{*'\eta}^{-1} E_{\theta}) 
= p_{*'\eta} D_{\theta}$. 
For $\lambda \in \hat{\mathbb{C}}^{\times}_{\eta, W}$, 
$D_{\theta} |_{\eta} 
= p_{*' \eta}^{-1} E_{\theta} 
= R_{\lambda*} (p_{*'R_{\lambda}^{-1} \eta}^{-1} E_{\theta}) 
= R_{\lambda*} (D_{\theta} |_{R_{\lambda}^{-1} \eta})$ 
by $p_{*'} \circ R_{\lambda *} = p_{*'}$. 
It is obvious that 
$R_{\lambda^{-1}}^{*} (\alpha_{D_{\theta}} |_{\eta}) 
= \alpha_{D_{\theta}} |_{R_{\lambda^{-1}} \eta}$ 
and 
$R_{\lambda*} (\alpha_{D_{\theta}}^{\perp} |_{R_{\lambda^{-1}} \eta}) 
= \alpha_{D_{\theta}}^{\perp} |_{\eta}$. 
For $X \in \Gamma_{\rm loc}(D_{\theta})$, 
$[{\cal E}_{P}, X] |_{\eta} 
= \lim_{\lambda \to 1} 
\frac{1}{\lambda - 1} 
(R_{\lambda}^{*} (X |_{R_{\lambda} \eta}) - X |_{\eta}) 
= \lim_{\lambda \to 1} 
\frac{1}{\lambda - 1} 
(R_{\lambda*}^{-1} (X |_{R_{\lambda} \eta}) - X |_{\eta}) 
\in D_{\theta} |_{\eta.}$ 
Hence, 
$[{\cal E}_{P}, X] \in \Gamma_{\rm loc}(D_{\theta})$. 
By 
${\cal E}_{P} \in p_{*'}^{-1}(0) \subseteq D_{\theta}$ 
and 
$\alpha_{D_{\theta}}({\cal E}_{P}, X) 
= \varpi_{D_{\theta}}([{\cal E}_{P}, X]) = 0$ 
($X \in \Gamma_{\rm loc}(D_{\theta}))$, 
one has that 
${\cal E}_{P} |_{\eta} \in \alpha_{D_{\theta}}^{\perp}$ 
as required. 
By H.~Cartan's formula (cf. \cite[p.29]{OK1977}), 
$\flat ({\cal E}_{P} |_{W}) 
= - \iota_{{\cal E}_{P}} d \theta  
= - {\cal L}_{{\cal E}_{P}} \theta  
+ d (\iota_{{\cal E}_{P}} \theta) 
= - \epsilon_{\theta} \cdot \theta$ 
as required. 
Take any $\sigma \in \Gamma_{\rm loc}(W)$ and $\xi \in V_{\sigma}$. 
For any $Y \in E_{\theta} |_{\xi}$, 
there exists 
$X \in D_{\theta} |_{\sigma'(\xi)}$ 
such that $Y = p_{*'} X$. 
For 
$X_{1} := {\sigma}_{*} Y \in T' \sigma(V_{\sigma})$, 
$p_{*'} X_{1} = Y$ 
and 
$X_{1} - X \in p_{*'}^{-1}(0) \subseteq D_{\theta}$, 
so that $X_{1} \in D_{\theta}$. 
Then 
$0 = \theta(X_{1}) = \theta({\sigma}_{*} Y) 
= \gamma_{\sigma}^{\theta} Y$. 
Hence, 
$E_{\theta} \subseteq 
\{ Y \in T' V_{\sigma} {;~} \gamma_{\sigma}^{\theta} Y = 0 \}$. 
Conversely, for any $Y \in T' V_{\sigma}$ such as 
$\gamma_{\sigma}^{\theta} Y = 0$, 
put $X := {\sigma}_{*} Y \in T' W$. 
Then 
$\theta(X) = \gamma_{\sigma}^{\theta} Y = 0$ 
and $Y = p_{*'} X \in p_{*'} D_{\theta} = E_{\theta}$. 
Hence, 
$E_{\theta} = 
\{ Y \in T' V_{\sigma} {;~} \gamma_{\sigma}^{\theta} Y = 0 \}$. 
{\it 
Claim. 
${T'} M |_{\xi} \ne E_{\theta} |_{\xi}$ 
($\xi \in p(W)$):} 
{\it If not}, 
${T'} M |_{\xi} = E_{\theta} |_{\xi}$ 
for some $\xi \in p(W)$, 
so that 
${T'} P |_{\eta} 
= p_{*'\eta}^{-1} E_{\theta} 
= D_{\theta} |_{\eta} \ne {T'} P |_{\eta}$ 
for $\eta \in p^{-1}(\xi)$, 
{\it a contradiction}. 
Then 
$0 \ne \gamma_{\sigma}^{\theta} |_{\xi} 
\in \Gamma_{M}({T'} V_{\sigma}) |_{\xi}$ 
for any $\xi \in V_{\sigma'}$, so that 
$E_{\theta}$ 
is a holomorphic distribution of complex codimension one on $p(W)$. 
\\ \ 
(b) 
For $\sigma \in \Gamma_{\rm loc}(W)$, 
take any $\xi \in V_{\sigma}$. 
By the claim (a), there exists
$Y_{0} \in (T' M \backslash E_{\theta}) |_{\xi} \backslash \{ 0 \}$. 
Take any $Y_{1}, Y_{2} \in E_{\theta} |_{\xi}$. 
There exist a restriction 
$\sigma_{\xi} := \sigma |_{V_{\sigma_{\xi}}} 
\in \Gamma_{\rm loc}(V_{\sigma})$ 
of $\sigma$ such as $\xi \in V_{\sigma_{\xi}}$, 
$\widetilde{Y}_{0} \in \Gamma_{V_{\sigma_{\xi}}}(T' M \backslash E_{\theta})$ 
and $\widetilde{Y}_{1}, \widetilde{Y}_{2} \in \Gamma_{V_{\sigma_{\xi}}}(E_{\theta})$ 
such that $\xi \in V_{\sigma_{\xi}}$ and 
$Y_{j} = \widetilde{Y}_{j} |_{\xi}$ ($j \in \{ 0, 1, 2 \}$). 
Put 
$g: V_{\sigma_{\xi}} \to \mathbb{C}; 
\widetilde{\xi} \mapsto 
\gamma_{\sigma}^{\theta} ([\widetilde{Y}_{1}, \widetilde{Y}_{2}])/ 
\gamma_{\sigma}^{\theta} (\widetilde{Y}_{0}) |_{\widetilde{\xi}.}$ 
Then 
$[\widetilde{Y}_{1}, \widetilde{Y}_{2}] - g \cdot \widetilde{Y}_{0} 
\in \Gamma_{V_{\sigma_{\xi}}}(E_{\theta})$ 
and 
$d \gamma_{\sigma}^{\theta}(\widetilde{Y}_{1}, \widetilde{Y}_{2}) 
= \widetilde{Y}_{1} \gamma_{\sigma}^{\theta}(\widetilde{Y}_{2}) 
- \widetilde{Y}_{2} \gamma_{\sigma}^{\theta}(\widetilde{Y}_{1}) 
- \gamma_{\sigma}^{\theta}([\widetilde{Y}_{1}, \widetilde{Y}_{2}]) 
= - g \cdot \gamma_{\sigma}^{\theta}(\widetilde{Y}_{0})$. 
Hence, 
\[
- d \gamma_{\sigma}^{\theta}(\widetilde{Y}_{1}, \widetilde{Y}_{2}) \cdot 
\varpi_{E_{\theta}}(\widetilde{Y}_{0})/\gamma_{\sigma}(\widetilde{Y}_{0}) 
= g \cdot \varpi_{E_{\theta}}(\widetilde{Y}_{0}) 
= \varpi_{E_{\theta}}([ \widetilde{Y}_{1}, \widetilde{Y}_{2} ]) 
= \alpha_{E_{\theta}}(\widetilde{Y}_{1}, \widetilde{Y}_{2}), 
\]
which gives the first equation. 
For $j \in \{ 0, 1, 2 \}$ and 
$(\widetilde{\xi}, \lambda) \in V_{\sigma_{\xi}} 
\times \breve{\mathbb{C}}^{\times}_{\sigma(\widetilde{\xi}), W}$, 
put 
$\widetilde{X}_{j} |_{\phi_{\sigma} (\widetilde{\xi}, \lambda)} 
:= R_{\lambda*} (\sigma_{* \widetilde{\xi}} (\widetilde{Y}_{j}))$. 
Then  
$p_{*'} \widetilde{X}_{j} = \widetilde{Y}_{j}$, 
$\widetilde{X}_{0} \in \Gamma_{V_{\sigma_{\xi}}}(T' P \backslash D_{\theta})$ 
and 
$\widetilde{X}_{k} \in \Gamma_{V_{\sigma_{\xi}}}(D_{\theta})$ 
($k \in \{ 1, 2 \}$) 
because of 
$R_{\lambda*} (D_{\theta} |_{\sigma(\widetilde{\xi})}) 
= D_{\theta} |_{R_{\lambda} \sigma(\widetilde{\xi})}$ 
by (a). 
Put 
$h: P |_{V_{\sigma_{\xi}}} \to \mathbb{C}; 
\eta \mapsto 
\theta ([\widetilde{X}_{1}, \widetilde{X}_{2}])/ \theta(\widetilde{X}_{0})$. 
Then 
$d \theta (\widetilde{X}_{1}, \widetilde{X}_{2}) 
= - \theta ([\widetilde{X}_{1}, \widetilde{X}_{2}]) 
= - h \cdot \theta(\widetilde{X}_{0})$ 
and that 
$[\widetilde{X}_{1}, \widetilde{X}_{2}] - h \cdot \widetilde{X}_{0} 
\in \Gamma_{P |_{V_{\sigma_{\xi}}}}(D_{\theta})$, 
so that 
$\alpha_{D_{\theta}}(\widetilde{X}_{1}, \widetilde{X}_{2}) 
= h \cdot \varpi_{D_{\theta}}(\widetilde{X}_{0}) 
= - d \theta (\widetilde{X}_{1}, \widetilde{X}_{2}) \cdot 
\varpi_{D_{\theta}}(\widetilde{X}_{0})/\theta(\widetilde{X}_{0})$. 
It follows from (a) that 
\begin{eqnarray*}
& &d \theta (\widetilde{X}_{1} + c_{1} {\cal E}_{P}, 
\widetilde{X}_{2} + c_{2} {\cal E}_{P}) 
= d \theta (\widetilde{X}_{1}, \widetilde{X}_{2}) 
\\
&=& - (\theta(\widetilde{X}_{0})/\varpi_{D_{\theta}}(\widetilde{X}_{0})) 
\cdot \alpha_{D_{\theta}}(\widetilde{X}_{1}, 
\widetilde{X}_{2}) 
\\
&=& - (\theta(\widetilde{X}_{0})/\varpi_{D_{\theta}}(\widetilde{X}_{0})) 
\cdot \alpha_{D_{\theta}}(\widetilde{X}_{1} + c_{1} {\cal E}_{P}, 
\widetilde{X}_{2} + c_{2} {\cal E}_{P}), 
\end{eqnarray*} 
which gives the second equation because of  
$D_{\theta} |_{\phi_{\sigma}(\widetilde{\xi}, \lambda)} 
= p_{*'\phi_{\sigma}(\widetilde{\xi}, \lambda)}^{-1} E_{\theta} 
= R_{\lambda*} \sigma_{*\widetilde{\xi}} E_{\theta} 
+ \mathbb{C} {\cal E}_{P} |_{\phi_{\sigma}(\widetilde{\xi}, \lambda)} 
~(\widetilde{\xi} \in V_{\sigma}, 
\lambda \in \breve{\mathbb{C}}^{\times}_{\sigma(\widetilde{\xi}), W})$ 
by the claim (a). 
Let 
$Y_{1} \in \alpha_{E_{\theta}}^{\perp} |_{\xi.}$ 
By 
$p_{*'\sigma(\xi)} \widetilde{X}_{k} = Y_{k} 
\in E_{\theta} |_{\xi}$ ($k \in \{ 1, 2 \}$), 
\[
d \theta (\widetilde{X}_{1}, \widetilde{X}_{2} + c_{2} {\cal E}_{P}) 
|_{\sigma(\xi)} 
= d \theta (\widetilde{X}_{1}, \widetilde{X}_{2}) 
|_{\sigma(\xi)} 
= \sigma^{*} d \theta (Y_{1}, Y_{2}) 
= d \gamma_{\sigma}(Y_{1}, Y_{2}) = 0 
\]
for any $Y_{2} \in E_{\theta}$. 
Then 
$\widetilde{X}_{1} |_{\sigma(\xi)} \in \alpha_{D_{\theta}}^{\perp}$ 
and 
$Y_{1} = p_{*'} 
(\widetilde{X}_{1} |_{\sigma(\xi)}) 
\in p_{*'\sigma(\xi)} \alpha_{D_{\theta}}^{\perp}$. 
Hence, 
$\alpha_{E_{\theta}}^{\perp} |_{\xi} \subseteq 
p_{*'\sigma(\xi)} \alpha_{D_{\theta}.}^{\perp}$ 
Conversely, let 
$X := (\widetilde{X}_{1} + c_{1} {\cal E}_{P}) |_{\sigma(\xi)} 
\in \alpha_{D_{\theta}}^{\perp} |_{\sigma(\xi)}$. 
Then 
$0 = d \theta (X, \widetilde{X}_{2}) |_{\sigma(\xi)} 
= d \theta (\widetilde{X}_{1}, \widetilde{X}_{2}) |_{\sigma(\xi)} 
= \sigma^{*} d \theta (Y_{1}, Y_{2}) 
= d \gamma_{\sigma}(Y_{1}, Y_{2})$, 
so that 
$Y_{1} = p_{*'} X \in \alpha_{E_{\theta}}^{\perp} |_{\xi}$. 
Hence, 
$\alpha_{E_{\theta}}^{\perp} |_{\xi} 
\supseteq p_{*'\sigma(\xi)} \alpha_{D_{\theta}.}^{\perp}$ 
Then 
$\alpha_{E_{\theta}}^{\perp} |_{\xi} 
= p_{*'\sigma(\xi)} \alpha_{D_{\theta}}^{\perp} 
= p_{*'\phi_{\sigma}(\xi, \lambda)} (R_{\lambda*} \alpha_{D_{\theta}}^{\perp}) 
= p_{*'\phi_{\sigma(\xi, \lambda)}} \alpha_{D_{\theta}}^{\perp}$ 
for $\lambda \in \mathbb{C}^{\times}$ by the claim (a). 
Hence, $\alpha_{E_{\theta}}^{\perp} |_{p(\eta)} 
= p_{*\eta} \alpha_{D_{\theta}}^{\perp}$ 
for $\eta \in W$.  
\\ \
(c) 
$\{ X \in T' W; d \theta ({\cal E}_{P}, X) = 0 \} 
= \{ X \in T' W; \epsilon_{\theta} \cdot \theta (X) = 0 \} 
= D_{\theta}$ 
because of the claim (a) and $\epsilon_{\theta} \ne 0$, 
so that 
$(d \theta)^{\perp} \subseteq D_{\theta}$. 
By the claim (b), 
$\alpha_{D_{\theta}}^{\perp} \supseteq 
(d \theta)^{\perp} \cap D_{\theta} = (d \theta)^{\perp}$. 
Take any $\eta \in W$. 
By 
$D_{\theta} |_{\eta} \ne T' W |_{\eta}$, 
there exists $X_{0} \in T' W |_{\eta}$ such as $\theta(X_{0}) \ne 0$, 
so that 
$d \theta ({\cal E}_{P}, X_{0}) 
= \epsilon_{\theta} \cdot \theta (X_{0}) \ne 0$. 
Hence, 
$(d \theta)^{\perp} |_{\eta} \cap \mathbb{C} {\cal E}_{P} |_{\eta} 
= \{ 0 \}$. 
Put 
$\hat{\eta}: \alpha_{D_{\theta}}^{\perp} |_{\eta} \to \mathbb{C}; 
X \mapsto d \theta (X, X_{0})$, 
which is a $\mathbb{C}$-linear map such that 
$\hat{\eta}^{-1}(0) \supseteq (d \theta)^{\perp} |_{\eta}$. 
Conversely, 
$\hat{\eta}^{-1}(0) \subseteq (d \theta)^{\perp} |_{\eta}$ 
by the claim (b) and 
$T' W |_{\eta} = D_{\theta} |_{\eta} + \mathbb{C} \widetilde{X}_{0} |_{\eta}$. 
Hence, 
$\hat{\eta}^{-1}(0) = (d \theta)^{\perp} |_{\eta}$. 
Because of 
$\hat{\eta}({\cal E}_{P}) = \epsilon_{\theta} \cdot \theta (X_{0}) \ne 0$, 
${\rm Im}( \hat{\eta} ) = \mathbb{C}$. 
Then 
${\rm dim}_{\mathbb{C}} 
\alpha_{D_{\theta}}^{\perp} |_{\eta} 
= {\rm dim}_{\mathbb{C}} \hat{\eta}^{-1}(0) 
+ {\rm dim}_{\mathbb{C}} {\rm Im}( \hat{\eta} ) 
= {\rm dim}_{\mathbb{C}} 
((d \theta)^{\perp} |_{\eta} \oplus \mathbb{C} {\cal E}_{P} |_{\eta})$, 
so that 
$\alpha_{D_{\theta}}^{\perp} |_{\eta} 
= (d \theta)^{\perp} |_{\eta} \oplus \mathbb{C} {\cal E}_{P} |_{\eta}$. 
\\ \ 
(d) 
Because of the claims (b) and (c), the correspondence 
$\breve{p}_{\eta}$ 
is well-defined as a $\mathbb{C}$-linear surjection. 
Moreover, 
$\breve{p}_{\eta}^{-1}(0) 
= (d \theta)^{\perp} |_{\eta} \cap p_{*'\eta}^{-1}(0) 
= (d \theta)^{\perp} |_{\eta} \cap \mathbb{C} {\cal E}_{P} |_{\eta} 
= \{ 0 \}$ 
by (c), 
so that 
$\breve{p}_{\eta}$ is a $\mathbb{C}$-linear isomorphism, 
as required. 
\\ \
(e) 
Take any $\eta \in W$. 
Put 
$U := \breve{\mathbb{C}}^{\times}_{\eta, W}$ 
and a holomorphic map 
$g: U \to (T' W)^{*} |_{\eta}; 
\lambda \to g(\lambda) := 
R_{\lambda}^{*} (\theta |_{R_{\lambda} \eta}) 
- \lambda^{\ell} \cdot \theta |_{\eta}$ 
into $(T' W)^{*} |_{\eta} \cong \mathbb{C}_{{\rm n}_{W}}$. 
For $\lambda \in U$, 
$\ell \cdot R_{\lambda}^{*} (\theta |_{R_{\lambda} \eta}) 
= R_{\lambda}^{*}(({\cal L}_{{\cal E}_{P}} \theta) |_{R_{\lambda} \eta}) 
= R_{\lambda}^{*}( \left. 
\frac{d}{d \lambda_{1}} R_{\lambda_{1}}^{*} 
\theta |_{R_{\lambda_{1} \lambda} \eta} \right|_{\lambda_{1} = 1}) 
= \left. \lambda \cdot 
\frac{d}{d \lambda'} R_{\lambda}^{*} 
\theta |_{R_{\lambda'} \eta} \right|_{\lambda' = \lambda} 
= \lambda \cdot \frac{d}{d \lambda} R_{\lambda}^{*} 
\theta |_{R_{\lambda} \eta}$, 
so that 
$\lambda \cdot \frac{d}{d \lambda} g(\lambda) 
= \ell \cdot g(\lambda)$. 
Since $W |_{\xi}$ 
is connected and open in $P |_{\xi}$ ($\xi \in p(W)$), 
$U$ is a connected open neighbourhood of $1$ in $\mathbb{C}^{\times}$ 
by Lemma 2.2 (B). 
By Lemma 2.2 (A) for $U \ni 1$ and $(g, m, k) := (g, 1, \ell)$, 
one has that 
$g(\lambda) = g(1) \cdot \lambda^{\ell} 
= 0 \cdot \lambda^{\ell} = 0$ ($\lambda \in U$) 
as required. 
\\ \ 
(B) 
Since 
$\Pi_{P_{E}}$, $p_{E*}$, $\varpi_E$ and $<\!\!\!*,~*\!\!\!>_{L_{E}}$ 
are holomorphic, so is $\theta_{E}$. 
For $\lambda \in \mathbb{C}^{\times}$, 
$R_{\lambda}^{*} \theta_{E} 
= <\!\! \Pi_{P_{E}} \circ R_{\lambda*}, 
(\varpi_E \circ p_{E*}) \circ R_{\lambda*} \!\!>_{L_{E}} 
= <\!\! R_{\lambda} \circ \Pi_{P_{E}}, \varpi_E \circ p_{E*} \!\!>_{L_{E}} 
= \lambda \cdot \theta_{E}$. 
For $(\eta, X) \in P_{E} \times {T'}_{\eta} P_{E}$, 
$X \in D_{\theta_{E}}$ 
if and only if 
$\theta_{E}(X) = 0$, 
i.e., 
$p_{E*\eta} X \in E |_{p_{E}(\eta)}  
\ne {T'} M |_{p_{E}(\eta)} = p_{E*\eta} {T'} P_{E}$, 
i.e., 
$X \in p_{E*'\eta}^{-1} E \ne {T'} P_{E} |_{\eta}$. 
Hence, 
$D_{\theta_{E}} |_{\eta} = p_{E*'\eta}^{-1} E \ne {T'} P_{E} |_{\eta}$. 
Moreover, 
$p_{E*'} D_{\theta_{E}} = p_{E*'} (p_{E*'\eta}^{-1} E) = E$ 
since $p_{E*'\eta}$ is surjective. 
Then 
$(P_{E}, \theta_{E})$ is a pre-symplectification of degree $1$ 
on $M$ related to $E$. 
By the claim (A) (d), $(P_{E}, \theta_{E})$ 
is a symplectification of degree $1$ on $M$ related to $E$. 
~\qed 
\smallskip \\ \
Let 
$(W, \theta)$ 
be a local symplectification of degree 
$\epsilon_{\theta} \in \mathbb{C}^{\times}$ 
in a $\mathbb{C}^{\times}$-principal bundle $P$ 
on a complex manifold $M$ with the projection 
$p: P \to M$. 
By $(d \theta)^{\perp} = 0$, 
$\flat 
= \flat_{d \theta}: {T'} P \to {T'}^{*} P; 
v \mapsto \flat(v) = - \iota_{v} d \theta$ 
is a holomorphic vector bundle isomorphism. 
Then 
\[
\natural := \natural_{d \theta}: 
{T'}^{*} P \to {T'} P; \omega \mapsto \natural \omega;~
\omega =: \flat (\natural \omega) = -\iota_{\natural \omega} d \theta 
\]
is well-defined as the inverse map of $\flat_{d \theta}$. 
By means of the Lie derivative ${\cal L}_{X}$ 
by a (1,0)-vector field $X$ 
(cf. \cite[p.29]{KN1963}, \cite[p.28]{OK1977}), 
put 
$a(W, \theta) := \{ X \in a(W); {\cal L}_{X} \theta = 0 \}$. 
Let 
$H^{\theta}_{f} := \natural (d f)$ 
be {\it the Hamiltonian vector field of} $f \in {\cal O}(W)$. 
For $\delta \in \mathbb{C}$, put 
\[
\widetilde{H}^{\delta}_{\theta}: 
{\cal O}^{\delta}(W) \to 
a^{\delta}(W) := \{ H^{\theta}_{f}; f \in {\cal O}^{\delta} (W) \}; 
f \mapsto H^{\theta}_{f}. 
\]
For 
$f_{i} \in {\cal O}(W)$ $(i \in \{ 1, 2 \})$, 
put 
{\it the Poisson bracket} 
$\{ f_{1} \mid f_{2} \} 
:= d \theta (H^{\theta}_{f_{1}}, H^{\theta}_{f_{2}})$ 
(cf. \cite[(2.5)]{OK1977}). 
The local symplectification $(W, \theta)$ of degree $\epsilon_{\theta}$ 
in $P$ on $M$ 
is said to be {\it infinitesimally homogeneous} 
if and only if 
$T' W |_{\eta} = \{ X |_{\eta}; X \in a(W, \theta) \}$ 
at each $\eta \in W$. 
For a subset $S$ of $W$, a map 
$X: S \to T' W; \eta \mapsto X |_{\eta}$ 
such as $X |_{\eta} \in T' W |_{\eta}$ 
($\eta \in S$) 
is said to be {\it projectable} 
when there exists a map $Y: p(S) \to T' M$ such that 
$Y |_{\xi} \in T' M |_{\xi}$ 
($\xi \in p(S)$) and 
$p_{*\eta} X = Y |_{p(\eta)}$ 
for all $\eta \in S$. In this case, $Y$ is uniquely determined by $X$. 
\smallskip \\ \
\textsc{Lemma 4.1.} 
(A) 
${\cal E}_{P} |_{W} 
= - \epsilon_{\theta} \cdot \natural \theta \in \Gamma_{W}(T' W)$. 
\\ \
(B) 
{\it 
${\cal O}^{\delta}(W) 
= \{ f \in {\cal O}(W){;~} 
\theta(H^{\theta}_{f}) = \frac{\delta}{\epsilon_{\theta}} f \} 
\subseteq \{ f \in {\cal O}(W){;~} 
{\cal L}_{H^{\theta}_{f}} \theta 
= \frac{\delta - \epsilon_{\theta}}{\epsilon_{\theta}} d f \}$ 
for $\delta \in \mathbb{C}$, 
and that 
$a^{\epsilon_{\theta}}(W) = a(W, \theta) 
= \{ X \in a(W); H^{\theta}_{\theta(X)} = X \}$. 
} 
\\ \
(C) 
{\it 
If 
$\delta \in \mathbb{C}^{\times}$, 
then 
$\widetilde{H}^{\delta}_{\theta}$ 
is a $\mathbb{C}$-linear isomorphism with the inverse map 
\[
\widetilde{\theta}^{\delta}: 
a^{\delta}(W) \to {\cal O}^{\delta}(W); 
X \mapsto \frac{\epsilon_{\theta}}{\delta} \cdot \theta(X). 
\]
}
(D) 
{\it 
${\cal O}^{\delta}(W) 
\subseteq 
\{ f \in {\cal O}(W){;~} 
[{\cal E}_{P}, H^{\theta}_{f}] 
= (\delta - \epsilon_{\theta}) H^{\theta}_{f} \}$. 
Especially, 
$a(W, \theta) \subseteq \{ X \in a(W); [{\cal E}_{P}, X] = 0 \}$. 
If 
$f_{i} \in {\cal O}^{\delta_{i}}(W)$ 
with 
$\delta_{i} \in \mathbb{C}$ $(i \in \{ 1, 2 \})$, 
then 
$\{ f_{1} \mid f_{2} \} 
\in {\cal O}^{\delta_{1} + \delta_{2} - \epsilon_{\theta}}(W)$. 
} 
\\ \
(E) 
(a) 
{\it 
Let $X: W \to T' W; \eta \mapsto X |_{\eta}$ 
be a projectable vector field. 
If $X \in a(W, \theta)$, then $p_{*} X \in a(p(W), E_{\theta})$. 
}
\\ \
(b) 
{\it 
Let $S$ be a subset of $W$ 
such that $S |_{p(\eta)}$ is a connected open subset of $W |_{p(\eta)}$ 
for all $\eta \in S$. 
For a map 
$X: S \to T' W; \eta \mapsto X |_{\eta} \in T' W |_{\eta}$ 
and $\eta \in S$, put 
$\hat{\eta}_{X}: 
\hat{\mathbb{C}}^{\times}_{\eta, S} \to T' W |_{\eta}; 
\lambda \mapsto \hat{\eta}_{X} (\lambda) := 
R_{\lambda *} (X |_{R_{\lambda}^{-1} \eta}) \in T' W |_{\eta}$. 
For any $\eta \in S$, 
assume that 
$\left. \frac{d}{d \lambda} 
\hat{\eta}_{X} (\lambda) \right|_{\lambda = 1} = 0$. 
Then $X$ is projectable. 
}
\\ \
(c) 
{\it 
Assume that the fiber $W |_{p(\eta)}$ 
is connected for every $\eta \in W$. 
Then a map 
\[
\widetilde{p}_{*}: 
a(W, \theta) \to a(p(W), E_{\theta}); 
X \mapsto p_{*} X
\]
is well-defined as a complex Lie algebra homomorphism. 
}
\\ \
(d) 
(S.~Kobayashi \cite[p.30, Thm.I.7.1]{Ks1972}) 
{\it 
Let $(M, E)$ be a c-manifold. 
Put 
\[
\widetilde{\varpi}_{E}: a(M, E) \to \Gamma_{M}(L_{E}); 
Y \mapsto \varpi_{E}(Y)
\]
with $\varpi_{E}(Y): M \to L_{E}; \xi \mapsto \varpi_{E}(Y |_{\xi})$. 
Then 
$\widetilde{\varpi}_{E}$ is a $\mathbb{C}$-linear isomorphism. 
}
\\ \ 
~(e) 
{\it 
$\widetilde{p}_{E*}: a^{1}(P_{E}) \to a(M, E); 
X \mapsto p_{*} X$ 
is a complex Lie algebra isomorphism 
with the inverse map 
${\cal P}_{\theta_{E}} := \widetilde{H}^{1}_{\theta_{E}} 
\circ \widetilde{\iota}_{L_{E}} \circ \widetilde{\varpi}_{E}: 
a(M, E) \to a^{1}(P_{E})$. 
}
\smallskip \\ \
{\it Proof}. 
(A) 
By Lemma 4.0 (A) (a), 
$\flat ({\cal E}_{P} |_{W}) = - \epsilon_{\theta} \cdot \theta$, 
so that 
${\cal E}_{P} |_{W} = - \epsilon_{\theta} \cdot \sharp \theta$. 
\\ \
(B) 
For any $f \in {\cal O}(W)$, 
$\theta(H^{\theta}_{f}) = (- \iota_{\natural \theta} d \theta) H^{\theta}_{f} 
= (\iota_{H^{\theta}_{f}} d \theta) (\natural \theta) 
= - (d f) (\natural \theta) 
= - (\natural \theta) f = {\cal E}_{P} f / \epsilon_{\theta}$ 
because of the claim (A). 
Hence, 
$\theta(H^{\theta}_{f}) = \frac{\delta}{\epsilon_{\theta}} f$ 
if and only if $f \in {\cal O}^{\delta} (W)$. 
For any $f \in {\cal O}(W)$, 
${\cal L}_{H^{\theta}_{f}} \theta 
= \iota_{H^{\theta}_{f}} d \theta + d (\theta (H^{\theta}_{f})) 
= - d f + d (\theta(H^{\theta}_{f}))$, 
which equals 
$\frac{- \epsilon_{\theta} + \delta}{\epsilon_{\theta}} d f$ 
when $f \in {\cal O}^{\delta}(W)$, 
so that 
${\cal L}_{H^{\theta}_{f}} \theta = 0$ 
when $f \in {\cal O}^{\epsilon_{\theta}}(W)$. 
Hence, 
$a^{\epsilon_{\theta}}(W) \subseteq a(W, \theta)$. 
On the other hand, 
$a^{\epsilon_{\theta}}(W) 
= \{ H^{\theta}_{f}; \theta(H^{\theta}_{f}) = f \in {\cal O}(W) \} 
\supseteq 
\{ X \in a(W);~H^{\theta}_{\theta(X)} = X \}$. 
For any $X \in a(W, \theta)$, 
$d (\theta(X)) + \iota_{X} d \theta 
= {\cal L}_{X} \theta = 0$, so that 
$X = H^{\theta}_{\theta(X)}$. 
Hence, 
$a(W, \theta) \subseteq 
\{ X \in a(W);~H^{\theta}_{\theta(X)} = X \}$, 
so that 
$a^{\epsilon_{\theta}}(W) = a(W, \theta) 
= \{ X \in a(W); H^{\theta}_{\theta(X)} = X \}$. 
\\ \ 
(C) 
Because of the step (A), 
$\widetilde{\theta}^{\delta} \circ \widetilde{H}^{\delta}_{\theta} 
= {\rm id}_{{\cal O}^{\delta}(W)}$, 
so that $\widetilde{H}^{\delta}_{\theta}$ is injective. 
By definition, 
$\widetilde{H}^{\delta}_{\theta}$ is surjective. 
Then 
$\widetilde{H}^{\delta}$ is a $\mathbb{C}$-linear isomorphism 
with the inverse map $\breve{\theta}^{\delta}$. 
\\ \
(D) 
For 
$f \in {\cal O}(W)$, 
${\cal L}_{{\cal E}_{P}} d f 
= - {\cal L}_{{\cal E}_{P}} (\iota_{H^{\theta}_{f}} d \theta) 
= - \iota_{H^{\theta}_{f}} ({\cal L}_{{\cal E}_{P}} d \theta) 
- [ {\cal L}_{{\cal E}_{P}}, \iota_{H^{\theta}_{f}}] d \theta 
= - \iota_{H^{\theta}_{f}} (\epsilon_{\theta} \cdot d \theta) 
- \iota_{{\cal L}_{{\cal E}_{P}} H^{\theta}_{f}} d \theta 
= \epsilon_{\theta} \cdot d f 
- \iota_{[{\cal E}_{P}, H^{\theta}_{f}]} d \theta$ 
by the standard formula \cite[Prop.I.3.10 (b)]{KN1963}. 
If $f \in {\cal O}^{\delta}(W)$, then 
${\cal L}_{{\cal E}_{P}} d f 
= \delta \cdot (d f)$, so that 
$- \iota_{[{\cal E}_{P}, H^{\theta}_f]} d \theta 
= (\delta - \epsilon_{\theta}) d f 
= - \iota_{(\delta - \epsilon_{\theta}) H^{\theta}_{f}} d \theta$. 
By $(d \theta)^{\perp} = 0$, 
$[{\cal E}_{P}, H^{\theta}_f] 
= (\delta - \epsilon_{\theta}) H^{\theta}_{f}$, 
the first assertion. 
By the claim (B) for $\delta := \epsilon_{\theta}$, 
the second assertion follows. 
If $f_{i} \in {\cal O}^{\delta_{i}}(W)$, then 
${\cal E}_{P} \{ f_{1} \mid f_{2} \} 
= {\cal E}_{P} d \theta (H^{\theta}_{f_{1}}, H^{\theta}_{f_{2}}) 
= ({\cal L}_{{\cal E}_{P}} d \theta) (H^{\theta}_{f_{1}}, H^{\theta}_{f_{2}}) 
+ d \theta ({\cal L}_{{\cal E}_{P}} H^{\theta}_{f_{1}}, H^{\theta}_{f_{2}}) 
+ d \theta (H^{\theta}_{f_{1}}, {\cal L}_{{\cal E}_{P}} H^{\theta}_{f_{2}}) 
= (\epsilon_{\theta} + (\delta_{1} - \epsilon_{\theta}) + (\delta_{2} - \epsilon_{\theta})) 
d \theta (H^{\theta}_{f_{1}}, H^{\theta}_{f_{2}}) 
= (\delta_{1} + \delta_{2} - \epsilon_{\theta}) \{ f_{1} \mid f_{2} \}$, 
so that 
$\{ f_{1} \mid f_{2} \} 
\in {\cal O}^{\delta_{1} + \delta_{2} - \epsilon_{\theta}}(W)$. 
\\ \
(E) 
(a) 
Let $X \in a(W)$ be projectable: 
There exists a map $Y: p(W) \to T' M; \xi \mapsto Y |_{\xi}$ 
such that $p_{*\eta} X = Y |_{p(\eta)}$ for all $\eta \in W$. 
Since $p$ gives a holomorphic submersion from $W$ 
to the open subset $p(W)$ of $M$, 
each $\xi \in p(W)$ admits some open set $V$ of $p(W)$ 
and $\sigma \in \Gamma_{V}(W)$ such that $\xi \in V$ 
by the constant rank map theorem 
(e.g. \cite[p.86, Thm.I.15.5]{Hs1978}). 
Then $Y |_{\xi} = p_{*\sigma(\xi)} X$, 
which is holomorphic on $\xi \in M$, 
so that $Y \in a(M)$ and $p_{*} X = Y$. 
Assume that $X \in a(W, \theta)$. 
Take any $\sigma \in \Gamma_{\rm loc}(P)$ 
and $Y' \in \Gamma_{V_{\sigma}} (E_{\theta})$. 
Put 
$X' \in \Gamma_{V_{\sigma}} (D_{\theta})$ 
as 
$X' |_{\phi_{\sigma}(\xi, \lambda)} 
:= R_{\lambda*} (\sigma_{* \xi} Y')$ 
for $(\xi, \lambda) \in V_{\sigma} \times \mathbb{C}^{\times}$. 
Then $p_{*} X' = Y'$, 
$\iota_{[X, X']} \theta 
= [ {\cal L}_{X}, \iota_{X'} ] \theta 
= {\cal L}_{X} (\iota_{X'} \theta) - \iota_{X'} ({\cal L}_{X} \theta) 
= {\cal L}_{X} 0 - \iota_{X'} 0 
= 0$ 
and 
$[X, X'] \in \Gamma_{\rm loc}(D_{\theta})$, 
so that 
$[Y, Y'] = p_{*} [X, X'] \in \Gamma_{\rm loc}(E_{\theta})$. 
Hence, 
$Y \in a(M, E_{\theta})$, 
as required. 
\\ \
(b) 
Take any $\eta \in S$. 
Put 
$U := \hat{\mathbb{C}}^{\times}_{\eta, W}$. 
By the assumption on $S$ and Lemma 2.2 (B), 
$U$ is a connected open subset of $\mathbb{C}^{\times}$ 
such that $1 \in U$. 
Put $h: U \to T' W |_{\eta}; 
h(\lambda) := \hat{\eta}_{X}(\lambda) - \hat{\eta}_{X}(1)$, 
so that $h(1) = 0$. 
At $\lambda = 1$, 
$\left. \frac{d}{d \lambda} 
\hat{\eta}_{X} (\lambda) \right|_{\lambda = 1} 
= - {\cal L}_{{\cal E}_{P}} X |_{\eta} = 0$. 
At any $\lambda \in U$, 
\[
\lambda \frac{d}{d \lambda} h = \frac{d}{d \lambda_{1}} 
R_{\lambda*} R_{\lambda_{1} *} 
(X |_{R_{\lambda_{1}}^{-1} R_{\lambda}^{-1} \eta}) |_{\lambda_{1} = 1} 
= - R_{\lambda *} ({\cal L}_{{\cal E}_{P}} X) |_{R_{\lambda}^{-1} \eta}) 
= 0 = 0 \cdot h(\lambda). 
\]
Let $e_{1}, \ldots, e_{{\rm n}_{W}}$ 
be a $\mathbb{C}$-linear basis of 
$T' W |_{\eta}$ by which put 
$\sum_{i = 1}^{{\rm n}_{W}} h_{i} e_{i} := h$. 
By Lemma 2.2 (A) for $U$ with $(g, m, k) := (h_{i}, 1, 0)$, 
$h_{i}(\lambda) = h_{i}(1) \cdot \lambda^{0} = 0 \cdot 1 = 0$ 
for $\lambda \in U$ ($i \in \{ 1, \ldots, {\rm n}_{W} \}$), 
so that $h(\lambda) = 0$ ($\lambda \in U$). 
Hence, 
$R_{\lambda *} (X |_{R_{\lambda}^{-1} \eta}) = X |_{\eta}$, 
so that 
$p_{*} (X |_{R_{\lambda}^{-1} \eta}) 
= p_{*} (R_{\lambda *}^{-1} (X |_{\eta})) 
= p_{*} (X |_{\eta})$ ($\lambda \in U$). 
For any $\xi \in p(S)$ with some $\eta \in S$ such as $\xi = p(\eta)$, 
put $Y |_{\xi} := p_{*\eta} X |_{\eta} \in T' M$. 
Then $Y: p(S) \to T' M; \xi \mapsto Y |_{\xi}$ 
is well-defined as a map since any $\eta' \in W$ 
such as $p(\eta') = p(\eta)$ admits some $\lambda \in U$ 
such that $\eta' = R_{\lambda}^{-1} \eta$. 
Hence, $X$ is projectable, as required. 
\\ \
(c) 
For any $X \in a(W, \theta)$, 
${\cal L}_{{\cal E}_{P}} X = [{\cal E}_{P}, X] 
= 0$ by the claim (D). 
Because of the claim (b), 
$X$ is projectable. 
By the claim (a), 
$\widetilde{p}_{*}: a(W, \theta) \to a(M, E_{\theta}); 
X \mapsto p_{*} X$ 
is well-defined as a map, 
i.e., 
any $X, X' \in a(W, \theta)$ admit some 
$Y, Y' \in a(M, E_{\theta})$ 
such as $p_{*} X = Y$ and $p_{*} X' = Y'$. 
Then $p_{*} [X, X'] = [Y, Y']$. 
Since $\widetilde{p}_{*}$ is $\mathbb{C}$-linear, 
$\widetilde{p}_{*}$ is a complex Lie algebra homomorphism. 
\\ \
(d) 
{\it Claim 0. $\widetilde{\varpi}_{E}$ is injective:} 
In fact, let $Y \in a(M, E)$ satisfy that $\widetilde{\varpi}_{E} Y = 0$. 
For any $Y_{1} \in \Gamma_{\rm loc}(E)$, 
$\alpha_{E}(Y, Y_{1}) = \varpi_{E}([Y, Y_{1}]) 
= 0$, so that $Y = 0$ since $\alpha_{E}$ is non-degenerate. 
Since $\widetilde{\varpi}_{E}$ is $\mathbb{C}$-linear map, 
$\widetilde{\varpi}_{E}$ is injective, 
as required. 
\\ \
{\it Claim 1. 
$\widetilde{\varpi}_{E} \circ \widetilde{p}_{E*} 
\circ \widetilde{H}_{E}^{1} \circ \widetilde{\iota}_{L_{E}} 
= {\rm id}_{\Gamma_{M}(L_{E})}$, 
so that 
$\widetilde{\varpi}_{E}$ is surjective: 
} 
In fact, for any $s \in \Gamma_{M}(L_{E})$ 
and $\eta \in P_{L_{E}}$, 
$<\!\! \eta, 
\varpi_{E}(p_{E*\eta} H^{\theta_{E}}_{\iota_{L_{E}}(s)}) \!\!>_{L_{E}} 
= \theta_{E}(H^{\theta_{E}}_{\iota_{L_{E}}(s)}) |_{\eta} 
= \iota_{L_{E}}(s) |_{\eta}$ 
by the claim (B), 
which equals $<\!\! \eta, s(p(\eta)) \!\!>_{L_{E}}$, 
so that 
$\varpi_{E}(p_{E*\eta} H^{\theta_{E}}_{\iota_{L_{E}}(s)}) |_{p_{E}(\eta)} 
= s(p_{E}(\eta))$, 
which proves the first equation. 
By Lemmas 3.1 (A) and 4.1 (B), 
$\widetilde{\iota}_{L_{E}}$ 
and $\widetilde{H}_{\theta_{E}}^{1}$ 
are $\mathbb{C}$-linear isomorphisms, 
which proves the last assertion 
since ${\rm id}_{\Gamma_{M}(L_{E})}$ 
is bijective. 
\\ \
(e) By the claims (C), (D), (E) (d), 
the claim 1 of the step (d) and Lemma 3.1 (A), 
one has that 
$\widetilde{p}_{E*} = \widetilde{\varpi}_{E}^{-1} \circ {\rm id}_{\Gamma_{M}(L_{E})} 
\circ \widetilde{\iota}_{L_{E}}^{-1} \circ (\widetilde{H}^{1}_{\theta_{E}})^{-1} 
= {\cal P}_{\theta_{E}}^{-1}$, 
which is a $\mathbb{C}$-linear isomorphism. 
By the claim (c) for $(W, \theta) := (P_{E}, \theta_{E})$, 
$\widetilde{p}_{E*}$ is a complex Lie algebra isomorphism. 
~\qed
\smallskip \\ \
Let 
$(M, E)$ 
be a c-manifold with an $E$-related local symplectification 
$(W, \theta)$ of degree $\epsilon_{\theta}$ 
in a holomorphic principlal $\mathbb{C}^{\times}$-bundle $P$ 
on $M$ with the projection $p: P \to M$. 
Let 
$\kappa_{\theta}: 
W \to P_{E} := P_{L_{E}}; \eta \mapsto \kappa_{\theta} (\eta)$ 
be a map such that 
$p_{E}(\kappa_{\theta}(\eta)) = p(\eta)$ 
and 
$<\!\! \kappa_{\theta} (\eta), \varpi_{E}(p_{*\eta} X) \!\!>_{L_{E}} 
= \theta (X)$ 
for all $\eta \in W$ and $X \in {T'} P |_{\eta}$, 
which is said to be {\it the canonical map} of $(W, \theta)$. 
Especially, 
$\kappa_{\theta_{E}} = {\rm id}_{P_{E}}$ 
if 
$(W, \theta) = (P_{E}, \theta_{E})$. 
Then one has the following result. 
\smallskip \\ \
\textsc{Lemma 4.2.} 
(A) 
{\it 
$\kappa_{\theta}$ is well-defined as a holomorphic immersion 
satisfying that 
$\kappa_{\theta}^{*} \theta_{E} = \epsilon_{\theta} \cdot \theta$, 
so that $\kappa_{\theta}$ is locally biholomorphic 
and that 
$\kappa_{\theta*} {\cal E}_{P} = \epsilon_{\theta} \cdot {\cal E}_{P_{E}}$. 
}
\\ \
(B) 
{\it 
For any map 
$Y: p(W) \to T' M; \xi \mapsto Y |_{\xi}$ 
such as $Y |_{\xi} \in T' M |_{\xi}$, put 
$j_{\theta}(Y) := \iota_{L_{E}} (\varpi_{E}(Y)) \circ \kappa_{\theta}: 
W \to \mathbb{C}$. 
Put 
$\widetilde{j}_{\theta}: 
a(p(W), E) \to {\cal O}^{\epsilon_{\theta}}(W); Y \mapsto j_{\theta}(Y)$. 
Then 
$\widetilde{j}_{\theta}$ 
is well-defined as a $\mathbb{C}$-linear map, 
by which the following $\mathbb{C}$-linear map is well-defined: 
\[
{\cal P}_{\theta} 
:= \widetilde{H}^{\epsilon_{\theta}}_{\theta} \circ \widetilde{j}_{\theta}: 
a(p(W), E) \to a(W, \theta); 
Y \mapsto {\cal P}_{\theta}(Y) := H^{\theta}_{j_{\theta}(Y).}  
\] 
Let $X \in a(W, \theta)$ be projectable such that 
$p_{*} X = Y$ for some vector field $Y$ on $p(W)$. 
Then $Y \in a(p(W), E)$ and 
$\widetilde{j}_{\theta}(Y) = \theta (X) 
\in {\cal O}^{\epsilon_{\theta}}(W)$. 
}
\smallskip \\ \
{\it Proof}. 
(A) 
For $\eta \in W$, 
take some 
$X \in T_{\eta}' P \backslash D_{\theta} 
= T_{\eta'} P \backslash p_{*'}^{-1} E$. 
By $\theta(X) \ne 0$ and $\varpi_{E}(p_{*'} X) \ne 0$, 
there exists unique $\eta_{X} \in P_{E}$ such that 
$<\!\! \eta_{X}, \varpi_{E}(p_{*'} X) \!\!>_{L_{E}} 
= \theta (X)$. 
For any $X' \in T_{\eta}' P \backslash D_{\theta}$, 
put $c := \theta(X')/ \theta(X) \in \mathbb{C}$. 
Then 
$X' - c \cdot X \in D_{\theta} = p_{*'}^{-1} E$, 
so that 
$\varpi_{E}(p_{*'} X') = c \cdot \varpi_{E}(p_{*'} X)$. 
Hence, one has that 
$<\!\! \eta_{X}, \varpi_{E}(p_{*'} X') \!\!>_{L_{E}} 
= c \cdot \theta(X) = \theta (X') 
= <\!\! \eta_{X'}, \varpi_{E}(p_{*'} X') \!\!>_{L_{E}}$, 
so that $\eta_{X} = \eta_{X'}$, 
which does not depend on the choice of such $X$. 
Hence, a map 
$\kappa_{\theta}: W \to P_{E}; \eta \mapsto \eta_{X}$ 
is well-defined such that $p_{E} \circ \kappa_{\theta} = p |_{W}$. 
If 
$X \in \Gamma_{\rm loc}({T'} W \backslash D_{\theta})$, 
then 
$\theta(X |_{\eta})$ is a holomorphic function of $\eta$, 
so does $\eta_{X} (= \kappa_{\theta}(\eta))$ 
given as 
$<\!\! \eta_{X}, \varpi_{E}(p_{*'\eta} X) \!\!>_{L_{E}} 
= \theta (X |_{\eta})$. 
Hence, $\kappa_{\theta}$ is a holomorphic map. 
For $X \in {T'} W |_{\eta}$, 
$\theta_{E}(\kappa_{\theta *} X) 
= <\!\! \kappa_{\theta} (\eta), \varpi_{E}(p_{*'} X) \!\!>_{L_{E}} 
= \theta(X)$, i.e., 
$\kappa_{\theta}^{*} \theta_{E} = \theta$, 
so that $\kappa_{\theta}^{*} d \theta_{E} = d \theta$. 
For $X \in {T'} W |_{\eta}$, 
if $\kappa_{\theta*} X = 0$, 
then 
$0 = d \theta_{E}(\kappa_{\theta*} X, \kappa_{\theta*} Y) = d \theta (X, Y)$ 
for all $Y \in {T'} W |_{\eta}$, 
so that $X = 0$ because of $(d \theta)^{\perp} = 0$. 
Hence, $\kappa_{\theta*'\eta}$ is injective at $\eta \in W$, 
which is also surjective by 
${\rm n}_{W} = {\rm n}_{P} = {\rm n}_{P_{E}}$. 
Then $\kappa_{\theta}$ is locally biholomorphic. 
\\ \
Let $\sigma \in \Gamma_{\rm loc}(P)$ be such that 
there exist a system $(z_{1}, \ldots, z_{{\rm n}_{M}}): 
V_{\sigma} \to \mathbb{C}_{{\rm n}_{M}}$ of complex coordinates 
with some 
$Y_{0} \in \Gamma_{V_{\sigma}}({T'} M \backslash E)$. 
Put 
$X_{0} |_{\sigma(\xi)} := \sigma_{*\xi} Y_{0}$. 
Then 
$p_{*'} (X_{0} |_{\sigma(\xi)}) = Y_{0} |_{\xi} \not\in E$, 
so that 
$X_{0} |_{\sigma(\xi)} \not\in p_{*'}^{-1} E = D_{\theta}$. 
Put $W_{\sigma} := W \cap p^{-1}(V_{\sigma})$, on which 
$\widetilde{X}_{0} \in \Gamma_{W_{\sigma}}({T'} P \backslash D_{\theta})$ 
is defined as 
$\widetilde{X}_{0} |_{\phi_{\sigma}(\xi_{1}, \lambda_{1})} 
:= R_{\lambda_{1}*} (X_{0} |_{\sigma(\xi_{1})})$ 
for $(\xi_{1}, \lambda_{1}) \in V_{\sigma} 
\times \breve{\mathbb{C}}^{\times}_{\sigma(\xi_{1}), W}$, 
so that 
$R_{\lambda*\eta_{1}} \widetilde{X}_{0} = \widetilde{X}_{0}$ 
for $\eta_{1} \in W_{\sigma}$ and 
$\lambda \in \breve{\mathbb{C}}^{\times}_{\eta_{1}, W_{\sigma}}$ 
and that $\widetilde{X}_{0} \not\in D_{\theta}$ at all points. 
Fix any $\eta_{0} \in W_{\sigma}$ 
and put 
$\eta_{0}' := \kappa_{\theta}(\eta_{0}) \in \kappa_{\theta}(W_{\sigma})$. 
Since $\kappa_{\theta}$ is locally biholomorphic, 
there exist an open neighbourhod $U_{0}'$ of $\eta_{0}'$ in $P_{E}$ 
and an open neighbourhood $U_{0}$ of $\eta_{0}$ in $W_{\sigma}$ 
such that the restriction 
$\kappa_{\theta} |_{U_{0} \to U_{0}'}$ 
is biholomorphic. 
Put 
$f_{0}: U_{0}' = \kappa_{\theta}(U_{0}) \to \mathbb{C}^{\times}; 
\eta' \mapsto \theta(\widetilde{X}_{0} |_{\kappa_{\theta}^{-1}(\eta')}) 
= <\!\! \eta', 
\varpi_{E}(p_{*'} \widetilde{X}_{0}) \!\!>_{L_{E}}$, 
so that $f_{0} (R_{\lambda} \eta') = \lambda \cdot f_{0}(\eta')$ 
for any 
$(\eta', \lambda) \in U_{0}' \times 
\breve{\mathbb{C}}^{\times}_{\eta', U_{0}'.}$ 
Hence, 
$f_{0} \in {\cal O}^{1}(U_{0}')$. 
Take any 
$\eta_{1} = \phi_{\sigma}(\xi_{1}, \lambda_{1}) \in U_{0}$. 
Then 
\begin{eqnarray*}
(\kappa_{\theta*\eta_{1}} {\cal E}_{P}) f_{0} 
&=& {\cal E}_{P}|_{\eta_{1}} (f_{0} \circ \kappa_{\theta}) 
= \left. 
\frac{d}{d \lambda} f_{0}(\kappa_{\theta}(R_{\lambda} \eta_{1})) 
\right|_{\lambda = 1} 
= 
\left. 
\frac{d}{d \lambda} \theta(\widetilde{X}_{0} |_{R_{\lambda} \eta_{1}}) 
\right|_{\lambda = 1} 
\\
&=& \left. 
\frac{d}{d \lambda} \theta(R_{\lambda*\eta_{1}} \widetilde{X}_{0}) 
\right|_{\lambda = 1} 
= ({\cal L}_{{\cal E}_{P}} \theta) (\widetilde{X}_{0} |_{\eta_{1}}) 
= \epsilon_{\theta} \cdot \theta (\widetilde{X}_{0} |_{\eta_{1}}) 
\\
&=& \epsilon_{\theta} \cdot f_{0}(\kappa_{\theta} (\eta_{1})) 
= \epsilon_{\theta} \cdot ({\cal E}_{P_{E}} f_{0}) (\kappa_{\theta}(\eta_{1})) 
= \epsilon_{\theta} \cdot {\cal E}_{P_{E}} |_{\kappa_{\theta}(\eta_{1})} f_{0}. 
\end{eqnarray*}
For $i \in \{ 1, \ldots, {\rm n}_{M} \}$, 
put 
$f_{i}(\kappa_{\theta}(\eta_{1})) 
:= z_{i}(\xi_{1}) \cdot f_{0}(\kappa_{\theta}(\eta_{1}))$. 
Then 
\[
(\kappa_{\theta*\eta_{1}} {\cal E}_{P}) f_{i} 
= z_{i}(\xi_{1}) \cdot (\kappa_{\theta*\eta_{1}} {\cal E}_{P}) f_{0} 
= z_{i}(\xi_{1}) \cdot 
\epsilon_{\theta} \cdot {\cal E}_{P_{E}} |_{\kappa_{\theta}(\eta_{1})} f_{0}
= \epsilon_{\theta} \cdot {\cal E}_{P_{E}} |_{\kappa_{\theta}(\eta_{1})} f_{i}, 
\]
because of 
${\cal E}_{P_{E}} f_{i} = z_{i} \cdot {\cal E}_{P_{E}} f_{0}$. 
Put 
$F := ( f_{0}, f_{1}, \ldots, f_{{\rm n}_{M}}): 
U_{0}' \to \mathbb{C}_{{\rm n}_{W}}$. 
Then $d F = (d f_{0}, f_{0} d z_{1} + z_{1} d f_{0}, \ldots, 
f_{0} d z_{{\rm n}_{M}} + z_{{\rm n}_{M}} d f_{{\rm n}_{M}})$, 
so that 
\begin{eqnarray*}
& &\{ X \in T' U_{0}'; (d F) X = 0 \} 
\\
&=& \{ X \in T' U_{0}; d f_{0} X = d z_{i}(p_{E*} X) = 0~(i = 1, \ldots, {\rm n}_{M}) \} 
\\
&=& \{ X \in T' U_{0}; d f_{0} X = 0,~p_{E*} X = 0
~(i = 1, \ldots, {\rm n}_{M}) \}
\\ 
&=& 
\{ c \cdot {\cal E}_{P} |_{\eta'}; 
c \in \mathbb{C}, \eta' \in U_{0}', 
d f_{0}(c \cdot {\cal E}_{P} |_{\eta'}) = 0 \}
\\
&=& 
\{ c \cdot {\cal E}_{P} |_{\eta'}; 
c \in \mathbb{C}, \eta' \in U_{0}', 
c \cdot f_{0}(\eta') = 0 \}
\\
&=& 0. 
\end{eqnarray*}
Hence, $F$ gives a system of complex coordinates of 
the complex manifold $U_{0}'$ of complex dimension ${\rm n}_{M}$, 
so that any locally defined holomorphic function $f$ 
on $U_{0}'$ is described as some power series of 
$f_{i}$'s ($i \in \{ 0, \ldots, N \}$). 
Then 
$(\kappa_{\theta*\eta_{1}} {\cal E}_{P}) f 
= \epsilon_{\theta} \cdot 
{\cal E}_{P_{E}} |_{\kappa_{\theta}(\eta_{1})} f$. 
Hence, 
$\kappa_{\theta*\eta_{1}} {\cal E}_{P} 
= \epsilon_{\theta} \cdot {\cal E}_{P_{E}} 
|_{\kappa_{\theta}(\eta_{1})}$.
\\ \
(B)
By the claim (A), $j_{\theta}(Y)$ is well-defined as a map for any $Y$. 
Assume that $Y \in a(p(W), E)$. 
By Lemmas 4.1 (E) (d) and 3.1 (A), one has that 
$\iota_{L_{E}}(\varpi_{E}(Y)) \in {\cal O}^{1}(P_{E} |_{p_{E}(W)})$. 
Then $\widetilde{j}_{\theta} (Y) \in {\cal O}^{\epsilon_{\theta}} (W)$, 
because of 
\begin{eqnarray*}
& &{\cal E}_{P} (\widetilde{j}_{\theta}(Y)) 
= {\cal E}_{P} (\iota_{L_{E}}(\varpi_{E}(Y)) \circ \kappa_{\theta}) 
= ((\kappa_{\theta *} {\cal E}_{P}) \iota_{L_{E}}(\varpi_{E}(Y))) 
\circ \kappa_{\theta} 
\\
&=& (\epsilon_{\theta} \cdot {\cal E}_{P_{E}} \iota_{L_{E}}(\varpi_{E}(Y))) 
\circ \kappa_{\theta} 
= \epsilon_{\theta} \cdot \iota_{L_{E}}(\varpi_{E}(Y)) \circ \kappa_{\theta} 
= \epsilon_{\theta} \cdot \widetilde{j}_{\theta}(Y) 
\end{eqnarray*}
by (A). 
Hence, 
$\widetilde{j}_{\theta}$ is well-defined as a $\mathbb{C}$-linear map. 
By Lemma 4.1 (C), 
${\cal P}_{\theta} 
= \widetilde{H}^{\epsilon_{\theta}} \circ \widetilde{j}_{\theta}$ 
is well-defined as a $\mathbb{C}$-linear map. 
Let 
$X \in a(W, \theta)$ be projectable. 
Then 
$Y := p_{*} X \in a(p(W), E)$ by Lemma 4.1 (E) (a). 
For $\eta \in W$, 
\[
\widetilde{j}_{\theta} (Y) |_{\eta} 
= <\!\! \kappa_{\theta} (\eta),~\varpi_{E}(Y) \!\!>_{L_{E}} 
= <\!\! \kappa_{\theta} (\eta),~\varpi_{E}(p_{*\eta} X) \!\!>_{L_{E}}
= \theta (X) |_{\eta}, 
\]
as required. 
~\qed 
\smallskip \\ \
\textsc{Proposition 4.3.} 
{\it 
Let $(W, \theta)$ be an $E$-related local sysmplectification of degree 
$\epsilon_{\theta} \in \mathbb{C}^{\times}$ 
in a holomorphic principal $\mathbb{C}^{\times}$-bundle $P$ 
on a c-manifold $(M, E)$ with the projection 
$p: P \to M$ and the right action $R_{\lambda}$ 
of $\lambda \in \mathbb{C}^{\times}$ 
such that $W |_{\eta}$ is connected for all $\eta \in W$. 
Then one has the following two results.
} 
\\ \
(A) 
{\it 
$\widetilde{p}_{*} \circ {\cal P}_{\theta} = {\rm id}_{a(p(W), E)}$ 
and 
${\cal P}_{\theta} \circ \widetilde{p}_{*} = {\rm id}_{a(W, \theta)}$, 
so that $\widetilde{p}_{*}$ and ${\cal P}_{\theta}$ 
are isomorphisms of complex Lie algebras. 
Moreover, 
$\widetilde{j}_{\theta}$ 
is a $\mathbb{C}$-linear isomorphism onto 
${\cal O}^{\epsilon_{\theta}} (W)$ 
spliting $\widetilde{\theta}^{\epsilon_{\theta}}$ 
as  
$\widetilde{j}_{\theta} \circ \widetilde{p}_{*} 
= \widetilde{\theta}^{\epsilon_{\theta}}$ 
on $a(W, \theta)$. 
}
\\ \
(B) 
{\rm 
Assume that $W = P$ and $\epsilon_{\theta} = \ell$ 
for some 
$\ell \in \mathbb{C}^{\times} \cap \mathbb{Z}$. 
Then $\kappa_{\theta} \circ R_{\lambda} = \lambda^{\ell} \cdot \kappa_{\theta}$ 
$(\lambda \in \mathbb{C}^{\times})$ 
and that $\kappa_{\theta}$ is a covering map with $| \ell |$-sheets. 
}
\smallskip \\ \
{\it Proof.} 
(A) 
Take any $\eta \in W$ and $Y \in a(p(W), E)$. 
By Lemma 4.2 (B), ${\cal P}_{\theta} (Y) \in a(W, \theta)$. 
By Lemma 4.1 (E) (c), 
$\widetilde{p}_{*} ({\cal P}_{\theta} (Y)) \in a(p(W), E)$ 
since $W |_{\eta}$ is connected for all $\eta \in W$. 
Put 
$Y_{0} := \widetilde{p}_{*} ({\cal P}_{\theta} (Y)) - Y 
\in a(p(W), E_{\theta}) = a(p(W), E)$. 
By the first assertion of Lemma 4.1 (B) for $\delta = \epsilon_{\theta}$, 
\[
<\!\! \kappa_{\theta} (\eta),~
\varpi_{E}(\widetilde{p}_{*} ({\cal P}_{\theta} (Y))) \!\!>_{L_{E}} 
= \theta({\cal P}_{\theta}(Y)) |_{\eta} 
= \theta(H^{\theta}_{j_{\theta}(Y)}) |_{\eta}
= j_{\theta}(Y) |_{\eta} 
\]
which equals 
$<\!\! \kappa_{\theta} (\eta),~\varpi_{E}(Y) \!\!>_{L_{E}}$, 
so that $\varpi_{E}(Y_{0}) = 0$ on $p(W)$. 
By Lemma 4.1 (E) (d), $Y_{0} = 0$, i.e., 
$\widetilde{p}_{*} ({\cal P}_{\theta} (Y)) = Y$. 
Hence, $\widetilde{p}_{*} \circ {\cal P}_{\theta} = {\rm id}_{a(p(W), E)}$. 
\\ \
Conversely, take any $X \in a(W, \theta)$. 
By Lemma 4.1 (E) (c), 
$Y := \widetilde{p}_{*} X \in a(p(W), E)$ 
since $W |_{\eta}$ is connected for all $\eta \in W$. 
By Lemma 4.2 (B), 
$\widetilde{j}_{\theta}(Y) = \theta (X)$. 
By Lemma 4.1 (B), 
$X = H^{\theta}_{\theta (X)} = H^{\theta}_{\widetilde{j}_{\theta} (Y)} 
= {\cal P}_{\theta} (Y) = {\cal P}_{\theta} (\widetilde{p}_{*} X)$. 
Hence, ${\cal P}_{\theta} \circ \widetilde{p}_{*} = {\rm id}_{a(W, \theta)}$. 
Then 
$\widetilde{p}_{*}$ and ${\cal P}_{\theta}$ 
are $\mathbb{C}$-linear isomorphisms, 
which are isomorphisms of complex Lie algebras 
by Lemma 4.1 (E) (c). 
By Lemmas 4.1 (C), 
$\widetilde{j}_{\theta} 
= (\widetilde{H}^{\epsilon_{\theta}})^{-1} \circ {\cal P}_{\theta}: 
a(p(W), E) \to {\cal O}^{\epsilon_{\theta}}(W)$, 
which is then a $\mathbb{C}$-linear isomorphism. 
\\ \
(B) 
Put 
$\hat{\mu}_{\ell}: 
\mathbb{C}^{\times} \to \mathbb{C}^{\times}; 
\lambda \mapsto \lambda^{\ell}$. 
For $\sigma \in \Gamma_{\rm loc}(P)$, 
$(\xi, \lambda) \in V_{\sigma} \times \mathbb{C}^{\times}$, 
$\eta = R_{\lambda} \sigma(\xi)$ 
and $X \in T' P |_{\eta}$, 
\begin{eqnarray*}
<\!\! \kappa_{\theta} (\eta), \varpi_{E}(p_{*'} X) \!\!>_{L_{E}} 
&=& \theta (X) 
= \theta (R_{\lambda*} (R_{\lambda*}^{-1} X)) 
= \lambda^{\ell} \cdot 
\theta |_{\sigma (\xi)} (R_{\lambda*}^{-1} X) 
\\
&=& < \!\! \lambda^{\ell} \cdot \kappa_{\theta} (\sigma (\xi)), 
\varpi_{E}(p_{*'}X) \!\!>_{L_{E}}
\end{eqnarray*} 
because of Lemma 4.0 (A) (e). 
Then 
$\kappa_{\theta} \circ R_{\lambda} 
= \lambda^{\ell} \cdot \kappa_{\theta}$ 
($\lambda \in \mathbb{C}^{\times}$). 
By Lemma 2.3 (A), 
$\kappa_{\theta}$ is a covering map with $|\ell|$-sheets. 
~\qed 
\smallskip \\ \
A complex Lie group $G$ acting on a complex manifold $M$ 
is said to be a {\it holomorphic Lie transformation group} 
when the action map, 
$G \times M \to M; (g, \xi) \mapsto g \xi$, 
is holomorphic. 
Then one has the following result. 
\smallskip \\ \
\textsc{Lemma 4.4}. 
(A) (cf. \cite[p.74, Thm.2.4']{OK1977}) 
{\it 
Let 
$(P, \theta)$ 
be a symplectification of degree 
$\ell \in \mathbb{C}^{\times} \cap \mathbb{Z}$ 
on $M$. 
Put 
${\rm Aut}(P, \theta) 
:= \{ \beta \in {\rm Aut}(P){;}~\beta^{*} \theta = \theta \}$ 
and 
${\rm Aut}(P, d \theta) 
:= \{ \beta \in {\rm Aut}(P){;}~\beta^{*} d \theta = d \theta \}$. 
Then 
\begin{eqnarray*}
{\rm Aut}(P, \theta) 
&=& \{ \beta \in {\rm Aut}(P, d \theta); 
\beta_{*} {\cal E}_{P} = {\cal E}_{P} \} 
\\
&=& \{ \beta \in {\rm Aut}(P, d \theta); 
R_{\lambda} \circ \beta 
= \beta \circ R_{\lambda}~(\lambda \in \mathbb{C}^{\times}) \}. 
\end{eqnarray*}
Let $(W, \theta)$ 
be a local symplectification of degree $\epsilon_{\theta} \ne 0$ 
in $P$ on $M$, 
and put 
$a(W, d \theta) := \{ X \in a(W){;~} {\cal L}_{X} d \theta = 0 \}$. 
Then 
\[
a(W, \theta) = \{ X \in a(W, d \theta){;~} [{\cal E}_{P}, X] = 0 \}. 
\]
}
\
(B) (cf. \cite[Thm.1.5]{NT1987}) 
{\it 
Let $(M, E)$ be a c-manifold with 
${\rm Aut}(M, E) = \{ \alpha \in {\rm Aut}(M); \alpha_{*} E = E \}$. 
The prolongation of $\alpha \in {\rm Aut}(M, E)$ 
is uniquely determied as an element 
${\cal P}(\alpha) \in {\rm Aut}(P_{E})$ 
such that 
${\cal P}(\alpha) \eta \in P_{E} |_{\alpha (p_{E} \eta)}$ 
and 
\[
<\!\! {\cal P}(\alpha) \eta, 
\varpi_{E}(X) \!\!>_{L_{E}} 
:= <\!\! \eta, 
\varpi_{E}(\alpha_{*}^{-1} X) \!\!>_{L_{E}} 
\]
for 
$\eta \in P_{E}$ and $X \in {T'} M |_{\alpha(p_{E}(\eta))}$, 
so that 
${\rm Aut}(M, E) \to {\rm Aut}(P_{E}); 
\alpha \mapsto {\cal P}(\alpha)$ 
is well-defined as a group homomorphism such that 
$p_{E} \circ {\cal P} (\alpha) = \alpha \circ p_{E}$. 
Let 
$(P_{E}, \theta_{E})$ 
be the canonical symplectification 
in the sense of Lemma 4.0 (B). 
Then 
${\cal P}(\alpha) \in {\rm Aut}(P_{E}, \theta_{E})$ 
for $\alpha \in {\rm Aut}(M, E)$. 
Moreover, 
\[
{\cal P}: 
{\rm Aut}(M, E) \to 
{\rm Aut}(P_{E}, \theta_{E}); 
\alpha \mapsto {\cal P} (\alpha) 
\]
is well-defined as a group isomorphism. 
Assume that 
${\rm Aut}(M, E)$ 
has a structure of a complex Lie group and that 
${\rm Aut}(M, E) \times M \to M; (\alpha, \xi) \mapsto \alpha \xi$ 
is holomorphic. 
Then ${\rm Aut}(P_{E}, \theta_{E})$ 
has a structure of a complex Lie group 
such that ${\cal P}$ is a complex Lie group isomorphism, 
and that the inverse map 
${\cal P}^{-1}: {\rm Aut}(P_{E}, \theta_{E}) \to {\rm Aut}(M, E); 
\beta \mapsto \beta_{M}$ 
satisfies that $\beta_{M} \circ p_{E} = p_{E} \circ \beta$. 
Moreover, 
${\rm Aut}(P_{E}, \theta_{E}) \times P_{E} \to P_{E}; 
(\beta, \eta) \mapsto \beta \eta$ 
is holomorphic. 
}
\smallskip \\ \
{\it Proof}. 
(A0) 
Take any $\beta \in {\rm Aut}(P, \theta)$. 
Then 
$\beta^{*} d \theta 
= d (\beta^{*} \theta) 
= d \theta$, 
so that 
$d \theta(Y, \beta_{*} \natural \theta) 
= d \theta(\beta_{*}^{-1} Y, \natural \theta) 
= \theta(\beta_{*}^{-1} Y) 
= \theta(Y) = d \theta(Y, \natural \theta)$ 
for $Y \in T' P$. 
Hence, 
$\beta_{*} (\natural \theta) 
= \natural \theta$, 
so that 
$\beta_{*} {\cal E}_{P} = {\cal E}_{P}$ 
by Lemma 4.1 (A). 
For 
$(\lambda, \eta) \in \mathbb{C}^{\times} \times P_{E}$, 
one has that 
\begin{eqnarray*}
& &\lambda \frac{d}{d \lambda} 
R_{\lambda}^{-1} 
(\beta (R_{\lambda} \eta)) 
= \left. 
\frac{d}{d \lambda_{1}} 
(R_{\lambda \lambda_{1}}^{-1} 
(\beta (R_{\lambda \lambda_{1}} \eta))) 
\right|_{\lambda_{1} = 1} 
\\ 
&=&
- {\cal E}_{P} |_{R_{\lambda}^{-1} (\beta (R_{\lambda} \eta))} 
+ R_{\lambda*}^{-1} \beta_{*} ({\cal E}_{P} |_{R_{\lambda} \eta}) 
\\ 
&=& - {\cal E}_{P} |_{R_{\lambda}^{-1} (\beta (R_{\lambda} \eta))} 
+ R_{\lambda*}^{-1} {\cal E}_{P} |_{\beta (R_{\lambda} \eta)} 
= 0,
\end{eqnarray*} 
so that 
$\frac{d}{d \lambda} 
R_{\lambda}^{-1} 
(\beta (R_{\lambda} \eta)) 
= 0$ 
and that 
$R_{\lambda}^{-1} \circ \beta \circ R_{\lambda} =\beta$. 
Conversely, take any $\beta \in {\rm Aut}(P, d \theta)$ 
such that 
$R_{\lambda} \circ \beta 
= \beta \circ R_{\lambda} 
(\lambda \in \mathbb{C}^{\times})$, 
so that 
$\beta_{*} {\cal E}_{P} = {\cal E}_{P}$. 
By Lemma 4.1 (A), 
$\beta_{*} \natural \theta = \natural \theta$, 
so that 
$\theta(\beta_{*}^{-1} Y) 
= d \theta(\beta_{*}^{-1} Y, \natural \theta) 
= d \theta(Y, \beta_{*} \natural \theta) 
= d \theta(Y, \natural \theta) 
= \theta(Y)$. 
\\ \
(A1) $a(W, \theta) \subseteq \{ X \in a(W, d \theta){;~} 
[{\cal E}_{P}, X] = 0 \}$, because of Lemma 4.1 (D). 
Conversely, 
assume that $X \in a(W, d \theta)$ and $[{\cal E}_{P}, X] = 0$. 
By Lemma 4.1 (A), 
$\iota_{{\cal E}_{P}} d \theta 
= - \epsilon_{\theta} \cdot \iota_{\natural \theta} d \theta 
= \epsilon_{\theta} \cdot \theta$. 
Then 
$\epsilon_{\theta} \cdot {\cal L}_{X} \theta 
= {\cal L}_{X} (\iota_{{\cal E}_{P}} d \theta) 
= ([{\cal L}_{X}, \iota_{{\cal E}_{P}}] 
+ \iota_{{\cal E}_{P}} \circ {\cal L}_{X}) d \theta 
= \iota_{[X, {\cal E}_{P}]} d \theta + \iota_{{\cal E}_{P}} 0 
= \iota_{0} d \theta = 0$, 
so that $X \in a(W, \theta)$, 
as required. 
\\ \
(B0) 
For any $\alpha \in {\rm Aut}(M, E)$, 
${\cal P}(\alpha) \in {\rm Aut}(P_{E}, \theta_{E})$ 
because of 
\[
\theta_{E}({\cal P}(\alpha)_{*\eta} Y) 
= <\!\!{\cal P}(\alpha)\eta, 
\varpi_{E*} (\alpha_{*} (p_{E*} Y))\!\!>_{L_{E}} 
= <\!\!\eta, \varpi_{E*}(p_{E*} Y) \!\!>_{L_{E}} 
= \theta_{E}(Y)
\]
for all $\eta \in P_{E}$ and $Y \in {T'} P_{E} |_{\eta}$. 
By definition of ${\cal P}$, 
$\alpha \circ p_{E} = p_{E} \circ {\cal P}(\alpha)$ 
Hence, 
${\cal P}$ 
is injective. 
\\ \ 
(B1) 
Assume that 
$\beta \in {\rm Aut}(P_{E}, \theta_{E})$. 
By the claim (A), 
$\beta \circ R_{\lambda} = R_{\lambda} \circ \beta$ 
for all $\lambda \in \mathbb{C}^{\times}$. 
Put $\beta_{M} (p_{E} \eta) := p_{E} (\beta \eta)$ 
for $\eta \in P_{E}$. 
If $\eta' \in P_{E}$ and $p_{E} \eta' = p_{E} \eta$, 
then $\eta' = R_{\lambda} \eta$ for some $\lambda \in \mathbb{C}^{\times}$, 
so that 
$\beta_{M} (p_{E} \eta') 
= p_{E} (\beta (R_{\lambda} \eta)) 
= p_{E} (R_{\lambda} (\beta \eta)) 
= p_{E} (\beta \eta)$. 
Hence, 
$\beta_{M}: M \to M; p_{E} \eta \mapsto \beta_{M} (p_{E} \eta)$ 
is well-defined as a map such that 
$\beta_{M} \circ p_{E} = p_{E} \circ \beta$, 
which is holomorphic with respect to 
the local trivializations of $P$. 
By 
$\beta_{M*} E 
= \beta_{M*} (p_{E*} D_{\theta_{E}}) 
= p_{E*} (\beta_{*} D_{\theta_{E}}) 
= p_{E*} D_{\theta_{E}} = E$, 
$\beta_{M} \in {\rm Aut}(M, E)$. 
For $\eta \in P_{E}$ and $X \in {T'} P_{E} |_{\eta}$, 
put 
$Y := \beta_{M*} (p_{E*} X) \in {T'} M |_{\beta_{M}(p_{E}(\eta))}$. 
Then 
\begin{eqnarray*}
<\!\! \beta (\eta), \varpi_{E} (Y) \!\!>_{L_{E}} 
&=& <\!\! \beta (\eta), \varpi_{E} (\beta_{M*}(p_{E*} X)) \!\!>_{L_{E}}
\\
&=& <\!\! \beta (\eta), \varpi_{E} (p_{E*} (\beta_{*} X)) \!\!>_{L_{E}}
\\
&=& \theta_{E} (\beta_{*} X) 
= \theta_{E} (X) 
= <\!\! \eta, \varpi_{E} (p_{E*} X) \!\!>_{L_{E}}
\\ 
&=& <\!\! \eta, \varpi_{E}(\beta_{M*}^{-1} Y)\!\!>_{L_{E}} 
= <\!\! {\cal P}(\beta_{M}) \eta,~\varpi_{E}(Y)\!\!>_{L_{E}}, 
\end{eqnarray*}
so that 
$\beta = {\cal P}(\beta_{M})$. 
Hence, 
${\cal P}$ is surjective. 
\\ \
(B2) 
${\rm Aut}(M, E)$ is a complex Lie group 
by assumption. 
${\cal P}$ is a group isomorphism by the steps (B0, B1). 
Then ${\rm Aut}(P_{E}, \theta_{E})$ 
has a structure of complex Lie group 
such that ${\cal P}$ is a complex Lie group isomorphism. 
The right hand side of the equation, 
$<\!\! {\cal P}(\alpha) \eta, 
\varpi_{E}(X) \!\!>_{L_{E}} 
:= <\!\! \eta, 
\varpi_{E}(\alpha_{*}^{-1} X) \!\!>_{L_{E}}$ 
($\eta \in P_{E}, X \in \Gamma_{\rm loc} ({T'} M)$), 
is holomorphic on $(\alpha, \eta) \in {\rm Aut}(M, E) \times P_{E}$, 
since the action map, 
${\rm Aut}(M, E) \times M \to M; (\alpha, \xi) \mapsto \alpha \xi$, 
is holomorphic, 
so that the left hand side is holomorphic on 
$(\alpha, \eta) \in {\rm Aut}(M, E) \times P_{E}$. 
Hence, 
the action map: 
${\rm Aut}(M, E) \times P_{E} \to P_{E}; 
(\alpha, \eta) \mapsto {\cal P}(\alpha) \eta$, 
is holomorphic. 
Since ${\cal P}$ is biholomorphic, the action map, 
${\rm Aut}(P_{E}, \theta_{E}) \times P_{E} \to P_{E}; 
(\beta, \eta) \mapsto \beta \eta 
= {\cal P}({\cal P}^{-1}(\beta)) \eta$ 
is holomorphic.~\qed 
\medskip \\ \
\textsc{Lemma 4.5} 
(cf. \cite[(3.2)]{Bwm1961}). 
{\it 
Let $(M, E)$ be a connected compact c-manifold of ${\rm n}_{M} > 0$. 
Then 
${\rm Aut}(P_{E}, \theta_{E})$ 
is a holomorphic Lie transformation group on $P_{E}$ 
such that $a(P_{E}, \theta_{E})$ is isomorphic to 
the complex Lie algebra of all left invariant vector field on 
${\rm Aut}(P_{E}, \theta_{E})$ 
such that 
${\rm n}_{{\rm Aut}(P_{E}, \theta_{E})} 
= {\rm n}_{a(P_{E}, \theta_{E})} < \infty$. 
}
\\ \
{\it Proof}. 
Let $a(M)$ (resp. $a(M)^{r}$) 
be the Lie algebra of holomorphic (1,0) (resp. real) 
vector fields on $M$. 
Let $J$ be the almost complex structure of $M$ 
on the real tangent bundle $T M$ and the complexification 
$T M \otimes_{\mathbb{R}} \mathbb{C}$. 
By $J [X, Y] = [X, J Y]$ ($X, Y \in a(M)^{r}$), 
the correspondence: 
$(a(M)^{r}, J) \to (a(M), \sqrt{-1} \cdot); 
Y \mapsto Y' := \frac{1}{2}(Y - \sqrt{-1} J Y)$ 
is a complex Lie algebra isomorphism, 
the inverse map of which is recovered as 
$Y' \mapsto (Y')^{r} := Y' + Y'' = Y$ 
with $Y'' := \frac{1}{2}(Y + \sqrt{-1} J Y)$. 
Since $M$ is compact, each  
$Y \in a(M)^{r}$ 
generates 1-parameter subgroup 
$\{ \alpha_{Y} (t) {;~} t \in \mathbb{R} \}$ 
of ${\rm Aut}(M)$ 
such that 
$\left. \frac{d}{d t} \right|_{t = 0} 
\alpha_{Y} (t) \xi = Y |_{\xi}$ ($\xi \in M$) 
(cf. \cite[p.35, Prop.I.1.6]{KN1963}). 
By \cite{BM1947}, 
${\rm Aut}(M)$ is a holomorphic Lie transformation group on $M$ 
such that $a(M)^{r}$ is finitely dimensional 
(cf. S.~Kobayashi \cite[Thm.III.1.1]{Ks1972}). 
The Lie algebra 
${\it Lie}~{\rm Aut}(M)$ 
of all left invariant vector fields on 
${\rm Aut}(M)$ in the sense of \cite{Cc1946} 
is isomorphic to $a(M)$ by 
${\it Lie}~{\rm Aut}(M) \to a(M); X \mapsto Y := (- X^{+})'$; 
$X^{+} |_{\xi} := \left. \frac{d}{dt} \right|_{t = 0} ({\rm exp} t X) \xi$ 
($\xi \in M$) 
as complex Lie algebras 
(cf. \cite[2.120]{Bal1987}, 
\cite[Thm.II.3.4]{Hs1978}). 
Since $a(M, E)$ is a complex Lie subalgebra of $a(M)$, 
${\rm Aut}(M, E)$ 
is a complex Lie subgroup of 
${\rm Aut}(M)$ 
corresponding to $a(M, E)$ as a closed subgroup of ${\rm Aut}(M)$ 
in compact-open topology by which the action map ${\rm Aut}(M) \times T' M \to T' M$ 
is continuous (cf. \cite[Lem.IV.14.2, Cor.]{Cc1946}, \cite[Thm.II.2.3]{Hs1978}). 
By Lemma 4.4, 
${\rm Aut}(P_{E}, \theta_{E})$ 
is a holomorphic Lie transformation group 
on $P_{E}$ such that 
${\cal P}: {\rm Aut}(M, E) \to {\rm Aut}(P_{E}, \theta_{E})$ 
is a complex Lie group isomorphism satisfying  
$p_{E} ({\cal P}(\alpha) \eta) = \alpha (p_{E} \eta)$ 
$(\alpha \in {\rm Aut}(M, E), \eta \in P_{E})$. 
Especially, 
${\rm n}_{\mathfrak{g}'} = {\rm n}_{a(M, E)} 
\leq {\rm n}_{a(M)} = {\rm n}_{a(M)^{r}}< \infty$. 
\qed 
\smallskip \\ \
\textsc{Lemma 4.6}. 
{\it 
Let $(M, E)$ be a compact c-manifold such that 
$L_{E}$ is immersive. 
Then $M$ is Fano, so that each connected component of $M$ 
is simply connected.
}
\\ \
{\it Proof}. 
By a theorem of S.~Kobayashi \cite{Ks1959}, 
$K_{M}^{*} \cong L_{E}^{\otimes {\rm k}_{M}}$ 
for ${\rm k}_{M} \in \mathbb{N}$ 
such as ${\rm n}_{M} = 2 {\rm k}_{M} - 1$. 
By Lemma 3.2, $K_{M}^{*}$ is immersive. 
By Proposition 3.3 (D), $M$ is Fano, 
so that $M$ is simply connected. 
\qed
\smallskip \\ \
{\bf 5. 
The adjoint varieties. 
} 
Let $V$ be a $\mathbb{C}$-linear space 
with the group 
$GL_{V} \mathbb{C}$ of all $\mathbb{C}$-linear transformations 
and a $\mathbb{C}$-bilinear form $S$. 
Put 
$S^{\perp} 
:= \{ X \in V {;~} S(X, Y) = 0~(Y \in V) \}$. 
Let ${\cal G}$ be a complex Lie algebra 
with the product 
${\cal G}^{\times 2} \to {\cal G}; 
(X, Y) \mapsto [X, Y]$,  
${\rm ad} X: {\cal G} \to {\cal G}; Y \mapsto [X, Y]$.  
Put 
$B_{\cal G}(X, Y) := {\rm trace}(({\rm ad} X) \circ ({\rm ad} Y))$ 
and 
${\rm Aut}({\cal G}) 
:= \{ g \in GL_{{\cal G}} \mathbb{C} {;~} 
g [X, Y] = [gX, gY] (X, Y \in {\cal G}) \}$. 
By definition, 
{\it a $\mathbb{Z}$-gradation of c-type for} ${\cal G}$ 
is a direct-sum decomposition 
${\cal G} = \oplus_{i \in \mathbb{Z}} {\cal G}_{i}$   
by some 
$\mathbb{C}$-linear subspaces ${\cal G}_{i}$ of ${\cal G}$ 
such that 
$[{\cal G}_{i}, {\cal G}_{j}] 
\subseteqq {\cal G}_{i+j}$ ($i, j \in \mathbb{Z}$),  
${\cal G}_{k} = 0$ ($k \geq 3$) and ${\rm n}_{{\cal G}_{2}} = 1$, 
where ${\cal G}_{i}$ for $i \leq 1$ is not assumed to be nonzero. 
It is known by \cite[4.2.Thm.]{Wja1965}, N.~Tanaka \cite{Tn1976} 
and H.~Asano \cite{KOY1999} 
that any simple complex Lie algebra 
$\mathfrak{g}$ of ${\rm n}_{\mathfrak{g}} > 1$. 
admits a $\mathbb{Z}$-gradation of c-type by means of the highest root. 
\\ \
Let 
$\mathfrak{g}$ be a simple complex Lie algebra of ${\rm n}_{\mathfrak{g}} > 1$. 
Let fix any $\mathbb{Z}$-gradation of c-type 
$\mathfrak{g} = \oplus_{i \in \mathbb{Z}} \mathfrak{g}_{i}$ 
and any non-zero element $e_{2}$ of $\mathfrak{g}_{2}$, 
by which put 
$\hat{\mathfrak{g}} 
:= \mathfrak{g}_{-2} \oplus [\mathfrak{g}_{-2}, \mathfrak{g}_{2}] \oplus \mathfrak{g}_{2}$, 
$\mathfrak{b}_{i, j} 
:= \{ X \in \mathfrak{g}_{i}{;~} 
B_{\mathfrak{g}}(X, Y) = 0~(Y \in \mathfrak{g}_{j}) \}$ 
and 
$\mathfrak{g}_{i, j} := \{ X \in \mathfrak{g}_{i}{;~} 
[X, Y] = 0~(Y \in \mathfrak{g}_{j}) \}$ 
for $i, j \in \mathbb{Z}$. 
By virtue of \cite{Tn1976}, 
one has the following results. 
\vspace{1mm} 
\\ \
\textsc{Lemma 5.0.} 
(A) 
{\it 
For 
$X = \sum_{i \in \mathbb{Z}} X_{i}$ 
and 
$Y = \sum_{i \in \mathbb{Z}} Y_{i}$ 
such as 
$X_{i}, Y_{i} \in \mathfrak{g}_{i}$, 
one has that 
$B_{\mathfrak{g}}(X, Y) 
= \sum_{i+j = 0} B_{\mathfrak{g}}(X_{i}, Y_{j})$, 
so that 
$(B_{\mathfrak{g}} 
|_{ \mathfrak{g}_{0}^{\times 2}})^{\perp} 
= \{ 0 \}$.  
For $i \in \mathbb{Z}$, 
$\mathfrak{b}_{i, -i} = 0$, 
${\rm n}_{\mathfrak{g}_{-i}} = {\rm n}_{\mathfrak{g}_{i}}$ 
and $\mathfrak{g}_{0, -i} = \mathfrak{g}_{0, i}$, 
so that 
$\mathfrak{g} = \oplus_{i = -2}^{2} \mathfrak{g}_{i}$ 
with unique 
$e_{-2} \in \mathfrak{g}_{-2}$ 
such that $\mathfrak{g}_{-2} = \mathbb{C} e_{-2}$ 
and $B(e_{2}, e_{-2}) = 1$.  
There exists unique 
$h_{0} \in \mathfrak{g}_{0}$ 
such that 
$\mathfrak{g}_{i} = 
\{ X \in \mathfrak{g}{;~} ({\rm ad} h_{0}) X 
= i \cdot X \}$ $(i \in \mathbb{Z})$, 
so that 
$B_{\mathfrak{g}} (h_{0}, h_{0}) = 8 + 2 {\rm n}_{\mathfrak{g}_{1}} \geq 8$ 
and $\mathfrak{g}_{0} = \mathbb{C} h_{0} \oplus \mathfrak{g}_{0,2}$. 
} 
\\ \
(B) 
{\it 
Fix any 
$j \in \{ -1, 1 \}$. 
Then 
$\mathfrak{g}_{0, j} = \mathfrak{g}_{0, -j} = 0$ 
if $\mathfrak{g}_{j} \ne 0$, 
and 
$\mathfrak{g} 
= \hat{\mathfrak{g}} 
= \mathbb{C} e_{-2} + \mathbb{C} h_{0} + \mathbb{C} e_{2}$ 
if $\mathfrak{g}_{j} = 0$. 
For ${\rm ad}_{j} (e_{-2j}): 
\mathfrak{g}_{j} \to \mathfrak{g}_{-j}; 
X \mapsto [e_{-2 j}, X]$, 
${\rm ad}_{j} (e_{-2j})$ 
is a $\mathbb{C}$-linear isomorphism, so that 
$\mathfrak{g}_{j, j} = 0$. 
Put 
\[
<\!\! *,~* \!\!>_{j}: 
\mathfrak{g}_{j}^{\times 2} \to \mathbb{C}; 
(Y', Y'') \mapsto <\!\! Y', Y'' \!\!>_{j} 
:= B_{\mathfrak{g}} (e_{-2j}, [Y', Y'']).
\]  
Then 
$<\!\! *,~* \!\!>_{j}^{\perp} = \{ 0 \}$ 
and 
$<\!\! [X, Y'], Y'' \!\!>_{j} 
+ <\!\! Y', [X, Y''] \!\!>_{j} = 0$ 
for any 
$X \in \mathfrak{g}_{0, 2}; 
Y', Y'' \in \mathfrak{g}_{j}$, 
so that 
$\mathfrak{g}_{0, 2} 
= \{ X \in \mathfrak{g}_{0} {;~}  
B_{\mathfrak{g}} (X, \mathbb{C} h_{0}) = 0 \}$, 
$(B_{\mathfrak{g}} 
|_{\mathfrak{g}_{0, 2}^{\times 2}})^{\perp} 
= \{ 0 \}$ 
and 
$\mathbb{C} h_{0} = \{ X \in \mathfrak{g}_{0} {;~} 
B_{\mathfrak{g}} (X, \mathfrak{g}_{0, 2}) = 0 
\}$. 
}
\\ \ 
(C) 
{\it 
Put 
$c_{0} 
:= \frac{2}{B_{\mathfrak{g}} (h_{0}, h_{0})} > 0$. 
Then 
$[e_{2}, e_{-2}] = c_{0} \cdot h_{0}$, 
$[h_{0}, e_{2j}] = 2j e_{2j}$ $(j \in \{ -1, 1 \})$ 
and 
$\hat{\mathfrak{g}} 
= \mathbb{C} e_{-2} \oplus \mathbb{C} h_{0} \oplus \mathbb{C} e_{2} 
\cong sl_{2} \mathbb{C}$, 
a simple complex Lie algebra of type $A_{1}$. 
} 
\\ \ 
(D) 
{\it 
Put 
$\mathfrak{n} := \mathfrak{g}_{1} \oplus \mathfrak{g}_{2}$, 
$\mathfrak{g}_{2'} := \{ X \in \mathfrak{g} {;~} 
[\mathfrak{n}, X] \subseteq \mathbb{C} X \}$ 
and 
$\mathfrak{n}' := \{ X \in \mathfrak{g} {;~} 
[\mathfrak{n}, X] \subseteq \mathfrak{g}_{2}, 
[\mathfrak{g}_{2}, X] \subseteq \{ 0 \} \}$. 
Then 
$\mathfrak{g}_{2'} = \mathfrak{g}_{2}$ 
and 
$\mathfrak{n}' = \mathfrak{n}$. 
} 
\vspace{1mm} 
\\ \
{\it Proof}. 
(A) 
By \'E. Cartan, 
$B_{\mathfrak{g}}^{\perp}$ 
is a solvable ideal of $\mathfrak{g}$, 
which is $\{ 0 \}$ 
since $\mathfrak{g}$ is simple with 
${\rm n}_{\mathfrak{g}} > 1$ 
(e.g. M.~Goto \& F.D.~Grosshans \cite[(2.1.3), Proof]{GG1978}). 
If 
$X_{i} \in \mathfrak{g}_{i}$ 
and 
$X_{j} \in \mathfrak{g}_{j}$ 
with 
$i + j \ne 0$, 
then 
$({\rm ad} X_{i}) \circ ({\rm ad} X_{i'}) 
\mathfrak{g}_{k} \subseteq 
\mathfrak{g}_{i+ i' + k}$ 
with $i + i' + k \ne k$, 
so that 
$B_{\mathfrak{g}} (X_{i}, X_{i'}) 
= 0$. 
Hence, 
$B_{\mathfrak{g}}(X, Y) 
= \sum_{i+j = 0} B_{\mathfrak{g}} (X_{i}, Y_{j})$, 
so that 
$B_{\mathfrak{g}} (\mathfrak{g}_{0}, 
\oplus_{i \ne 0} \mathfrak{g}_{i}) = 0$ 
and 
$(B_{\mathfrak{g}} 
|_{\mathfrak{g}_{0}^{\times 2}})^{\perp} 
= \mathfrak{g}_{0} \cap 
B_{\mathfrak{g}}^{\perp} 
= \{ 0 \}$. 
For $i \in \mathbb{Z}$, 
$\{ 0 \} = \mathfrak{g}_{i} \cap 
B_{\mathfrak{g}}^{\perp} 
= \mathfrak{b}_{i, -i}$, 
i.e., 
the {map} 
$\flat_{i}: 
\mathfrak{g}_{i} \to \mathfrak{g}_{-i}^{*}$ 
such as 
$<\!\! \flat_{i} (X), Y \!\!>_{\mathfrak{g}_{-i}} 
= B_{\mathfrak{g}}(X, Y)$ 
is injective, so that 
${\rm n}_{\mathfrak{g}_{i}} 
\leq {\rm n}_{\mathfrak{g}^{*}_{-i}} 
= {\rm n}_{\mathfrak{g}_{-i}}$. 
Hence, 
${\rm n}_{\mathfrak{g}_{i}} 
= {\rm n}_{\mathfrak{g}_{-i}}$. 
For $X \in \mathfrak{g}_{0, i}$, 
$0 = B_{\mathfrak{g}} 
(\mathfrak{g}_{-i}, [X, \mathfrak{g}_{i}]) 
= B_{\mathfrak{g}} ([\mathfrak{g}_{-i}, X], \mathfrak{g}_{i})$, 
so that $X \in \mathfrak{g}_{0, -i}$. 
By $\mathfrak{g}_{-2}^{*} 
\ni \flat_{2}(e_{2}) \ne 0$ 
and 
${\rm n}_{\mathfrak{g}_{-2}} 
= {\rm n}_{\mathfrak{g}_{2}} = 1$, 
there exists uniquely 
$e_{-2} \in \mathfrak{g}_{-2}$ 
such as 
$<\!\! \flat_{2}(e_{2}),~e_{-2} \!\!>_{\mathfrak{g}_{-2}} = 1$. 
After \cite[p.139]{Tn1976}, put 
$f: \mathfrak{g} \to \mathfrak{g}; 
\sum_{i = -2}^{2} X_{i} \mapsto 
\sum_{i = -2}^{2} i \cdot X_{i}$. 
Then 
$\mathfrak{g}_{i} = \{ X \in \mathfrak{g}{;~} 
f (X) = i \cdot X \}$ 
($i \in \mathbb{Z}$) 
and  
$f ([X, Y]) = [ f (X), Y ] + [ X, f (Y) ]$ 
($X, Y \in \mathfrak{g}$). 
Since $\mathfrak{g}$ is simple, 
there exists unique $h_{0} \in \mathfrak{g}$ 
such that $f = {\rm ad} h_0$, 
so that 
$h_{0} \in \mathfrak{g}_{0}$ 
by $f (h_{0}) = [ h_{0}, h_{0} ] = 0$. 
Moreover, 
$B_{\mathfrak{g}} (h_{0}, h_{0}) 
= {\rm trace}({\rm ad} h_{0})^2 
= \sum_{i = -2}^2 i^2 {\rm n}_{\mathfrak{g}_{i}} 
= 8 + 2 {\rm n}_{\mathfrak{g}_{1}}$. 
Note that 
$h_{0} \in \mathfrak{g}_{0} \backslash \mathfrak{g}_{0, 2}$, 
and that 
$\mathfrak{g}_{0, 2} \to \mathfrak{g}_{0} \to \mathfrak{g}_{2} \to \{ 0 \}$ 
is an exact sequence (cf. H.~Kaji \& O.Y.\cite{KY2004}), 
so that  
${\rm n}_{\mathfrak{g}_{0}} 
= 1 + {\rm n}_{\mathfrak{g}_{0,2}}$ 
and 
$\mathfrak{g}_{0} 
= \mathbb{C} h_{0} \oplus \mathfrak{g}_{0,2}$. 
\\ \ 
(B) 
After \cite[Proof of Lem.3.1]{Tn1976}, put 
$\mathfrak{a} := [\mathfrak{g}_{-j}, 
\mathfrak{g}_{-j}] 
\oplus \mathfrak{g}_{-j} 
\oplus [\mathfrak{g}_{-j}, \mathfrak{g}_{j}] 
\oplus \mathfrak{g}_{j} 
\oplus [\mathfrak{g}_{j}, \mathfrak{g}_{j}]$. 
If 
$\mathfrak{g}_{j} \ne \{ 0 \}$, 
then 
$\mathfrak{a} \ne \{ 0 \}$ 
and  
$\mathfrak{a} = \mathfrak{g}$ 
since $\mathfrak{a}$ is an ideal of $\mathfrak{g}$. 
In this case, 
for $X \in \mathfrak{g}_{0, j} = \mathfrak{g}_{0, -j}$, 
$\{ 0 \} = [X, \mathfrak{a}] = [X, \mathfrak{g}]$, 
so that $\mathbb{C} X = \{ 0 \}$ 
as a proper ideal of $\mathfrak{g}$, 
hence 
$\mathfrak{g}_{0, j} = \mathfrak{g}_{0, -j} = \{ 0 \}$. 
If 
$\mathfrak{g}_{j} = \{ 0 \}$, then 
$\mathfrak{g} 
= \oplus_{i = -2, 0, 2} \mathfrak{g}_{i}$, 
so that 
$\hat{\mathfrak{g}}$ 
and 
$\hat{\mathfrak{g}} + \mathbb{C} h_{0}$ 
are non-zero ideals of $\mathfrak{g}$. 
In this case, 
$\mathfrak{g} = \hat{\mathfrak{g}} 
= \hat{\mathfrak{g}} + \mathbb{C} h_{0} 
= \mathbb{C} e_{-2} + \mathbb{C} h_{0} + \mathbb{C} e_{2}$. 
To prove that ${\rm ad}_{j}(e_{-2j})$ 
is an isomorphism, it is enough to show that 
$({\rm ad} e_{-2j}) \mathfrak{g}_{j} 
= \mathfrak{g}_{-j}$, 
which is trivial if $\mathfrak{g}_{j} = \{ 0 \}$   
because of 
${\rm n}_{\mathfrak{g}_{-j}} 
= {\rm n}_{\mathfrak{g}_{j}} = \{ 0 \}$. 
When $\mathfrak{g}_{j} \ne \{ 0 \}$, 
by 
\cite[Proof of Lem.3.1]{Tn1976}, 
put 
$\mathfrak{b} 
:= \mathfrak{g}_{-2 j} \oplus 
[\mathfrak{g}_{-2 j}, \mathfrak{g}_{j}] 
\oplus 
[[\mathfrak{g}_{-2 j}, \mathfrak{g}_{j}], \mathfrak{g}_{j}] 
\oplus 
[[[\mathfrak{g}_{-2 j}, \mathfrak{g}_{j}], \mathfrak{g}_{j}], 
\mathfrak{g}_{j}] 
\oplus 
[[[[\mathfrak{g}_{-2 j}, \mathfrak{g}_{j}], \mathfrak{g}_{j}], 
\mathfrak{g}_{j}], \mathfrak{g}_{j}]$, 
which is an ideal of 
$\mathfrak{g} = \mathfrak{a}$ 
because of 
$[ \mathfrak{g}_{\pm j}, \mathfrak{b}] \subseteq \mathfrak{b}$. 
By the claim (A), 
$\{ 0 \} \ne \mathfrak{g}_{-2 j} 
\subseteq \mathfrak{b}$, 
so that 
$\mathfrak{g} = \mathfrak{b}$ 
and that 
$({\rm ad} e_{-2j}) \mathfrak{g}_{j} 
= [\mathfrak{g}_{-2 j}, \mathfrak{g}_{j}] 
= \mathfrak{g}_{- j}$, 
as required. 
Then 
$B_{\mathfrak{g}} 
([X, \mathfrak{g}_{j}], e_{-2 j}) 
= B_{\mathfrak{g}} 
(X, [\mathfrak{g}_{j}, e_{-2 j}]) 
= B_{\mathfrak{g}} (X, \mathfrak{g}_{-j})$ 
($X \in \mathfrak{g}_{j}$), 
so that 
$\mathfrak{g}_{j, j} \subseteq 
\mathfrak{b}_{j, -j} = \{ 0 \}$.  
Hence, 
$\mathfrak{g}_{j, j} = \{ 0 \}$ 
and 
$<\!*,*\!>_{j}^{\perp} = \{ 0 \}$, 
as required. 
For any $X \in \mathfrak{g}_{0, 2} 
= \mathfrak{g}_{0, -2}$, 
$<\!\![X, Y'], Y''\!\!>_{j} 
+ <\!\!Y', [X, Y'']\!\!>_{j} 
= B_{\mathfrak{g}} 
(e_{-2j}, [[X, Y'], Y''] + [Y', [X, Y'']]) 
= B_{\mathfrak{g}} (e_{-2j}, [X, [Y',Y'']]) 
= B_{\mathfrak{g}} ([e_{-2j}, X], [Y',Y'']) 
= 0$. 
If 
$\mathfrak{g}_{j} \ne \{ 0 \}$, 
then 
${\rm ad} X |_{\mathfrak{g}_{j}}$ 
is an anti-symmetric operator 
w.r.t. the symplectic form 
$<\!*,*\!>_{j}$ 
on $\mathfrak{g}_{j}$, 
so that 
${\rm trace}({\rm ad} X |_{\mathfrak{g}_{j}}) 
= 0$. 
If 
$\mathfrak{g}_{j} = \{ 0 \}$, 
then 
${\rm trace} ({\rm ad} X |_{\mathfrak{g}_{j}}) 
= {\rm trace} (0) = 0$. 
In both cases, 
$B_{\mathfrak{g}}(X, h_{0}) 
= \sum_{j = -1, 1} 
j \cdot {\rm trace}({\rm ad} X 
|_{\mathfrak{g}_{j}})  = 0$ 
and 
$B_{\mathfrak{g}}(\mathfrak{g}_{0, 2}, h_{0}) = 0$. 
By the claim (A), 
$B_{\mathfrak{g}}(h_{0}, h_{0}) \ne 0$, 
$(B_{\mathfrak{g}} |_{\mathfrak{g}_{0}^{\times 2}})^{\perp} 
= \{ 0 \}$ 
and  
$\mathfrak{g}_{0} 
= \mathbb{C} h_{0} \oplus \mathfrak{g}_{0, 2}$. 
Hence,  
$\mathfrak{g}_{0, 2} = 
\{ X \in \mathfrak{g}_{0} {;~} 
B_{\mathfrak{g}} (X, h_{0}) = 0 \}$, 
$(B_{\mathfrak{g}} 
|_{\mathfrak{g}_{0, 2}^{\times 2}})^{\perp}  
= \mathfrak{g}_{0, 2} 
\cap (B_{\mathfrak{g}} 
|_{\mathfrak{g}_{0}^{\times 2}})^{\perp} 
= \{ 0 \}$ 
and 
$\mathbb{C} h_{0} = \{ X \in \mathfrak{g}_{0} {;~}  
B_{\mathfrak{g}} 
(X, \mathfrak{g}_{0, 2}) = 0 \}$, 
as required. 
\\ \
(C) 
By the claim (A), 
$[e_{2}, e_{-2}] = c \cdot h_{0} + X$ 
for some $X \in \mathfrak{g}_{0, 2}$ 
and $c \in \mathbb{C}$.  
By the claim (B), 
$B_{\mathfrak{g}}(\mathfrak{g}_{0, 2}, X) 
= B_{\mathfrak{g}}(\mathfrak{g}_{0, 2}, [e_{2}, e_{-2}]) 
= B_{\mathfrak{g}} ([\mathfrak{g}_{0, 2}, e_{2}], e_{-2}) 
= B_{\mathfrak{g}} (\{ 0 \}, e_{-2}) = \{ 0 \}$, 
so that $X = 0$. 
Hence, 
$[e_{2}, e_{-2}] = c \cdot h_{0}$ 
for some $c \in \mathbb{C}$. 
By the claim (A), 
$c B_{\mathfrak{g}} (h_{0}, h_{0}) 
= B_{\mathfrak{g}} (h_{0}, [e_{2}, e_{-2}]) 
= B_{\mathfrak{g}} ([h_{0}, e_{2}], e_{-2}) = 2$ 
and $B_{\mathfrak{g}}(h_{0}, h_{0}) > 0$. 
Then 
$c = c_{0} > 0$ and 
$\hat{\mathfrak{g}} 
= \mathbb{C} \hat{e}_{-2} \oplus \mathbb{C} h_{0} \oplus \mathbb{C} \hat{e}_{2}$. 
Put 
$\hat{e}_{-2} := e_{-2} / c_{0}$ 
and $\hat{e}_{2} := e_{2}$, 
so that 
$[\hat{e}_{2}, \hat{e}_{-2}] = h_{0}$. 
By 
$[h_{0}, \hat{e}_{2 j}] = 2j \hat{e}_{2 j}$, 
$\hat{\mathfrak{g}} 
\cong  \mathbb{C} [\delta_{i-1, j}]_{2 \times 2} 
\oplus \mathbb{C} [ (-1)^{i+1} \delta_{i,j}]_{2 \times 2}  
\oplus  \mathbb{C} [\delta_{i+1, j}]_{2 \times 2} 
= sl_{2} \mathbb{C}$ of type $A_{1}$. 
\\ \ 
(D) 
For the first claim, it follows that 
$\mathfrak{g}_{2} \subseteq \mathfrak{g}_{2'}$ 
from the $\mathbb{Z}$-gradation of c-type. 
To prove that $\mathfrak{g}_{2'} \subseteq \mathfrak{g}_{2}$, 
take any $X = \sum_{i = -2}^{2} X_{i} \in \mathfrak{g}_{2'} \backslash \{ 0 \}$ 
($X_{i} \in \mathfrak{g}_{i}$), 
so that there exists $c \in \mathbb{C}$ such that 
$c X = [e_{2}, X] = \sum_{i = -2}^{0} [e_{2}, X_{i}] 
\in \oplus_{j \geq 0} \mathfrak{g}_{j}$. 
Suppose that $X_{-2} \ne 0$ or $X_{-1} \ne 0$. 
Then $c = 0$, 
so that $[e_{2}, X_{i}] = 0$ $(i \in \{ -2, -1, 0 \})$. 
By the claim (B), $X_{-1} = 0$. 
Hence $X_{-2} \ne 0$. 
$0 = B_{\mathfrak{g}}(h_{0}, 0) 
= B_{\mathfrak{g}}(h_{0}, [e_{2}, X_{-2}]) 
= 2 B_{\mathfrak{g}}(e_{2}, X_{-2})$ 
because of $[h_{0}, e_{2}] = 2 e_{2}$ 
(by the claim (A)). 
Because of $\mathfrak{b}_{-2, 2} = \{ 0 \}$ 
(by the claim (A)), $X_{-2} = 0$, 
a contradiction. 
Hence, $X_{-2} = X_{-1} = 0$, 
so that $X = \sum_{i = 0}^{2} X_{i}$. 
Then 
$c X = [e_{2}, X] = [e_{2}, X_{0}] \in \mathfrak{g}_{2}$. 
Suppose that $X_{0} \ne 0$. 
Then $c = 0$, so that $X_{0} \in \mathfrak{g}_{0, 2}$. 
If $\mathfrak{g}_{1} = \{ 0 \}$, 
then 
$\mathfrak{g}_{0} = \mathbb{C} h_{0}$ 
(by the claim (B)), 
so that 
$X_{0} \in \mathbb{C} h_{0} \cap \mathfrak{g}_{0, 2} = \{ 0 \}$ 
(by the claim (A)), a contradiction. 
Hence, $\mathfrak{g}_{1} \ne \{ 0 \}$, 
so that $\mathfrak{g}_{0, 1} = \{ 0 \}$ 
(by the claim (B)): 
By $X_{0} \ne 0$, there exists 
$Y_{1} \in \mathfrak{g}_{1} \subseteq \mathfrak{n}$ 
such that $[Y_{1}, X_{0}] \ne 0$. 
By $X \in \mathfrak{g}_{2'}$, 
there exists $c_{1} \in \mathbb{C}$ 
such that 
$c_{1} X = [Y_{1}, X] = \sum_{i = 0}^{1} [Y_{1}, X_{i}] \ne 0$. 
Then 
$c_{1} \ne 0$, $X = \sum_{i=0}^{1} [Y_{1}, X_{i}]/c_{1} \in \mathfrak{n}$ 
and $X_{0} = 0$, a contradiction. 
Hence, $X_{0} = 0$, 
so that $X = \sum_{i = 1}^{2} X_{i}$. 
If $X_{1} \ne 0$, 
then 
$[Y', X_{1}] \ne 0$ 
for some $Y' \in \mathfrak{g}_{1}$ 
(by the claim (B)). 
In this case, there exists $c' \in \mathbb{C}$ 
such that 
$c' X = [Y', X] = [Y', X_{1}] \ne 0$, 
so that $c' \ne 0$. 
Then 
$\sum_{i=1}^{2} X_{i} = X = [Y', X_{1}] /c' \in \mathfrak{g}_{2}$, 
so that $X_{1} = 0$, a contradiction. 
Hence, $X_{1} = 0$, 
so that $X = X_{2} \in \mathfrak{g}_{2}$, 
as required. 
For the last claim, it follows that 
$\mathfrak{n} \subseteq \mathfrak{n}'$ 
from the $\mathbb{Z}$-gradation of c-type. 
To prove that 
$\mathfrak{n}' \subseteq \mathfrak{n}$, 
take any $X = \sum_{i = -2}^{2} X_{i} \in \mathfrak{n}'$ 
$(X_{i} \in \mathfrak{g}_{i})$. 
By $X \in \mathfrak{n}'$, 
$0 = [e_{2}, X] = \sum_{i = -2}^{0} [e_{2}, X_{i}]$, 
so that $[e_{2}, X_{i}] = 0$ $(i \in \{ -2, -1, 0 \})$ 
because of $\mathbb{Z}$-gradation of c-type. 
Then 
$0 = B_{\mathfrak{g}}(h_{0}, 0) 
= B_{\mathfrak{g}}(h_{0}, [e_{2}, X_{-2}]) 
= 2 B_{\mathfrak{g}}(e_{2}, X_{-2})$, 
so that $X_{-2} = 0$ 
because of $\mathfrak{b}_{-2, 2} = \{ 0 \}$ 
(by the claim (A)). 
Moreover, 
$\{ 0 \} = B_{\mathfrak{g}}(\mathfrak{g}_{-1}, 0) 
= B_{\mathfrak{g}}(\mathfrak{g}_{-1}, [e_{2}, X_{-1}]) 
= B_{\mathfrak{g}}(\mathfrak{g}_{1}, X_{-1})$ 
because of $\mathfrak{g}_{1} = [\mathfrak{g}_{-1}, e_{2}]$ 
(by the claim (B)), 
so that $X_{-1} = 0$ 
because of $\mathfrak{b}_{-1, 1} = \{ 0 \}$ 
(by the claim (A)). 
Hence, 
$X = \sum_{i = 0}^{2} X_{i}$ 
and $X_{0} \in \mathfrak{g}_{0, 2}$. 
Suppose that $X_{0} \ne 0$. 
If $\mathfrak{g}_{1} = \{ 0 \}$, 
then $\mathfrak{g}_{0} = \mathbb{C} h_{0}$ 
(by the claim (B)), so that 
$X_{0} \in \mathbb{C} h_{0} \cap \mathfrak{g}_{0, 2} = \{ 0 \}$ 
(by the claim (A)), a contradiction. 
Hence, $\mathfrak{g}_{1} \ne \{ 0 \}$, 
so that $\mathfrak{g}_{0, 1} = \{ 0 \}$ 
(by the claim (B)): 
By $X_{0} \ne 0$, there exists $Y_{1} \in \mathfrak{g}_{1}$ 
such that $\mathfrak{g}_{1} \ni [Y_{1}, X_{0}] \ne 0$. 
Then 
$[Y_{1}, X] = \sum_{i = 0}^{2} [Y_{1}, X_{i}] \not\in \mathfrak{g}_{2}$, 
so that $X \not\in \mathfrak{n}'$, a contradiction. 
Hence, 
$X_{0} = 0$, 
so that $X = \sum_{i = 1}^{2} X_{i} \in \mathfrak{n}$, 
as required.
~\qed 
\smallskip \\ \
Put 
$G := {\rm Aut}(\mathfrak{g})$ 
with the identity connected component 
$G^{\circ} := {\rm Ad}(\mathfrak{g})$, 
the connected adjoint orbit 
$P_{e_{2}} := G^{\circ} \cdot e_{2}$ 
of $e_{2}$ and 
{\it the adjoint variety 
$X_\mathfrak{g} := p_\mathfrak{g}(P_{e_{2}})$ 
of $\mathfrak{g}$} 
defined for the fixed $\mathbb{Z}$-gradation of c-type 
and $0 \ne e_{2} \in \mathfrak{g}_{2}$. 
For any $c \in \mathbb{C}$, put 
\[
\mathfrak{W}_{c} := 
\{ X \in \mathfrak{g} \backslash \{ 0 \} 
{;~} ({\rm ad} X)^{2} Y 
= - 2 c \cdot B_{\mathfrak{g}}(X, Y) X~(Y \in \mathfrak{g}) \}. 
\]
\
\textsc{Lemma 5.1.} 
(A) 
{\it 
Let 
$h_{0} \in \mathfrak{g}_{0}$ 
be given as 
$\mathfrak{g}_{i} = \{ X \in \mathfrak{g}{;~} [h_{0}, X] =i X \}$ 
$(i \in \mathbb{Z})$, by which put 
$c_{\mathfrak{g}} := \frac{2}{B_{\mathfrak{g}}(h_{0}, h_{0})}$. 
Then 
$P_{e_{2}} = G \cdot e_{2} = \mathfrak{W}_{c_{\mathfrak{g}}}$, 
so that 
$X_{\mathfrak{g}}$ 
is a $G$-stable subvariety of $P_{\mathfrak{g}} \mathbb{C}$. 
Moreover, 
$\mathfrak{g} 
= \{ \sum_{j} v_{j} 
{;~} v_{j} \in P_{e_{2}} \}$, 
so that the affine action of $G^{\circ}$ 
on $P_{e_{2}}$ is almost effective, 
i.e., 
\[
\{ X \in \mathfrak{g} 
{;~} ({\rm ad} X) P_{e_{2}} = \{ 0 \} \} = \{ 0 \}. 
\]
}
\ 
(B) 
{\it 
Let $\mathfrak{g}'$ 
be a semisimple complex Lie algebra 
of ${\rm n}_{\mathfrak{g}'} > 1$ 
with a point 
$e' \in \mathfrak{g}' \backslash \{ 0 \}$ 
and 
$P' := {\rm Ad}(\mathfrak{g}') \cdot e'$ 
such that  
$\{ X \in \mathfrak{g}
{;~} ({\rm ad} X) P' = \{ 0 \} \} = \{ 0 \}$. 
Put 
$M' := \pi_{\mathfrak{g}'}(P')$. 
Assume that 
$P' = \pi_{\mathfrak{g}'}^{-1}(M')$, 
and that 
$M'$ is compact in $P_{\mathfrak{g}'} \mathbb{C}$. 
Then $\mathfrak{g}'$ is simple. 
Moreover, 
let $\mathfrak{g}' = \oplus_{i = -2}^{2} \mathfrak{g}_{i}$ 
be any $\mathbb{Z}$-gradation of c-type for $\mathfrak{g}'$ 
such that 
$\mathfrak{g}_{2} = \mathbb{C} e_{2}$. 
Then 
$P' = P_{e_{2}}$ and $M' = X_{\mathfrak{g}'}$ 
by this 
$\mathbb{Z}$-gradation of c-type for $\mathfrak{g}'$. 
}
\smallskip \\ \
{\it Proof}. 
(A) 
Put $c := c_{\mathfrak{g}}$. 
For $Y \in \mathfrak{g}$, 
put $k := B_{\mathfrak{g}} (e_{2}, Y) \in \mathbb{C}$. 
By Lemma 5.0 (A, C), 
$({\rm ad} e_{2})^{2} Y  
= ({\rm ad} e_{2})^{2} (k e_{-2}) 
= k [e_{2}, c h_{0}] 
= - 2 c k e_{2} 
= - 2 c B_{\mathfrak{g}}(e_{2}, Y) e_{2}$, 
so that 
$\mathfrak{W}_{c} \ni e_{2}$ 
and 
$P_{e_{2}} 
\subseteq G \cdot e_{2} 
\subseteq G \cdot \mathfrak{W}_{c} 
\subseteq \mathfrak{W}_{c}$. 
Conversely, by the method of H.~Freudenthal 
given in \cite{KOY1999}, put 
$V := \{ \sum_{j = 1}^{\ell} v_{j}{;~} 
v_{j} \in \mathfrak{W}_{c}, 
B_{\mathfrak{g}}(v_{j}, P_{e_{2}}) = \{ 0 \} 
~(j \in \{ 1, \cdots, \ell \}), \ell \in \mathbb{N} \}$, 
so that 
$G \cdot V \subseteq V$. 
Then 
$[ \mathfrak{g}, V ] \subseteq V$. 
Hence, 
$V$ is an ideal of $\mathfrak{g}$. 
Suppose that 
$V \ne \{ 0 \}$. 
Since $\mathfrak{g}$ is simple, 
$V = \mathfrak{g} \ni e_{-2} = \sum_{j} v_{j}$ 
with $v_{j} \in \mathfrak{W}_{c}$ 
such as 
$\{ 0 \} = B_{\mathfrak{g}} 
(v_{j}, P_{e_{2}}) 
\ni B_{\mathfrak{g}} (v_{j}, e_{2})$ 
for all $j$, so that 
$1 = B_{\mathfrak{g}} 
(e_{-2}, e_{2}) = \sum_{j} 
B_{\mathfrak{g}} (v_{j}, e_{2}) = 0$, 
a contradiction. 
Hence, 
$V = \{ 0 \}$. 
For any 
$Y_{1} \in \mathfrak{W}_{c}$, 
$0 \ne Y_{1} \not\in V = \{ 0 \}$, 
so that 
$B_{\mathfrak{g}} (Y_{1}, Y_{2}) \ne 0$ 
for some 
$Y_{2} \in P_{e_{2}}$. 
For $z \in \mathbb{C}$ and $i, j \in \{ 1, 2 \}$ 
such as $i \ne j$, put 
$Y_{i, j}(z) 
:= e^{z~{\rm ad} Y_{i}} Y_{j} 
= Y_{j} + z [ Y_{i}, Y_{j} ] 
- z^{2} c B( Y_{i}, Y_{j} ) Y_{i}$. 
There exist 
$z' \in \mathbb{C}$ such that 
$(z')^{2} c B_{\mathfrak{g}} (Y_{1}, Y_{2}) = 1$, 
so that 
$Y_{i, j}(z') 
= Y_{j} + z' [Y_{i}, Y_{j}] - Y_{i} 
= - Y_{j, i}(z')$. 
Then 
$\pi_{\mathfrak{g}}( Y_{2, 1}(z') ) 
= \pi_{\mathfrak{g}}( Y_{1, 2}(z') ) 
= \pi_{\mathfrak{g}}(
e^{(z'~{\rm ad} Y_{1})} Y_{2}) 
\in \pi_{\mathfrak{g}}(P_{e_{2}})$, 
so that 
$\pi_{\mathfrak{g}}(Y_{1}) 
= \pi_{\mathfrak{g}} 
( (e^{(z'~{\rm ad} Y_{2})})^{-1} Y_{2, 1}(z')) 
\in \pi_{\mathfrak{g}} (P_{e_{2}})$ 
and $k_{1} Y_{1} = g_{1} \cdot e_{2}$ 
for some $(k_{1}, g_{1}) 
\in \mathbb{C}^{\times} \times G^{\circ}$. 
Take 
$z_{1} \in \mathbb{C}$ 
such as 
$k_{1}^{-1} = e^{2 z_{1}}$. 
Then 
$P_{e_{2}} 
\ni g_{1} \cdot (e^{z_{1} {\rm ad} h_{0}} e_{2}) 
= g_{1} \cdot (e^{2 z_{1}} e_{2}) 
= k_{1}^{-1} g_{1} \cdot e_{2} 
= Y_{1}$. 
Hence, 
$P_{e_{2}} \supseteq \mathfrak{W}_{c}$. 
Then 
$P_{e_{2}} = \mathfrak{W}_{c}$, 
as required. 
\\ \
(B1)  
Since 
$\mathfrak{g}'$ 
is semisimple, 
there exist $r \in \mathbb{N}$ and simple ideals 
$\mathfrak{a}_j$ of $\mathfrak{g}'$ 
($1 \leq j \leq r$) 
such that 
$\mathfrak{g}' 
= \oplus_{j = 1}^{r} \mathfrak{a}_{j}$ 
(e.g. \cite[(2.1.4)]{GG1978}), 
so that 
${\rm Ad}(\mathfrak{g}') \cong 
\Pi_{j = 1}^{r} {\rm Ad}(\mathfrak{a}_{j})$. 
Take some $X \in P'$ and $X_{j} \in \mathfrak{a}_{j}$ 
such that $X = \oplus_{j = 1}^r X_{j}$. 
Then  
$P' = {\rm Ad}(\mathfrak{g}') X 
= \oplus_{j = 1}^{r} {\rm Ad}(\mathfrak{a}_{j}) \cdot X_{j}$. 
By $X \ne 0$, 
there exists 
$j_{1} \in \{ 1, \ldots, r \}$ 
such that $X_{j_{1}} \ne 0$. 
{\it Claim. $P' \subseteq \mathfrak{a}_{j_{1}}$:} 
In fact, put 
$\hat{X} := X - X_{j_{1}}$, 
so that the $\mathfrak{a}_{j_{1}}$-component of $\hat{X}$ equals $0$. 
{\it Suppose that $\hat{X} \ne 0$:} 
Because of $X \in P'$, 
for any $\lambda \in \mathbb{C}^{\times}$, 
$\lambda \cdot X 
\in \pi_{\mathfrak{g}'}^{-1} 
(\pi_{\mathfrak{g}'} X) \subseteq P' = \oplus_{j = 1}^{r} 
{\rm Ad}(\mathfrak{a}_{j}) \cdot X_{j}$. 
There exist 
$a_{\lambda, j} \in {\rm Ad}(\mathfrak{a}_{j})$ 
such that 
$\lambda X 
= \oplus_{j} a_{\lambda, j} X_{j}$. 
For $j \in \{ 1, \ldots, r \}$, put 
$b_{\lambda, j} := a_{\lambda, j_{1}}$ 
(if $j = j_{1}$) 
or 
${\rm id}_{\mathfrak{a}_{j}}$ 
(if $j \ne j_{1}$), 
and 
$b_{\lambda} := \Pi_{j = 1}^{r} b_{\lambda, j} 
\in {\rm Ad}(\mathfrak{g}')$. 
Then 
$P' \ni b_{\lambda} X 
= \lambda X_{j_{1}} + \oplus_{j \ne j_{1}} X_{j}$. 
Put 
$\lambda_{\ell} := 1 / \ell \in \mathbb{C}^{\times}$ 
for $\ell \in \mathbb{N}$. 
For each $j \in \{ 1, \ldots, r \}$, 
take a Euclidian metric on $\mathfrak{a}_{j}$, 
by which there exists the direct-sum metric on 
$\mathfrak{g}'$. 
Then 
$\lim_{\ell \rightarrow \infty} 
b_{\lambda_{\ell}} (X) 
= \lim_{\ell \rightarrow \infty} 
\lambda_{\ell} X_{j_{1}} 
+ \oplus_{j \ne j_{1}} X_{j} 
= \hat{X} \ne 0$ by assumption, 
which is contained in 
${\rm cl}_{\mathfrak{g}'}(\{ b_{\lambda} X{;~} \lambda \in \mathbb{C}^{\times} \}) 
\backslash \{ 0 \} 
\subseteq {\rm cl}_{\mathfrak{g'}}(P') \backslash \{ 0 \}$. 
Since $\pi_{\mathfrak{g}'}$ is continuous 
and $M'$ is closed as a compact subset of $P_{\mathfrak{g}'} \mathbb{C}$, 
$P' = \pi_{\mathfrak{g}'}^{-1}(M')$ is closed in 
$\mathfrak{g}' \backslash \{ 0 \}$, 
so that ${\rm cl}_{\mathfrak{g'}}(P') \backslash \{ 0 \} = P'$. 
Then $\hat{X} \in P' = \oplus_{j = 1}^{r} {\rm Ad}(\mathfrak{a}_{j}) X_{j}$. 
On the other hand, each element of 
$\oplus_{j = 1}^{r} {\rm Ad}(\mathfrak{a}_{j}) X_{j}$ 
has non-zero $\mathfrak{a}_{j_{1}}$-component 
because of $X_{j_{1}} \ne 0$, {\it a contradiction}. 
Hence, 
$\hat{X} = 0$ and $X = X_{j_{1}}$, 
so that $P' \subseteq \mathfrak{a}_{j_{1}}$, 
{\it as required}. 
For any $Y \in \mathfrak{a}_{j}$ with $j \ne j_{1}$, 
$({\rm ad} Y) P' \subseteq \mathfrak{a}_{j} \cap \mathfrak{a}_{j_{1}} = \{ 0 \}$, 
so that $Y = 0$ by the assumption of almost effectivity of 
${\rm Ad}(\mathfrak{g}')$ on $P'$. 
Hence, $\mathfrak{g}' = \mathfrak{a}_{j_{1}}$, 
which is simple as required. 
\\ \ 
(B2) 
Since $M'$ is a compact orbit of ${\rm Ad}(\mathfrak{g}')$ 
with respect to the projective adjoint action, 
$M'$ is a complex submanifold of $P_{\mathfrak{g}'} \mathbb{C}$. 
By W-L.~Chow \cite{Cwl1949}'s theorem (cf. \cite[IV.5.13]{FG2002}), 
$M'$ is a projective algebraic variety. 
By (B1), 
$\mathfrak{g}' = \mathfrak{a}_{j_{1}}$ 
is simple. 
For any $\mathbb{Z}$-gradation 
$\mathfrak{g}' = \oplus_{i = -2}^{2} \mathfrak{g}_{i}$ 
of c-type, 
put $\mathfrak{n} := \mathfrak{g}_{1} \oplus \mathfrak{g}_{2}$, 
$N := \{ g \in {\rm Aut}(\mathfrak{g}'); 
g \mathfrak{n} \subseteq \mathfrak{g}_{2}, g X = X (X \in \mathfrak{g}_{2}) \}$ 
and 
${\it Lie} N := \{ X \in \mathfrak{g}'; e^{t {\rm ad} X} \in N (t \in \mathbb{R}) \}$. 
Then 
$N$ is an algebraic group in $GL_{\mathfrak{g}'} \mathbb{C}$ 
and 
${\it Lie} N = \mathfrak{n}' = \mathfrak{n}$ by Lemma 5.0 (D), 
so that $N$ is a nilpotent algebraic group in 
$GL_{\mathfrak{g}} \mathbb{C}$. 
Let $N^{z}$ be the Zariski identity connected component of $N$ 
in $GL_{\mathfrak{g}'} \mathbb{C}$, 
and $N^{\circ}$ be the identity connected component of $N$ in usual toplology. 
Then $N^{z} = N^{\circ}$ 
(cf. \cite[p.28, Remark 2.6]{ON2015}). 
Hence, $N^{o}$ is a Zariski connected solvable 
complex algebraic group in the algebraic group 
${\rm Ad}(\mathfrak{g}')$, 
which acts on the projective algebraic variety $M'$, 
projective linearly (cf. \cite[1.2, 1.3, 1.4, 4.1, 5.6]{ON2015}). 
By \cite[III.10.4Thm., AG.7.4Thm.]{Ba1991}, 
there exists $\xi' \in M'$ such that 
$N^{o} \xi' = \{ \xi' \}$ 
(cf. \cite[Thm.11.15, Thm.11.9]{ON2015}). 
Take some 
$\eta' \in \pi_{\mathfrak{g}}^{-1} (\xi')$. 
Then $\eta' \ne 0$ 
and $[ \mathfrak{g}_{1} \oplus \mathfrak{g}_{2}, \eta'] 
\subseteqq \mathbb{C} \eta'$.  
By Lemma 5.0 (D), 
$\eta' \in \mathfrak{g}_{2} = \mathbb{C} e_{2}$. 
Then 
$e_{2} \in \mathbb{C} \eta' = \pi_{\mathfrak{g}'}^{-1}(\xi') 
\subseteq \pi_{\mathfrak{g}'}^{-1}(M') = P'$. 
Hence, 
$P' 
= {\rm Ad}(\mathfrak{g}') \cdot e_{2} 
= P_{e_{2}}$ 
and 
$X_{\mathfrak{g}'} 
= \pi_{\mathfrak{g}'} (P_{e_{2}}) 
= \pi_{\mathfrak{g}'} (P') = M'$. 
~\qed 
\smallskip \\ \
For  
$G = {\rm Aut}(\mathfrak{g})$,  
put 
$G_{e_{2}} := \{ g \in G {;~} 
g (e_{2}) = e_{2} \}$ 
and 
$G_{+} := \{ g \in G {;~} 
g (e_{2}) \in \mathbb{C} e_{2} \}$. 
For $g \in G$, put 
$R_{g}: G \to G; 
\alpha \mapsto \alpha \circ g$, 
$L_{g}: G \to G; 
\alpha \mapsto g \circ \alpha$ 
and 
$A_{g} := L_{g} \circ R_{g^{-1}}: 
G \to G$. 
For $X \in \mathfrak{g}$, put 
$X_{R} |_{\alpha} 
:= \left. \frac{d}{dt} \right|_{t=0} 
R_{e^{t({\rm ad} X)}} (\alpha) 
\in T_{\alpha} G$, 
the real tangent space of $G$ 
at $\alpha \in G$ 
in the sense of \cite{Cc1946}. 
Then 
$d L_{g} X_{R} = X_{R}$, 
$d A_{g} X_{R} = (g X)_{R}$ 
and 
$[ X_{R}, Y_{R} ] = [X, Y]_{R}$ 
for $X, Y \in \mathfrak{g}$ 
(cf. \cite[p.42, Prop.4.1]{KN1963}). 
On 
$G = {\rm Aut}(\mathfrak{g})$, 
a holomorphic (1,0)-form $\theta_{G}$ 
is well-defined such that 
$\theta_{G} (X_{R}) := B_{\mathfrak{g}} (X, e_{2})$ 
for all $X \in \mathfrak{g}$. 
By 
\cite[p.55, Example 5.1]{KN1963}, 
there exist unique complex manifold structure 
on the topological coset space 
$G/G_{e_{2}}$ (resp. $G/G_{+}$) 
such that $G$ is a holomorphic fiber bundle on 
$G/G_{e_{2}}$ (resp. $G/G_{+}$) with 
$\pi_{e_{2}}: G \to G/G_{e_{2}}; 
\alpha \mapsto \alpha G_{e_{2}}$ 
(resp. 
$\pi_{+}: G \to G/G_{+}; 
\alpha \mapsto \alpha G_{+}$) 
as the holomorphic projection. 
There exists a holomorphic (1,0)-form 
$\breve{\theta}$ on $G/G_{e_{2}}$ 
such that 
$\theta_{G} = \pi_{e_{2}}^{*} \breve{\theta}$. 
Put 
$\breve{p}: 
G/G_{e_{2}} \to G/G_{+}; 
\pi_{e_{2}}(\alpha) \mapsto \pi_{+}(\alpha)$ 
such as 
$\pi_{+} = \breve{p}~\circ~\pi_{e_{2}}$,  
by which $G/G_{e_{2}}$ 
is a holomorphic principal 
$G_{+}/G_{e_{2}}$-bundle on $G/G_{+}$, 
where $G_{e_{2}}$ 
is a normal subgroup of $G_{+}$ 
because of 
$\beta^{-1} G_{e_{2}} \beta = G_{e_{2}}$ 
($\beta \in G_{+}$). 
Put 
$\Phi_{e_{2}}: 
G/G_{e_{2}} \to P_{e_{2}}; 
\pi_{e_{2}}(\alpha) 
\mapsto \alpha (e_{2})$ 
and 
$\varphi_{e_{2}}: 
G/G_{+} \to X_{\mathfrak{g}}; 
\pi_{+}(\alpha) \mapsto 
\pi_{\mathfrak{g}}(\alpha (e_{2}))$, 
as holomorphic injective immersions. 
\smallskip \\ \
\textsc{Proposition 5.2}. 
(A) 
{\it 
$P_{e_{2}}$ 
is a complex regular submanifold of $\mathfrak{g}$ 
biholomorphic to $G/G_{e_{2}}$ by $\Phi_{e_{2}}$, 
so that there exists unique holomorphic (1,0)-form 
$\theta_{\mathfrak{g}}$ 
on $P_{e_{2}}$ 
such that 
$\breve{\theta} 
= \Phi_{e_{2}}^{*} \theta_{\mathfrak{g}}$, 
so that 
$\theta_{G} = 
\pi_{e_{2}}^{*} (\Phi_{e_{2}}^{*} \theta_{\mathfrak{g}})$. 
Moreover, 
$X_{\mathfrak{g}}$ 
is a compact complex regular submanifold of 
$P_{\mathfrak{g}} \mathbb{C}$ 
bihomorphic to $G/G_{+}$ 
by $\varphi_{e_{2}}$. 
}
\\ \ 
(B) 
{\it 
${\rm n}_{X_{\mathfrak{g}}} 
= 1 + {\rm n}_{\mathfrak{g}_{-1}}$, 
which equals $1$ if and only if 
${\rm rank}(\mathfrak{g}) = 1$. 
When ${\rm n}_{X_{\mathfrak{g}}} = 1$, 
$\mathfrak{g}$ is of type $A_{1}$ 
with $X_{\mathfrak{g}} = P_{1} \mathbb{C}$ 
up to biholomorphisms. 
}
\\ \ 
(C) 
{\it 
$(P_{e_{2}}, \theta_{\mathfrak{g}})$ 
is a symplectifcation of degree $1$ on $X_{\mathfrak{g}}$ 
with the right action 
$R_{\lambda}$ 
$(\lambda \in \mathbb{C}^{\times})$ 
defined by the restriction 
of the complex scalar multiple 
$\lambda \cdot$ on 
$\mathfrak{g} \backslash \{ 0 \}$ 
and the projection 
$p: P_{e_{2}} \to X_{\mathfrak{g}}; 
\eta \mapsto \pi_{\mathfrak{g}}(\eta)$ 
defined by the canonical projetion 
$\pi_{\mathfrak{g}}: 
\mathfrak{g} \backslash \{ 0 \} 
\to P_{\mathfrak{g}} \mathbb{C}$, 
so that 
$(X_{\mathfrak{g}}, E_{\theta_{\mathfrak{g}}})$ 
is a c-manifold. 
}
\\ \
(D) 
{\it 
$L_{E_{\theta_{\mathfrak{g}}}}$ 
is very ample, 
so that $X_{\mathfrak{g}}$ is Fano 
and simply connected. 
Moreover, 
$X_{\mathfrak{g}}$ is the projective orbit of 
every highest root vector of 
$\mathfrak{g}$ by $G^{\circ} = {\rm Ad}(\mathfrak{g})$ 
in $P_{\mathfrak{g}} \mathbb{C}$, 
which does not depend on the choice of $\mathbb{Z}$-gradation of c-type 
for $\mathfrak{g}$. 
}
\\ \ 
(E) 
{\it 
$\delta_{\mathfrak{g}}: 
G = {\rm Aut}(\mathfrak{g}) 
\to {\rm Aut}(P_{e_{2}}, 
\theta_{\mathfrak{g}}); 
g \mapsto g |_{P_{e_{2}}}$ 
is well-defined as an injective group homomorphism, 
which is also an immersion by \cite{Cc1946}. 
}
\smallskip \\ \
{\it Proof}. 
(A) 
By Lemma 5.1 (A),  
$\Phi_{e_{2}}(G/G_{e_{2}}) = P_{e_{2}} 
= \mathfrak{W}_{c_{\mathfrak{g}}}$, 
which is a closed subset of a locally compact Hausdorff space 
$\mathfrak{g} \backslash \{ 0 \}$, 
so that 
$P_{e_{2}}$ is a locally compact Hausdorff space. 
By \cite[p.121, Thm.II.3.2]{Hs1978}, 
$\Phi_{e_{2}}$ 
is a homeomorphism onto $P_{e_{2}}$, 
so that 
$P_{e_{2}}$ 
is a complex regular submanifold of $\mathfrak{g}$ 
biholomorphic to $G/G_{e_{2}}$ by $\Phi_{e_{2}}$. 
Then there exists unique holomorphic (1,0)-form 
$\theta_{\mathfrak{g}}$ 
on $P_{e_{2}}$ 
such that 
$\breve{\theta} 
= \Phi_{e_{2}}^{*} \theta_{\mathfrak{g}}$. 
Similarly, one obtains the second assertion. 
\\ \
(B) 
Put 
${\it Lie}~G_{+} := \{ X \in \mathfrak{g} {;~} 
e^{ t {\rm ad} X} \in G_{+} (t \in \mathbb{R}) \} 
= \{ X \in \mathfrak{g} {;~} 
[e_{2}, X] \in \mathfrak{g}_{2} \}$.   
By Lemma 5.0 (B, C), 
${\rm ad}_{-1}(e_{2})$ is an isomorphism 
and $[e_{2}, e_{-2}] \not\in \mathfrak{g}_{2}$. 
Then 
${\it Lie}~G_{+} 
= \oplus_{j = 0}^{2} \mathfrak{g}_{j}$, 
so that 
${\rm n}_{X_{\mathfrak{g}}} = 
{\rm n}_{\oplus_{j = -2}^{-1} \mathfrak{g}_{j}} 
= 1 + {\rm n}_{\mathfrak{g}_{-1}}$. 
Hence, 
${\rm n}_{X_{\mathfrak{g}}} = 1$ 
if and only if 
$\mathfrak{g}_{-1} = \{ 0 \}$ 
if and only if 
$\mathfrak{g} = \hat{\mathfrak{g}} \cong sl_{2} \mathbb{C}$ 
if and only if 
${\rm rank}(\mathfrak{g}) = 1$, 
so that 
$X_{\mathfrak{g}} = P_{1} \mathbb{C}$ 
by \cite{KOY1999}. 
\\ \
(C) 
By Lemma 5.0 (A), 
$D_{\theta_{G}} 
= \cup_{\alpha \in G} \{ X_{R} |_{\alpha} {;~} 
X \in \oplus_{j = -1}^{2} \mathfrak{g}_{j} \}$ 
and 
\begin{eqnarray*}
(d \theta_{G} |_{D_{\theta_{G}}})^{\perp} 
&=& \cup_{\alpha \in G} \{ X_{R} |_{\alpha} {;~} 
X \in \oplus_{j = -1}^{2} \mathfrak{g}_{j}, 
B_{\mathfrak{g}} ([X, \oplus_{j = -1}^{2} \mathfrak{g}_{j}], e_{2}) = 0 \} 
\\
&=& \cup_{\alpha \in G} \{ X_{R} |_{\alpha} {;~} 
X \in \oplus_{j = -1}^{2} \mathfrak{g}_{j}, 
B_{\mathfrak{g}} 
(X, \oplus_{j = 1}^{2} \mathfrak{g}_{j}) = 0 \} 
\\
&=& \cup_{\alpha \in G} \{ X_{R} |_{\alpha} {;~} 
X \in \oplus_{j = 0}^{2} \mathfrak{g}_{j} \}. 
\end{eqnarray*}
For $g \in G$, 
$L_{g}^{*} \theta_{G} = \theta_{G}$ 
and 
$(R_{g}^{*} \theta_{G}) X_{R} 
= (A_{g^{-1}}^{*} \theta_{G}) X_{R} 
= B_{\mathfrak{g}} (X, g (e_{2}))$, 
so that 
$R_{\beta}^{*} \theta_{G} = \theta_{G}$ 
($\beta \in G_{e_{2}}$). 
For $g \in G$, put 
$\breve{g}: 
G/G_{e_{2}} \to G/G_{e_{2}}; 
\pi_{e_{2}}(\alpha) 
\mapsto \pi_{e_{2}}(L_{g} (\alpha))$ 
such that 
$\pi_{e_{2}} \circ L_{g} 
= \breve{g} \circ \pi_{e_{2}}$, 
which is biholomorphic by \cite[Thm.II.4.2]{Hs1978}, 
so that 
$\breve{g}^{*} \breve{\theta} = \breve{\theta}$ 
and 
$g \circ \Phi_{e_{2}} 
= \Phi_{e_{2}} \circ \breve{g}$. 
By Lemma 5.0 (A), 
$\phi: G_{+}/ G_{e_{2}} \to \mathbb{C}^{\times}; 
\beta G_{e_{2}} \mapsto 
B_{\mathfrak{g}} (e_{-2}, \beta (e_{2}))$ 
is a complex Lie group isomorphism. 
Any 
$\lambda \in \mathbb{C}^{\times}$ 
admits unique 
$g_{\lambda} \in G_{+}$ 
such that 
$\phi(\pi_{e_{2}}(g_{\lambda})) = \lambda$, 
so that 
$g_{\lambda}(e_{2}) = \lambda \cdot e_{2}$ 
and 
$R_{g_{\lambda}}^{*} \theta_{G} 
= \lambda \cdot \theta_{G}$. 
For $\lambda \in \mathbb{C}^{\times}$, 
put 
$\breve{R}_{\lambda}: 
G/G_{e_{2}} \to G/G_{e_{2}}; 
\pi_{e_{2}}(g) \mapsto 
\pi_{e_{2}} (R_{g_{\lambda}}(g))$ 
such as 
$\breve{R}_{\lambda} \circ \pi_{e_{2}} 
= \pi_{e_{2}} \circ R_{g_{\lambda}}$, 
by which $G/G_{e_{2}}$ 
is a holomorphic principal 
$\mathbb{C}^{\times}$-bundle on $G/G_{+}$ 
with the projection $\breve{p}$.  
Because of 
$\pi_{e_{2}}^{*} 
(\breve{R}_{\lambda}^{*} \breve{\theta}) 
= R_{g_{\lambda}}^{*} 
(\pi_{e_{2}}^{*} \breve{\theta}) 
= R_{g_{\lambda}}^{*} \theta_{G} 
= \lambda \cdot \theta_{G} 
= \pi_{e_{2}}^{*} 
(\lambda \cdot \breve{\theta})$, 
$\breve{R}_{\lambda}^{*} \breve{\theta} 
= \lambda \cdot \breve{\theta}$, 
i.e., 
$(G/G_{e_{2}}, \breve{\theta})$ 
{\it is a pre-symplectification 
of degree $1$ on} $G/G_{+}$. 
Then 
$P_{e_{2}}$ 
is a holomorphic principl $\mathbb{C}^{\times}$-bundle 
on $X_{\mathfrak{g}}$ 
such as 
$p \circ \Phi_{e_{2}} 
= \varphi_{e_{2}} \circ \breve{p}$: 
In fact, 
for $\lambda \in \mathbb{C}^{\times}$ 
and $\alpha \in G$, 
$R_{\lambda} (\Phi_{e_{2}}(\pi_{e_{2}}(\alpha))) 
= \lambda \cdot (\alpha (e_{2})) 
= \alpha (\lambda \cdot e_{2}) 
= (\alpha \circ g_{\lambda}) e_{2} 
= \Phi_{e_{2}} 
(\pi_{e_{2}} (R_{g_{\lambda}} \alpha)) 
= \Phi_{e_{2}} (\breve{R}_{\lambda} 
(\pi_{e_{2}}(\alpha)))$, 
i.e.,  
$R_{\lambda} \circ \Phi_{e_{2}} 
= \Phi_{e_{2}} \circ \breve{R}_{\lambda}$. 
Moreover, 
$(P_{e_{2}}, \theta_{\mathfrak{g}})$ 
is a pre-symplectification of degree $1$ 
on $X_{\mathfrak{g}}$. 
Because of  
$\theta = \pi_{e_{2}}^{*} (\Phi_{e_{2}}^{*} \theta_{\mathfrak{g}})$, 
one has that 
$p \circ \Phi_{e_{2}} \circ \pi_{e_{2}} 
= \varphi_{e_{2}} \circ \pi_{+}$, 
$p_{*'} (d \theta_{\mathfrak{g}} |_{D_{\theta_{\mathfrak{g}}}})^{\perp} 
= p_{*'} (\Phi_{e_{2}*}(\pi_{e_{2}*} 
(d \theta_{G} |_{D_{\theta_{G}}})^{\perp})) 
= \cup_{\alpha \in G} 
\{ \varphi_{e_{2}*}(\pi_{+*} (X_{R} |_{\alpha})) 
{;~} X \in \oplus_{j = 0}^{2} \mathfrak{g}_{j} 
\} = 0$. 
By Lemma 4.0 (A) (d), 
$(P_{e_{2}}, \theta_{\mathfrak{g}})$ 
is a symplectification of degree $1$ 
on $X_{\mathfrak{g}}$, so that 
$(X_{\mathfrak{g}}, E_{\theta_{\mathfrak{g}}})$ 
is a c-manifold. 
\\ \
(D) 
By Proposition 4.3 (B), 
there exists a biholomorphic map 
$\kappa_{\theta_{\mathfrak{g}}}: 
P_{e_{2}} \to P_{E_{\theta_{\mathfrak{g}}}}$ 
such that 
$\kappa_{\theta_{\mathfrak{g}}} \circ R_{\lambda} 
= R_{\lambda} \circ 
\kappa_{\theta_{\mathfrak{g}}}$ 
($\lambda \in \mathbb{C}^{\times}$). 
Let $\{ \mathfrak{e}_{i}{;~} i \in \{ 0, \ldots, {\rm n}_{\mathfrak{g}} - 1 \} \}$ 
be a $\mathbb{C}$-linear base of $\mathfrak{g}$. 
For $j \in \{ 0, \ldots, {\rm n}_{\mathfrak{g}} - 1 \}$, 
put $z_{j}: \mathfrak{g} \to \mathbb{C}; 
\sum_{i} z_{i} \mathfrak{e}_{i} \mapsto z_{j}$, 
so that 
$z_{j} |_{P_{e_{2}}} \in {\cal O}^{1}(P_{e_{2}})$ 
and 
$f_{j} 
:= z_j \circ \kappa_{\theta_{\mathfrak{g}}}^{-1} 
\in {\cal O}^1(P_{E_{\theta_{\mathfrak{g}}}})$. 
By Lemma 3.1 (A) 
with $m = 1$, 
there exists 
$t_{j} \in \Gamma_{M}(L_{E_{\theta_{\mathfrak{g}}}})$ 
such that 
$\iota_{L_{E_{\theta_{\mathfrak{g}}}}}(t_{j}) 
= f_{j}$.  
Since 
$(f_{0}, \ldots, f_{{\rm n}_{\mathfrak{g}}-1}) 
= (z_{0}, \ldots, z_{{\rm n}_{\mathfrak{g}}-1}) 
\circ \kappa_{\theta_{\mathfrak{g}}}^{-1}: 
P_{L_{E_{\theta_{\mathfrak{g}}}}} \to \mathbb{C}_{{\rm n}_{\mathfrak{g}}}$ 
is a holomorphic injective immersion, 
one has that 
$L_{E_{\theta_{\mathfrak{g}}}}$ is very ample. 
By Lemma 4.6, the connected $X_{\mathfrak{g}}$ 
is a Fano manifold which is simply connected. 
Take a $\mathbb{Z}$-gradation 
$\mathfrak{g} = \oplus_{i \in \mathbb{Z}} \mathfrak{g}_{i}$ 
of c-type such that 
$\mathfrak{g}_{2} = \mathbb{C} e_{2} = \mathbb{C} e_{\rho}$ 
with a highest root vector $e_{\rho}$ of $\mathfrak{g}$. 
By Lemma 5.1 (A, B), 
$P_{e_{\rho}} = P_{e_{2}}$ 
so that 
$p_{\mathfrak{g}}(P_{e_{\rho}}) = p_{\mathfrak{g}}(P_{e_{2}})$, 
as required. 
\\ \ 
(E) 
For $g \in G$, 
$\delta_{\mathfrak{g}} (g) \circ \Phi_{e_{2}} 
= \Phi_{e_{2}} \circ \breve{g}$ 
and 
$\Phi_{e_{2}}^{*} 
(\delta_{\mathfrak{g}}(g)^{*} \theta_{\mathfrak{g}}) 
= \breve{g}^{*} (\Phi_{e_{2}}^{*} \theta_{\mathfrak{g}}) 
= \breve{g}^{*} \breve{\theta} 
= \breve{\theta} 
= \Phi_{e_{2}}^{*} \theta_{\mathfrak{g}}$, 
so that 
$\delta_{\mathfrak{g}}(g)^{*} \theta_{\mathfrak{g}} 
= \theta_{\mathfrak{g}}$. 
Then 
$\delta_{\mathfrak{g}}$ 
is well-defined as a complex Lie group homomorphism. 
If 
$\delta_{\mathfrak{g}}(g) = {\rm id}_{P_{e_{2}}}$, 
then 
$g |_{P_{e_{2}}} = {\rm id}_{P_{e_{2}}}$, 
so that 
$g = {\rm id}_{\mathfrak{g}}$ 
because of 
$\mathfrak{g} 
= \{ \sum_{j} v_{j} {;~} v_{j} \in P_{e_{2}} \}$ 
by Lemma 5.1 (A). 
Also for $X \in \mathfrak{g}$ with $[ X, Y ] = 0$ ($Y \in P_{e_{2}}$), 
$[ X, Y ] = 0$ ($Y \in \mathfrak{g}$), 
so that $\mathbb{C} X = \{ 0 \}$ as an ideal of a simple $\mathfrak{g}$.
Hence, $\delta_{\mathfrak{g}}$ is an injective Lie group homomorphism, 
which is an immersion (cf. \cite[Lem.IV.9.1, Cor.1]{Cc1946}). 
~\qed 
\smallskip \\ \
\textsc{Lemma 5.3}. 
{\it 
Let $\mathfrak{g}$ and $\mathfrak{g}'$ be isomorphic simple complex Lie algebras 
of ${\rm n}_{\mathfrak{g}} = {\rm n}_{\mathfrak{g}'} > 1$. 
Then 
$(X_{\mathfrak{g}}, E_{\theta_{\mathfrak{g}}}) 
\cong (X_{\mathfrak{g}'}, E_{\theta_{\mathfrak{g}'}})$.
} 
\smallskip \\ \
{\it Proof}. 
Let $f: \mathfrak{g} \to \mathfrak{g}$ be the isomorphism. 
Let 
$\mathfrak{g} = \sum_{i \in \mathbb{Z}} \mathfrak{g}_{i}$ 
be a $\mathbb{Z}$-gradation of c-type with $0 \ne e_{2} \in \mathfrak{g}_{2}$, 
so is $\mathfrak{g}' = \sum_{i \in \mathbb{Z}} \mathfrak{g}_{i}'$ 
with $0 \ne e_{2}' \in \mathfrak{g}_{2}'$ 
for $\mathfrak{g}_{i}' := f(\mathfrak{g}_{i})~(i \in \mathbb{Z})$ 
and $e_{2}' := f(e_{2})$. 
Put 
$\widetilde{f}: G := {\rm Aut}(\mathfrak{g}) \to G' := {\rm Aut}(\mathfrak{g}'); 
\beta \mapsto f \beta f^{-1}$, 
which is a Lie group isomorphism such that 
$\widetilde{f}(G_{e_{2}}) = G_{e_{2}'}$. 
Put 
$\breve{f}: G/G_{e_{2}} \to G'/G_{e_{2}'}'; 
\beta G_{e_{2}} \to \widetilde{f}(\beta) G_{e_{2}'}'$, 
so that 
$\pi_{e_{2}'} \circ \widetilde{f} 
= \breve{f} \circ \pi_{e_{2}'}$ 
and 
$\Phi_{e_{2}'} \circ \breve{f} 
= f \circ \Phi_{e_{2}}$. 
Then $f(P_{e_{2}}) = \widetilde{f}(G) f(e_{2}) 
= P_{e_{2}'}$. 
Because of 
$\widetilde{f}_{*} X_{R} = f(X)_{R}$ 
($X \in \mathfrak{g}$) 
and $f^{*} B_{\mathfrak{g}'} = B_{\mathfrak{g}}$, 
one has that 
$\widetilde{f}^{*} \theta_{G'} = \theta_{G}$. 
By Proposition 5.2 (A), 
$\pi_{e_{2}}^{*} (\Phi_{e_{2}}^{*} \theta_{\mathfrak{g}}) 
= \theta_{G} 
= \widetilde{f}^{*} \theta_{G'} 
= \widetilde{f}^{*} (\pi_{e_{2}'}^{*} (\Phi_{e_{2}'}^{*} \theta_{\mathfrak{g}'})) 
= \pi_{e_{2}}^{*} (\breve{f}^{*} (\Phi_{e_{2}'}^{*} \theta_{\mathfrak{g}'})) 
= \pi_{e_{2}}^{*} (\Phi_{e_{2}}^{*} (f^{*} \theta_{\mathfrak{g}'}))$, 
so that 
$\theta_{\mathfrak{g}} = f^{*} \theta_{\mathfrak{g}'}$. 
Then 
$f_{*} D_{\theta_{\mathfrak{g}}} = D_{\theta_{\mathfrak{g}'}}$. 
Put 
$[f]: P_{\mathfrak{g}} \mathbb{C} \to P_{\mathfrak{g}'} \mathbb{C}: 
p_{\mathfrak{g}}(X) \to p_{\mathfrak{g}'}(f(X))$, 
which is a biholomorphism such that 
$[f] X_{\mathfrak{g}} = X_{\mathfrak{g}'}$ 
and 
$[f]_{*} E_{\theta_{\mathfrak{g}}} = E_{\theta_{\mathfrak{g}'}}$. 
\qed 
\smallskip \\ \
{\bf 6. 
Proof of Theorems 1.1 and 1.2.} 
\smallskip \\ \
\textsc{Lemma 6.0.} 
{\it 
Let $(M, E)$ be a c-manifold with the c-line bundle $L_{E}$ 
and the canonical symplectification $(P_{E}, \theta_{E})$ 
in Lemma 4.0 (B). 
Then 
$L_{E}$ is immersive if and only if 
$T_{\eta}' P_{E} 
= \{ X |_{\eta}{;~} X \in {\cal G} \}$ 
$(\eta \in P_{E})$ 
for some finite-dimensional $\mathbb{C}$-linear subspace 
${\cal G}$ of $a(P_{E}, \theta_{E})$. 
In this case, there exists 
a finite-dimensional $\mathbb{C}$-linear subspace 
${\cal A}$ of $a(M, E)$ 
such that 
at each $\xi \in M$ 
\[
T_{\xi}' M =\{ Y |_{\xi}{;~} Y \in {\cal A} \}. 
\]
}
\
{\it Proof}. 
By Lemma 3.1 (B), 
$L_{E}$ is immersive 
if and only if there exist $N \in \mathbb{N}$ 
and $\mathbb{C}$-linearly independent 
$f_{j} \in {\cal O}^{1}(P_{E})$ ($j \in \{ 0, \ldots, N \}$) 
such that $\Phi := (f_{0}, \ldots, f_{N}): P_{E} \to \mathbb{C}_{N+1}$ 
is an immersion. 
Put 
${\cal G} := \{ \sum_{j = 0}^{N} c_{j} H^{\theta_{E}}_{f_{j}} 
{;}~ c_{j} \in \mathbb{C} \}$. 
Then 
${\rm n}_{\cal G} = N < \infty$ 
and 
${\cal G} \subseteq a(P_{E}, \theta_{E})$ 
because of Lemma 4.1 (B) for 
$\delta = \epsilon_{\theta} = \epsilon_{\theta_{E}}= 1$. 
By natural complex coordinates 
$(z_{0}, \ldots, z_{N})$ of $\mathbb{C}_{N+1}$, 
$f_{j} = \Phi^{*} z_{j}$ ($j \in \{ 0, \ldots, N \}$). 
For $\eta \in P_{E}$, 
$\Phi_{*\eta}: T_{\eta}' P_{E} 
\to T_{\Phi(\eta)}' \mathbb{C}_{N+1}$ 
is injective if and only if the dual map 
$\Phi_{*\eta}^{*}: 
{T'}_{\Phi(\eta)}^{*} \mathbb{C}_{N+1} \to {T'}_{\eta}^{*} P_{E}$
is surjective, that is, 
$T_{\eta}' P_{E} 
= \natural (\Phi_{*\eta}^{*} {T'}_{\Phi(\eta)}^{*} \mathbb{C}_{N+1}) 
= \{ (\natural \circ \Phi_{*\eta}^{*}) \sum_{j = 0}^{N} c_{j} d z_{j}
{;~} c_{j} \in \mathbb{C} \} 
= \{ \sum_{j = 0}^{N} c_{j} H^{\theta_{E}}_{f_{j}} |_{\eta} 
{;~} c_{j} \in \mathbb{C} \} 
= \{ X |_{\eta} {;~} X \in {\cal G} \}$ 
by 
$H^{\theta_{E}}_{f_{j}} 
= \natural (d f_{j}) = \natural ((\Phi_{*})^{*} d z_{j})$, 
as required. 
By Lemma 4.1 (E) (e), put 
${\cal A} := \widetilde{p}_{E*} {\cal G} \subseteq a(M, E)$ 
of 
${\rm n}_{\cal A} = N < \infty$. 
For $\xi = p(\eta)$, 
$T_{\xi}' M = p_{*} T_{\eta}' P_{E} 
= p_{*} \{ X |_{\eta} {;~} X \in {\cal G} \} 
= \{ Y |_{\xi}{;~} Y \in {\cal A} \}$. 
\qed 
\smallskip \\ \
Let $(M, E)$ be a c-manifold 
with the canonical symplectification 
$(P_{E}, \theta_{E})$ of Lemma 4.0 (B). 
For any $\mathbb{C}$-linear subspace ${\cal G}$ of $a(P_{E}, \theta_{E})$ 
with the dual space ${\cal G}^{*}$ and the map {\it ev}: 
${\cal G}^{*} \times {\cal G} \to \mathbb{C}; 
(\omega, X) \mapsto <\! \omega, X \!>_{\cal G} 
:= \omega (X)$, put {\it a ${\cal G}$-moment map} 
as follows (cf. \cite{Sa2001}): 
\[
\mu_{{\cal G}}: P_{E} \to {\cal G}^{*}; 
\eta \mapsto \mu_{{\cal G}}(\eta); 
<\!\! \mu_{{\cal G}}(\eta), X \!\!>_{\cal G} 
= \theta_{E} (X |_{\eta})~~(X \in {\cal G}). 
\] 
For $\beta \in {\rm Aut}(P_{E})$, 
$d \beta: 
\Gamma_{P_{E}}(T' P_{E}) 
\to \Gamma_{P_{E}}(T' P_{E}); X \mapsto (d \beta) X 
:= \beta_{*} (X |_{\beta^{-1}})$ 
is bijective. 
For 
$\beta \in {\rm Aut}(P_{E}, \theta_{E})$ 
such as 
$(d \beta) {\cal G} = {\cal G}$, 
the restriction of $d \beta$ to ${\cal G}$ 
is denoted as 
$d_{\cal G} \beta: {\cal G} \to {\cal G}; 
X \mapsto (d \beta) X$, 
which is bijective. 
Then a map 
\[
d_{\cal G}^{*} \beta: {\cal G}^{*} \to {\cal G}^{*}; 
\omega \mapsto (d_{\cal G}^{*} \beta) \omega; 
<\!\! (d_{\cal G}^{*} \beta) \omega, X \!\!>_{\cal G} 
= <\!\! \omega, (d \beta)^{-1} X \!\!>_{\cal G}
\]
is well-defined. 
For a $\mathbb{C}$-bilinear form $B$ on ${\cal G}$, 
a $\mathbb{C}$-linear map is defined as 
\[
\flat_{B}: 
{\cal G} \to {\cal G}^{*}; 
X \mapsto \flat_{B} X; 
<\!\!\flat_{B} X, Y\!\!>_{\cal G} 
:= B(X,~Y)~(Y \in {\cal G}). 
\]
Assume that $\flat_{B}$ is bijective, 
i.e., 
$B^{\perp} = \{ 0 \}$ 
for 
$B^{\perp} := 
\{ X \in {\cal G};~B(X, Y) = 0~(Y \in {\cal G}) \}$. 
Put 
$\mu^{\natural}_{B} := \flat_{B}^{-1} \circ \mu_{{\cal G}}$, 
{\it the $B$-dual moment map}. 
\smallskip \\ \
\textsc{Lemma 6.1}. 
{\it 
Let 
$\beta \in {\rm Aut}(P_{E}, \theta_{E})$ 
be such that $(d \beta) {\cal G} = {\cal G}$ 
and $(d \beta)^{*} B = B$. 
Then one has the following three results. 
}
\\ \
(A) 
{\it 
$\mu_{{\cal G}} \circ R_{\lambda} 
= \lambda \cdot \mu_{{\cal G}}$ 
and 
$\mu^{\natural}_{B} \circ R_{\lambda} 
= \lambda \cdot \mu^{\natural}_{B}$ 
for any $\lambda \in \mathbb{C}^{\times}$. 
}
\\ \
(B) 
{\it 
$\mu_{{\cal G}} \circ \beta 
= d_{\cal G}^{*} \beta \circ \mu_{{\cal G}}$, 
$d^{*}_{\cal G} \beta \circ \flat_{B} 
= \flat_{B} \circ d_{\cal G} \beta$ 
and 
$\mu^{\natural}_{B} \circ \beta 
=  d_{\cal G} \beta \circ \mu^{\natural}_{B.}$ 
}
\\ \ 
(C) 
{\it 
Assume that 
${\rm n}_{{\cal G}} = N + 1 < \infty$ 
and 
$\{ H^{\theta_{E}}_{f_{j}}{;~} j \in \{ 0, \ldots, N \} \}$ 
is any $\mathbb{C}$-linear basis of ${\cal G}$ 
with some 
$\{ f_{j}{;~} j \in \{ 0, \ldots, N \}\} 
\subseteq {\cal O}^{1}(P_{E})$ 
and the dual basis 
$\{ \omega_{j}{;~} j \in \{ 0, \ldots, N \} \}$ 
of ${\cal G}^{*}$. 
Then 
$\mu_{{\cal G}} 
= \sum_{j = 0}^{N} f_{j} \cdot \omega_{j}$. 
} 
\smallskip \\ \
{\it Proof}. 
(A) 
$< \!\! \mu_{{\cal G}}(R_{\lambda} \eta), X \!\! >_{{\cal G}} 
= \theta_{E} (X |_{R_{\lambda} \eta}) 
= \theta_{E} ( R_{\lambda*} (X |_{\eta})) 
= \lambda \cdot \theta_{E} (X |_{\eta}) 
= \lambda \cdot <\!\! \mu_{{\cal G}}(\eta), X \!\!>_{{\cal G}} 
= <\!\! \lambda \cdot \mu_{{\cal G}}(\eta), X \!\!>_{{\cal G}}$ 
for $(\eta, \lambda) \in P_{E} \times \mathbb{C}^{\times}$. 
Hence, 
$\mu_{{\cal G}} \circ R_{\lambda} 
= \lambda \cdot \mu_{{\cal G}}$ 
and 
$\mu^{\natural}_{B} \circ R_{\lambda} 
= \flat_{B}^{-1} \circ \mu_{{\cal G}} \circ R_{\lambda} 
= \flat_{B}^{-1} \circ (\lambda \cdot \mu_{{\cal G}}) 
= \lambda \cdot \flat_{B}^{-1} \circ \mu_{{\cal G}} 
= \lambda \cdot \mu^{\natural}_{B}$. 
\\ \
(B) 
Because of 
$\beta^{*} \theta_{E} = \theta_{E}$, 
for $(\eta, X) \in P_{E} \times {\cal G}$, 
\begin{eqnarray*}
& &<\!\! \mu_{{\cal G}} (\beta \eta),~X \!\!>_{\cal G} 
= \theta_{E} (X |_{\beta \eta}) 
= \theta_{E} (\beta_{*} \beta_{*}^{-1} (X |_{\beta \eta})) 
= \theta_{E} (\beta_{*}^{-1} (X |_{\beta \eta})) 
\\
&=& \theta_{E} ( d \beta^{-1} X ) |_{\eta} 
= <\!\! \mu_{{\cal G}}(\eta),~(d \beta)^{-1} X \!\!>_{\cal G}  
= <\!\! (d^{*}_{{\cal G}} \beta) \mu_{{\cal G}}(\eta),~X \!\!>_{\cal G}. 
\end{eqnarray*}
Hence, 
$\mu_{{\cal G}} \circ \beta 
= d^{*}_{{\cal G}} \beta \circ \mu_{{\cal G}}$. 
Because of 
$(d \beta)^{*} B = B$, 
for $X, Y \in {\cal G}$, 
\begin{eqnarray*}
<\!\! d^{*}_{{\cal G}} \beta (\flat_{B} X),~Y \!\!>_{\cal G} 
&=& <\!\! \flat_{B} X,~(d \beta)^{-1} Y \!\!>_{\cal G} 
= B(X,~(d \beta)^{-1} Y) 
\\
&=& B((d \beta) X,~Y) 
= <\!\! \flat_{B}((d \beta) X),~Y \!\!>_{\cal G}, 
\end{eqnarray*}
i.e.,  
$d^{*}_{{\cal G}} \beta \circ \flat_{B} 
= \flat_{B} \circ d_{\cal G} \beta$. 
Then 
$\mu^{\natural}_{B} \circ \beta 
= \flat_{B}^{-1} \circ \mu_{{\cal G}} \circ \beta 
= \flat_{B}^{-1} \circ d^{*}_{{\cal G}} \beta \circ \mu_{{\cal G}} 
= d_{\cal G} \beta \circ \flat_{B}^{-1} \circ \mu_{{\cal G}} 
= d_{\cal G} \beta \circ \mu^{\natural}_{B}$. 
\\ \
(C) 
For $\eta \in P_{E}$, 
$\mu_{{\cal G}}(\eta) \in {\cal G}^{*} 
= \oplus_{i = 1}^{N} \mathbb{C} \omega_{i}$, 
so that 
$\mu_{{\cal G}}(\eta)  
= \sum_{j = 0}^{N} g_{j}(\eta) \cdot \omega_{j}$ 
for some $g_{j}(\eta) \in \mathbb{C}$. 
Because of Lemma 4.1 (B) 
for $\delta = \epsilon_{\theta} = \epsilon_{\theta_{E}} = 1$, 
\[
f_{j}(\eta) 
= \theta_{E} (H^{\theta_{E}}_{f_{j}} |_{\eta}) 
= <\!\! \mu_{{\cal G}}(\eta), H^{\theta_{E}}_{f_{j}} \!\!>_{\cal G} 
= g_{j}(\eta). 
\]
Hence, 
$\mu_{{\cal G}}(\eta)  
= \sum_{j = 0}^{N} f_{j}(\eta) \cdot \omega_{j}$. 
~\qed 
\smallskip \\ \
\textsc{Lemma 6.2}. 
{\it 
$L_{E}$ is immersive (resp. very ample) 
if and only if 
$\mu_{{\cal G}}$ is an immersion (resp. an injective immersion) 
for some finite-dimensional $\mathbb{C}$-linear subspace 
${\cal G}$ of $a(P_{E}, \theta_{E})$. 
This is the case for ${\cal G} := a(P_{E}, \theta_{E})$ 
if ${\rm n}_{a(P_{E}, \theta_{E})} < \infty$. 
} 
\vspace{1mm}
\\ \
{\it Proof}. 
By Lemma 3.1 (B), 
$L_{E}$ is immersive (resp. very ample) 
if and only if there exist $N \in \mathbb{N}$ 
and $\mathbb{C}$-linearly independent 
$f_{j} \in {\cal O}^{1}(P_{E})$ ($j \in \{ 0, \ldots, N \}$) 
such that 
$\Phi := (f_{0}, \ldots, f_{N}): P_{E} \to \mathbb{C}_{N+1}$ 
is an immersion (resp. injective immersion). 
Then 
$\{ f_{0}, \ldots, f_{N} \}$ 
is a $\mathbb{C}$-linear basis of 
${\cal G} := \{ \sum_{j = 0}^{N} c_{j} H^{\theta_{E}}_{f_{j}} 
{;}~ c_{j} \in \mathbb{C} \}$. 
Let 
$\{ \omega_{0}, \ldots, \omega_{N} \}$ 
be the dual basis of ${\cal G}^{*}$. 
By Lemma 6.1 (C), 
$\mu_{{\cal G}} = \sum_{k = 0}^{N} f_{k} \cdot \omega_{k}$, 
which is an immersion (resp. injection) if and only if 
$\Phi: P_{E} \to \mathbb{C}_{N+1}$ 
is an immersion (resp. injection), 
by which the first assertion follows. 
The last assertion follows from the first one.~\qed
\smallskip \\ \
\textsc{Proposition 6.3}. 
{\it 
Let $(M, E)$ be a connected compact c-manifold of ${\rm n}_{M} > 0$. 
Assume that ${\cal G}$ is a complex Lie subalgebra of $a(P_{E}, \theta_{E})$ 
such that 
$T_{\eta}' P_{E} = \{ X |_{\eta}; X \in {\cal G} \}$ 
at each $\eta \in P_{E}$. 
Then one has the following results. 
}
\smallskip \\ \
(A) 
{\it 
$M$ is a simply connected k\"ahlerian manifold. 
${\cal G}$ is a semisimple complex Lie algebra of ${\rm n}_{\cal G} > 1$, 
so that $P_{\cal G} \mathbb{C}$ is well-defined 
and that $\mu^{\natural}_{B_{\cal G}}: P_{E} \to {\cal G}$ 
is well-defined with respect to the Killing form $B_{\cal G}$ of ${\cal G}$ 
as a holomorphic immersion. 
} 
\smallskip \\ \
(B) 
{\it 
${\cal G}$ 
is a simple complex Lie algebra of ${\rm n}_{\cal G} > 1$, 
and that $\mu^{\natural}_{B}: P_{E} \to {\cal G}$ 
is an injection such that 
$\mu^{\natural}_{B}(P_{E}) = {\rm Ad}({\cal G}) e_{2}$ 
with respect to a $\mathbb{Z}$-gradation 
${\cal G} 
= \oplus_{i = -2}^{2} {\cal G}_{i}$ 
of c-type such as ${\cal G}_{2} = \mathbb{C} e_{2}$, 
so that the image restriction 
$\hat{\mu}^{\natural}_{B}: 
(P_{E}, \theta_{E}) \to 
({\rm Ad}({\cal G}) e_{2}, \theta_{\cal G}); 
\eta \mapsto \mu^{\natural}_{B}(\eta)$ 
of $\mu^{\natural}_{B}$ 
is an isomorphism, 
which induces an isomorphism of complex Lie groups 
as follows:} 
\[
\widetilde{\mu}^{\natural}_{B}: 
{\rm Aut}(P_{E}, \theta_{E}) \to 
{\rm Aut}({\rm Ad}({\cal G}) e_{2}, \theta_{\cal G}); 
\beta \mapsto \hat{\mu}^{\natural}_{B} 
\circ \beta \circ (\hat{\mu}^{\natural}_{B})^{-1}. 
\]
\
(C) 
{\it 
Both homomorphisms of complex Lie groups defined as 
\begin{eqnarray*}
& &d_{\cal G}: 
{\rm Aut}(P_{E}, \theta_{E}) \to {\rm Aut}({\cal G}); 
\beta \mapsto d_{\cal G} \beta 
:= d \beta |_{\cal G}; {\it and} 
\\
& &
\delta_{\cal G}: {\rm Aut}({\cal G}) 
\to {\rm Aut}({\rm Ad}({\cal G}) e_{2}, 
\theta_{\cal G}); 
g \mapsto g |_{{\rm Ad}({\cal G}) e_{2}} 
\end{eqnarray*}
are isomorphisms such that 
$\widetilde{\mu}^{\natural}_{\cal B} 
= \delta_{\cal G} \circ d_{\cal G}$ 
and 
$(d_{\cal G} \beta) \circ \mu^{\natural}_{B_{\cal G}} 
= \mu^{\natural}_{B_{\cal G}} \circ \beta$ 
for $\beta \in {\rm Aut}(P_{E}, \theta_{E})$. 
}
\smallskip \\ \
(D) 
{\it 
There exists a holomorphic injective immersion 
$\nu_{\cal G}: M \to P_{\cal G} \mathbb{C}$ 
such that 
$\nu_{\cal G} \circ p_{E} 
= \pi_{\cal G} \circ \mu^{\natural}_{\cal B}$, 
$\nu_{\cal G} M = X_{\cal G}$ 
and 
$\nu_{{\cal G}*} E = E_{\theta_{\cal G}}$. 
Moreover, 
$d_{\cal G} \circ {\cal P}: 
{\rm Aut}(M, E) \to {\rm Aut}({\cal G})$ 
is an isomorphism of complex Lie groups such that 
$[(d_{\cal G} \circ {\cal P}) \alpha] (\nu_{\cal G} \xi) 
= \nu_{\cal G}(\alpha \xi)$ 
for 
$(\alpha, \xi) \in {\rm Aut}(M, E) \times M$ 
with respect to the projective action 
$[(d_{\cal G} \circ {\cal P}) \alpha]$ 
of 
$(d_{\cal G} \circ {\cal P}) \alpha \in {\rm Aut}({\cal G})$ 
on $P_{\cal G} \mathbb{C}$. 
}
\smallskip \\ \
{\it Proof}. 
(A0) 
By Lemma 6.0, $L_{E}$ is immersive, 
so that there exists a holomorphic immersion 
$\varphi: M \to P_{N} \mathbb{C}$ 
for some $N \in \mathbb{N}$. 
Then $M$ is a connected compact k\"ahlerian manifold. 
By Lemma 4.3, 
$M$ is simply connected. 
By Lemma 4.5, 
${\rm Aut}(P_{E}, \theta_{E})$ 
is a complex Lie group of 
${\rm n}_{a(P_{E}, \theta_{E})} < \infty$. 
There exists a complex Lie subgroup $G$ 
of ${\rm Aut}(P_{E}, \theta_{E})$ 
corresponding to ${\cal G}$. 
By the assumption 
$T_{\eta}' P_{E} = \{ X |_{\eta}; X \in {\cal G} \}$ 
($\eta \in P_{E}$), 
one has that 
${\rm n}_{\cal G} \geq {\rm n}_{P_{E}} > 1$ 
and that $P_{E} = G \eta$ 
for $\eta \in P_{E}$ since $P_{E}$ is connected. 
\\ \
(A1) 
By Proposition 4.3 (A), 
$\widetilde{\mathfrak{g}} := \widetilde{p_{E*}} ({\cal G}) \subseteq a(M, E)$. 
For $\xi = p_{E}(\eta) \in M$, 
$T_{\xi}' M = p_{E*} T_{\eta}' P_{E} 
= \{ p_{E*\eta} X; X \in {\cal G} \} 
= \{ Y |_{\xi}; Y \in \widetilde{\mathfrak{g}} \}$. 
Let $\widetilde{G}$ be a complex Lie subgroup of 
${\rm Aut}(M, E)$ corresponding to $\widetilde{\mathfrak{g}}$. 
Then $M = \widetilde{G} \xi$ since $M$ is connected. 
\\ \
(A2) 
Let $\widetilde{\mathfrak{s}}$ 
be the solvable radical of 
$\widetilde{\mathfrak{g}}$ 
(cf. \cite[p.22]{GG1978}). 
{\it Claim} 
$\widetilde{\mathfrak{s}} = \{ 0 \}$ 
(cf. \cite[Prop.II, Proof]{Saj1973}): 
In fact, by \cite[Thm.IV.IV.1]{Cc1946}, 
there exists a connected Lie subgroup 
$\widetilde{S}$ of ${\rm Aut}(M, E)$ 
such that 
$\widetilde{\mathfrak{s}} 
= \{ Y \in a(M, E); \alpha_{Y}(t) \in \widetilde{S}~(t \in \mathbb{R}) \}$, 
so that $\widetilde{S}$ 
is a solvable connected complex Lie subgroup of the complex Lie group 
$\widetilde{G}$. 
Since $M$ is simply connected 
(by (A0)), 
there exists $\xi_{0} \in M$ 
such that $\widetilde{S} \xi_{0} = \xi_{0}$ 
by \cite[Prop.I]{Saj1973}. 
Take any $Y \in \widetilde{\mathfrak{s}}$ 
with the real part $Y^{r}$ 
by the notation given in the proof of Lemma 4.5. 
For any $t \in \mathbb{R}$, 
$\alpha_{Y^{r}}(t) \in \widetilde{S}$, 
so that 
$\alpha_{Y^{r}}(t) \xi_{0} = \xi_{0}$. 
For any $\xi_{1} \in M$, 
there exists 
$g_{1} \in \widetilde{G}$ 
such as 
$\xi_{1} = g_{1} \xi_{0}$ 
(by (A1)). 
Then 
$Y_{1} := d g_{1}^{-1} Y 
= {\rm Ad}(g_{1}^{-1}) Y \in \widetilde{\mathfrak{s}}$ 
and 
$\alpha_{Y^{r}}(t) (\xi_{1}) 
= \alpha_{Y^{r}}(t) ( g_{1} \xi_{0} ) 
= g_{1} (\alpha_{Y_{1}^{r}}(t) \xi_{0}) 
= g_{1} \xi_{0} = \xi_{1}$ 
($t \in \mathbb{R}$), 
so that 
$\left. Y^{r} \right|_{\xi_{1}} 
= \left. \frac{d}{dt} \right|_{t = 0} \alpha_{Y^{r}}(t) (\xi_{1}) 
= \left. \frac{d}{dt} \right|_{t = 0} \xi_{1} = 0$. 
Then $Y = (Y^{r})' = 0' = 0$. 
Hence, 
$\widetilde{\mathfrak{s}} = 0$, 
{\it as required}. 
Then 
$\widetilde{\mathfrak{g}}$ 
is semisimple (cf. \cite[(4.3)]{GG1978}), 
so is ${\cal G}$ (by Lemma 4.1 (E) (e)), 
so that $B_{\cal G}^{\perp} = \{ 0 \}$. 
Hence, 
$\mu^{\natural}_{B_{\cal G}} 
= \flat_{B_{\cal G}}^{-1} \circ \mu_{\cal G}$ 
is well-defined as a holomorphic immersion into 
${\cal G} \backslash \{ 0 \}$ 
with ${\rm n}_{\cal G} > 1$ 
by Lemma 6.2. 
\\ \
(D0) 
By Lemma 6.1 (A), a map 
$\nu_{\cal G}: M \to P_{\cal G} \mathbb{C}; 
p_{E}(\eta) \mapsto \pi_{\cal G}(\mu^{\natural}_{B_{\cal G}}(\eta))$ 
is well-defined such that 
$\nu_{\cal G} \circ p_{E} 
= \pi_{\cal G} \circ \mu^{\natural}_{B_{\cal G}}$, 
which is holomorphic. 
Then 
$\nu_{\cal G}$ is holomorphic since 
$p_{E}$ is a holomorphic submersion. 
{\it 
Claim. $\nu_{\cal G}$ is an immersion:} 
Let $X \in T_{\eta}' P_{E}$ be such that 
$\nu_{{\cal G}*} (p_{E *} X) = 0$. 
Then 
$\pi_{{\cal G} *}(\mu^{\natural}_{B_{\cal G} *} X) = 0$, 
so that there exists $c \in \mathbb{C}$ such that 
$\mu^{\natural}_{B_{\cal G} *} X 
= c \cdot \left. \frac{d}{d \lambda} \right|_{\lambda = 1} 
(\lambda \cdot \mu^{\natural}_{B_{\cal G}}(\eta)) 
= \mu^{\natural}_{B_{\cal G} *} 
(c \cdot \left. \frac{d}{d \lambda} \right|_{\lambda = 1} R_{\lambda} \eta)$ 
by Lemma 6.1 (A), so that 
$X = c \cdot \left. \frac{d}{d \lambda} \right|_{\lambda = 1} R_{\lambda} \eta$ 
since $\mu^{\natural}_{B_{\cal G} *\eta}$ is injective (by (A2)). 
Then $p_{E*} X = 0$, {\it as required}. 
Hence, 
$\nu_{\cal G}$ 
is a holomorphic immersion onto the image 
${\cal M} := \nu_{\cal G} (M) 
= \nu_{\cal G} (p_{E}(P_{E})) 
= \pi_{\cal G} (\mu^{\natural}_{B_{\cal G}}(P_{E}))$, 
which is a compact submanifold of $P_{\cal G} \mathbb{C}$ 
since $M$ is compact. 
\\ \ 
(B0) 
$\pi_{\cal G}^{-1}({\cal M}) 
= \pi_{\cal G}^{-1}(\nu_{\cal G} (M)) 
= \mathbb{C}^{\times} \cdot \mu^{\natural}_{B_{\cal G}}(P_{E}) 
= \mu^{\natural}_{B_{\cal G}}(R_{\mathbb{C}^{\times}} P_{E}) 
= \mu^{\natural}_{B_{\cal G}}(P_{E}) =: {\cal W}$ 
in ${\cal G}$. 
For $X \in {\cal G}$, 
$\alpha_{\widetilde{p_{E*}}(X)^{r}}(t) \in {\rm Aut}(M, E)$ 
is well-defined as well as the proof of Lemma 4.5. 
By Lemma 4.4 (B), put 
$\beta_{X}(t) := {\cal P} (\alpha_{\widetilde{p_{E*}}(X)^{r}}(t)) 
\in {\rm Aut}(P_{E}, \theta_{E})$ $(t \in \mathbb{R})$, 
which is the 1-parameter subgroup of ${\rm Aut}(P_{E}, \theta_{E})$ 
such that 
$p_{E} \circ \beta_{X}(t) 
= \alpha_{\widetilde{p_{E*}} (X)^{r}}(t) \circ p_{E}$. 
Then 
$p_{E*} \left. \frac{d}{d t} \right|_{t = 0} \beta_{X} (t) 
= p_{E*} (X^{r})$, 
so that 
$\widetilde{p_{E*}} (\left. \frac{d}{d t} \right|_{t = 0} \beta_{X} (t))' 
= \widetilde{p_{E*}} (X)$. 
Since $\widetilde{p_{E*}}$ is an isomorphism (by Proposition 4.3 (A)), 
$(\left. \frac{d}{d t} \right|_{t = 0} \beta_{X} (t))' = X$. 
For $t \in \mathbb{R}$, 
$d_{\cal G} \beta_{X}(t): {\cal G} \to {\cal G}; 
Y \mapsto (d \beta_{X}(t)) Y 
= \beta_{X}(t)_{*} (Y |_{\beta_{X}(t)^{-1}}) 
= {\rm Ad}(\beta_{X}(t)) Y$, 
so that 
${\cal G} \ni [X, Y] 
= (\left. \frac{d}{dt} \right|_{t=0} (d_{\cal G} \beta_{X}(t)) Y)'$. 
Taking $Y := \mu^{\natural}_{B_{\cal G}} (\eta)$, 
one has by Lemma 6.1 (B) that 
$[X, Y] = (\left. \frac{d}{dt} \right|_{t = 0} 
(d_{\cal G} \beta_{X}(t)) \mu^{\natural}_{B_{\cal G}} (\eta))' 
= \mu^{\natural}_{B_{\cal G}} 
(\left. \frac{d}{dt} \right|_{t = 0} \beta_{X}(t) \eta)'
= \mu^{\natural}_{B_{\cal G} *} (X |_{\eta})$. 
Moreover, 
\[
T_{\mu^{\natural}_{B_{\cal G}} (\eta)}' {\cal W} 
= \mu^{\natural}_{B_{\cal G} *} T_{\eta}' P_{E} 
= \{ \mu^{\natural}_{B_{\cal G} *} (X |_{\eta}){;}~X \in {\cal G} \} 
= \{ [X, \mu^{\natural}_{B_{\cal G}}(\eta)] {;~} X \in {\cal G} \}. 
\]
Since ${\cal W} = \mu^{\natural}_{B_{\cal G}}(P_{E})$ 
is connected, 
${\cal W} = {\rm Ad}({\cal G}) \eta'$ 
for $\eta' \in {\cal W}$. 
Moreover, 
\[
\{ X \in {\cal G} {;~}[X, {\cal W}] = \{ 0 \} \} 
= \{ X \in {\cal G} {;~}\mu^{\natural}_{B_{\cal G} *} (X |_{\eta}) = 0 
~(\eta \in P_{E}) \} = \{ 0 \} 
\] 
since $\mu^{\natural}_{B_{\cal G}}$ is an immersion. 
\\ \ 
(B1) 
Because of (A2), ${\cal G}$ is semisimple. 
By (D0, B0) with Lemma 5.1 (B), 
${\cal G}$ is a simple complex Lie algebra of ${\rm n}_{\cal G} > 1$, 
and that 
${\cal W} = {\rm Ad}({\cal G}) e_{2}$ 
and 
$\nu_{\cal G}(M) = {\cal M} 
= \pi_{\cal G}({\cal W}) = X_{\cal G}$ 
with respect to a $\mathbb{Z}$-gradation of c-type: 
${\cal G} = \oplus_{i = -2}^{2} {\cal G}_{i}$ 
such as ${\cal G}_{2} = \mathbb{C} e_{2}$. 
By Proposition 5.2 (D), 
$X_{\cal G}$ is connected and simply connected. 
By Lemma 2.0 (D), 
$\nu_{\cal G}: M \to X_{\cal G}$ 
is a covering map as a holomorphic immersion 
from a compact complex manifold $M$ onto a complex manifold 
$X_{\cal G}$ because of (D0). 
Since $M$ is connected, 
$\nu_{\cal G}$ is biholomorphic. 
{\it Claim. 
$\mu^{\natural}_{B_{\cal G}}$ is injective:} 
Let 
$\eta_{i} \in P_{E}$ ($i = 1, 2$) 
and 
$\mu^{\natural}_{B_{\cal G}}(\eta_{1}) 
= \mu^{\natural}_{B_{\cal G}}(\eta_{2})$. 
Then 
$\nu_{\cal G}(p_{\cal G}(\eta_{1})) 
= \pi_{\cal G}(\mu^{\natural}_{B_{\cal G}}(\eta_{1})) 
= \pi_{\cal G}(\mu^{\natural}_{B_{\cal G}}(\eta_{2})) 
= \nu_{\cal G}(p_{\cal G}(\eta_{2}))$. 
Since $\nu_{\cal G}$ is injective, 
$p_{\cal G}(\eta_{1}) = p_{\cal G}(\eta_{2})$, 
so that 
$\eta_{1} = \lambda \cdot \eta_{2}$ 
for some $\lambda \in \mathbb{C}^{\times}$. 
Then 
$\lambda \cdot \mu^{\natural}_{B_{\cal G}}(\eta_{2}) 
= \mu^{\natural}_{B_{\cal G}}(\lambda \cdot \eta_{2}) 
= \mu^{\natural}_{B_{\cal G}}(\eta_{2})$. 
Suppose that 
$\mu^{\natural}_{B_{\cal G}}(\eta_{2}) = 0$. 
Then 
$\mu^{\natural}_{B_{\cal G}*\eta_{2}} 
{\cal E}_{P_{E}} 
= \left. \frac{d}{d \lambda} \right|_{\lambda = 1} 
\mu^{\natural}_{{\cal B}} (\lambda \cdot \eta_{2}) 
= \left. \frac{d}{d \lambda} \right|_{\lambda = 1} 
\lambda \cdot \mu^{\natural}_{{\cal B}} (\eta_{2}) 
= 0$, so that ${\cal E}_{P_{E}} = 0$ 
since $\mu^{\natural}_{B_{\cal G}}$ is an immersion, 
{\it a contradiction}. 
Then  
$\mu^{\natural}_{B_{\cal G}}(\eta_{2}) \ne 0$ 
and 
$\lambda = 1$, 
so that $\eta_{1} = \eta_{2}$, 
{\it as required.} 
Hence, 
$\hat{\mu}^{\natural}_{B_{\cal G}}$ 
is biholomorphic. 
\\ \
(B2) 
Put 
$\eta_{2} := (\hat{\mu}^{\natural}_{\cal B})^{-1}(e_{2}) \in P_{E}$. 
For any $X \in {\cal G}$, 
$\theta_{\cal G} (\delta_{{\cal G} *} 
(\mu^{\natural}_{B_{\cal G} *} X |_{\eta_{2}})) 
= \theta_{\cal G}(\delta_{{\cal G} *} [X, e_{2}]) 
= \theta_{\cal G}((\Phi_{e_{2}} \circ \pi_{e_{2}})
_{* {\rm id}_{\cal G}} X_{R}) 
= \theta_{{\rm Ad}({\cal G})} (X_{R}) |_{{\rm id}_{\cal G}} 
= B_{\cal G}(e_{2}, X) 
= B_{\cal G}(\mu^{\natural}_{B_{\cal G}} (\eta_{2}), X)$ 
$= < \!\! \mu_{\cal G}(\eta_{2}), X \!\! >_{\cal G}$ 
$= \theta_{E} (X |_{\eta_{2}})$. 
By 
$T_{\eta_{2}}' P_{E} 
= \{ X |_{\eta_{2}} {;~} X \in {\cal G} \}$, 
$\theta_{\cal G} \circ \delta_{{\cal G} *} 
\circ \mu^{\natural}_{B_{\cal G} * \eta_{2}} 
= \theta_{E} |_{\eta_{2}}$. 
For $\beta \in G (\subseteq {\rm Aut}(P_{E}, \theta_{E}))$, 
$(d \beta^{-1}) {\cal G} = {\cal G}$ 
and $d_{\cal G} \beta^{-1} \in {\rm Aut}({\cal G})$. 
By Proposition 5.2 (E), 
$\delta_{\cal G} (d_{\cal G} \beta^{-1}) 
\in {\rm Aut}({\rm Ad}({\cal G}) e_{2}, \theta_{\cal G})$. 
By Lemma 6.1 (B), 
$\theta_E |_{\beta (\eta_{2})} 
= \theta_E |_{\eta_{2}} \circ 
\beta^{-1}_{* \beta (\eta_{2})}  
= \theta_{\cal G} \circ \delta_{{\cal G} *} 
\circ \mu^{\natural}_{B_{\cal G} * \eta_{2}} 
\circ \beta^{-1}_{* \beta(\eta_{2})} 
= \theta_{\cal G} \circ \delta_{{\cal G} *} 
\circ (d_{\cal G} \beta^{-1})_{*} 
\circ (\mu^{\natural}_{B_{\cal G}})_{* \beta (\eta_{2})}
= \theta_{\cal G} \circ (\delta_{\cal G}(d_{\cal G} \beta^{-1}))_{*} 
\circ (\mu^{\natural}_{B_{\cal G}})_{* \beta (\eta_{2})}
= \theta_{\cal G} \circ (\mu^{\natural}_{B_{\cal G}})_{* \beta (\eta_{2})}$. 
Because of 
$P_{E} = \{ \beta (\eta_{2}) {;~} \beta \in G \}$ 
(by (A0), 
$\theta_{E} = (\mu^{\natural}_{B_{\cal G}})^{*} 
\theta_{\cal G}$. 
Hence, 
$\hat{\mu}^{\natural}_{B_{\cal G}}$ is an isomorphism, 
and that 
$\widetilde{\mu}^{\natural}_{B_{\cal G}}$ 
is a complex Lie group isomorphism. 
\\ \
(D1) 
$E_{\theta_{\cal G}} 
= (\pi_{\cal G})_{*} (D_{\theta_{\cal G}}) 
= (\pi_{\cal G})_{*} ((\mu^{\natural}_{\cal B})_{*} D_{\theta_E}) 
= (\nu_{\cal G})_{*} (p_{E*} D_{\theta_E}) 
= (\nu_{\cal G})_{*} E$ 
by the steps (B2) and (D0). 
\\ \ 
(C) 
By Lemma 6.1 (B), 
$(d_{\cal G} \beta) \mu^{\natural}_{\cal B}(\eta) 
= \mu^{\natural}_{\cal B}(\beta \eta)$ 
for $\beta \in {\rm Aut}(P_{E}, \theta_{E})$ and $\eta \in P_{E}$, 
so that 
$\delta_{\cal G} (d_{\cal G} \beta) \circ \hat{\mu}^{\natural}_{B_{\cal G}} 
= \hat{\mu}^{\natural}_{B_{\cal G}} \circ \beta$. 
Hence, 
$\delta_{\cal G} \circ d_{\cal G} 
= \widetilde{\mu}^{\natural}_{B_{\cal G}}$, 
which is a complex Lie group isomorphism (by (B2)). 
Especially, 
$\delta_{\cal G}$ is surjective. 
By Proposition 5.2 (E), 
$\delta_{\cal G}$ is injective, 
so that 
$d_{\cal G} = 
\delta_{\cal G}^{-1} \circ \widetilde{\mu}^{\natural}_{B_{\cal G}}$ 
is a complex Lie group isomorphism. 
Then 
$\delta_{\cal G} 
= \widetilde{\mu}^{\natural}_{B_{\cal G}} \circ d_{\cal G}^{-1}$ 
is a complex Lie group isomorphism. 
The last claim of the claim (C) is proved in Lemma 6.1 (B). 
\\ \
(D2) 
By the step (C) and Lemma 4.4 (B), 
$d_{\cal G} \circ {\cal P}$ 
is a complex Lie group isomorphism. 
For $\xi = p_{E}(\eta)$ and $\beta = {\cal P}(\alpha)$, 
$\nu_{\cal G}(\alpha \xi) 
= \nu_{\cal G}(p_{E} (\beta \eta)) 
= \pi_{\cal G}(\mu^{\natural}_{B_{\cal G}} (\beta \eta)) 
= \pi_{\cal G} ((d_{\cal G} \beta) \mu^{\natural}_{B_{\cal G}} (\eta)) 
= [d_{\cal G} \beta] \nu_{\cal G}(\xi)$ 
by (D0) and Lemma 6.1 (B).~\qed 
\smallskip \\ \
{\it Proof of Theorem 1.1}. 
Let $(M, E)$ be a connected compact c-manifold of ${\rm n}_{M} > 0$. 
Put ${\cal G} := a(P_{E}, \theta_{E})$. 
Assume that $L_{E}$ be immersive. 
By Lemma 4.6, $M$ is Fano, 
so that the connected $M$ is simply connected. 
By Lemma 6.0, 
$T_{\eta}' P_{E} = \{ X |_{\eta}; X \in {\cal G} \}$ 
at each $\eta \in P_{E}$. 
Because of Proposition 6.3,  
${\cal G}$ 
is a simple complex Lie algebra of $\infty > {\rm n}_{\cal G} > 1$ 
and that 
$(M, E) \cong (X_{\cal G}, E_{\theta_{\cal G}})$. 
By Lemma 4.1 (B, E(e)), 
${\cal G}$ is isomorphic to $\mathfrak{g} = a(M, E)$, 
so that $\mathfrak{g}$ 
is a simple complex Lie algebra of $\infty > {\rm n}_{\mathfrak{g}} > 1$. 
By Lemma 5.3, 
$(M, E) \cong (X_{\mathfrak{g}}, E_{\theta_{\mathfrak{g}}})$. 
Then 
$L_{E} \cong L_{E_{\theta_{\mathfrak{g}}}}$, 
which is very ample by Proposition 5.2 (D). 
By \cite{Ks1959}, 
$K_{M}^{*} \cong L_{E}^{\otimes {\rm k}_{M}}$ 
for ${\rm k}_{M} \in \mathbb{N}$ 
such as ${\rm n}_{M} = 2 {\rm k}_{M} - 1$. 
By Lemma 3.2, $K_{M}^{*}$ is very ample. 
~\qed 
\smallskip \\ \
{\it Proof of Theorem 1.2}. 
(1) 
Put $(M, E) := (X_{\mathfrak{g}}, E_{\theta_{\mathfrak{g}}})$. 
By Proposition 5.2 (E), 
$\delta_{\mathfrak{g}*}: 
\mathfrak{g} \to a(P_{e_{2}}, \theta_{\mathfrak{g}}); 
X \mapsto 
\left. (\frac{d}{dt} \right|_{t = 0} e^{t {\rm ad} X} |_{P_{e_{2}}})'$ 
is a monomorphism of complex Lie algebras. 
By Proposition 5.2 (A), 
the action of 
${\rm Aut}(\mathfrak{g})$ 
on $P_{e_{2}}$ is transitive, 
so is the action of ${\rm Ad}(\mathfrak{g})$, 
so that 
$T_{\eta}' P_{e_{2}} = \{ {\rm ad} X |_{\eta}; X \in \mathfrak{g} \}$ 
at each $\eta \in P_{e_{2}}$. 
As well as the proof Proposition 5.2 (E), 
${\rm ad} \mathfrak{g} |_{P_{e_{2}}} \cong {\rm ad} \mathfrak{g} 
\cong \mathfrak{g}$ 
because $\mathfrak{g}$ is a simple complex Lie algebra of 
${\rm n}_{\mathfrak{g}} > 1$. 
With respect to these canonical isomorphisms and Proposition 6.3 (C), 
$\delta_{\mathfrak{g}}: 
{\rm Aut}(\mathfrak{g}) \to {\rm Aut}(P_{e_{2}}, \theta_{\mathfrak{g}})$ 
is an isomorphism of complex Lie groups. 
By Propositions 4.3 (B) and 5.2 (C), 
$\kappa_{\theta_{\mathfrak{g}}}: 
(P_{e_{2}}, \theta_{\mathfrak{g}}) 
\to (P_{E}, \theta_{E})$ 
is an isomorphism, 
so that 
$\widetilde{\kappa}: 
{\rm Aut}(P_{e_{2}}, \theta_{\mathfrak{g}}) 
\to {\rm Aut}(P_{E}, \theta_{E}); 
\beta \mapsto 
\kappa_{\theta_{\mathfrak{g}}} \circ \beta \circ \kappa_{\theta_{\mathfrak{g}}}^{-1}$ 
is an isomorphism of complex Lie groups. 
By Lemma 4.4 (B), 
there exists an isomorphism of complex Lie groups 
${\cal P}^{-1}: {\rm Aut}(P_{E}, \theta_{E}) 
\to {\rm Aut}(M, E); \beta \mapsto \beta_{M}$ 
such as $\beta_{M} \circ p_{E} = p_{E} \circ \beta$. 
Then ${\cal P}^{-1} \circ \widetilde{\kappa} \circ \delta_{\mathfrak{g}}: 
{\rm Aut}(\mathfrak{g}) \to {\rm Aut}(M, E)$ 
is an isomorphism of complex Lie groups. 
Take any $\beta \in {\rm Aut}(\mathfrak{g})$. 
For any $\xi \in M = X_{\mathfrak{g}}$, 
there exists $\eta \in P_{E}$ 
such that $\xi = p_{E} (\eta)$. 
Then 
${\cal P}^{-1}(\widetilde{\kappa}(\delta_{\mathfrak{g}}(\beta))) \xi 
= {\cal P}^{-1}(\widetilde{\kappa}(\delta_{\mathfrak{g}}(\beta))) 
p_{E} (\eta) 
= p_{E}(\widetilde{\kappa}(\delta_{\mathfrak{g}}(\beta)) \eta) 
= p_{\mathfrak{g}}(\delta_{\mathfrak{g}}(\beta) \eta) 
= p_{\mathfrak{g}}(\beta \eta) 
= [\beta] \xi$. 
\\ \
(2) By the claim (1), 
\[
{\rm Aut}(M, E)/{\rm Aut}(M, E)^{\circ} 
\cong {\rm Aut}(\mathfrak{g}) / {\rm Aut}(\mathfrak{g})^{\circ} 
= {\rm Aut}(\mathfrak{g}) / {\rm Ad}(\mathfrak{g}), 
\]
which is isomorphic to the automorphism group of 
the Dynkin diagram of $\mathfrak{g}$ 
(e.g. \cite[p.218, (5.3.7)]{GG1978}), 
which gives the result as required. 
~\qed 
\smallskip \\ \
By Theorem 1.2 (1),
J.A.~Wolf's proof of \cite[(2.5)]{Wja1965} 
is improved as follows. 
\\ \
\textsc{Corollary 6.4}. 
(A) (cf. \cite[(2.5)]{Wja1965}). 
{\it 
Let $\mathfrak{g}$ be a simple complex Lie algebra 
of ${\rm n}_{\mathfrak{g}} > 1$ 
with $\mathfrak{g}' := a(X_{\mathfrak{g}}, E_{\theta_{\mathfrak{g}}})$. 
Then $\mathfrak{g}'$ 
is isomorphic to $\mathfrak{g}$ 
as complex Lie algebras. 
}
\\ \
(B) 
{\it 
Let 
$\mathfrak{g}$ and $\mathfrak{g}'$ 
be simple complex Lie algebras of 
${\rm n}_{\mathfrak{g}} > 1$ and ${\rm n}_{\mathfrak{g}'} > 1$. 
Then 
$(X_{\mathfrak{g}}, E_{\theta_{\mathfrak{g}}}) 
\cong (X_{\mathfrak{g}'}, E_{\theta_{\mathfrak{g}'}})$ 
if and only if 
$\mathfrak{g}$ and $\mathfrak{g}'$ 
are isomorphic as complex Lie algebras. 
}
\smallskip \\ \
{\it Proof}. 
(A) By Theorem 1.2 (1), 
${\rm Aut}(X_{\mathfrak{g}}, E_{\theta_{\mathfrak{g}}}) 
\cong {\rm Aut}(\mathfrak{g})$. 
Then $\mathfrak{g}' \cong \mathfrak{g}$. 
\\ \
(B) ``If''-part is proved by Lemma 5.3. 
Assume that 
$(X_{\mathfrak{g}}, E_{\theta_{\mathfrak{g}}}) 
\cong (X_{\mathfrak{g}'}, E_{\theta_{\mathfrak{g}'}})$. 
Then 
${\rm Aut}(X_{\mathfrak{g}}, E_{\theta_{\mathfrak{g}}}) 
\cong {\rm Aut}(X_{\mathfrak{g}'}, E_{\theta_{\mathfrak{g}'}})$. 
By Theorem 1.2 (1), 
${\rm Aut}(\mathfrak{g}) \cong {\rm Aut}(\mathfrak{g}')$. 
Then $\mathfrak{g} \cong \mathfrak{g}'$.~\qed 
\smallskip \\ \
{\bf 7. 
Proof of Corollary 1.3}. 
\smallskip \\ \
Let $M$ be a compact Riemannian manifold 
with the real Lie algebra $i(M)$ of all infinitesimal isometries on $M$. 
For a subbundle $S$ of a tensor bundle $T^{r}_{s}(M)$ of type $(r, s)$ on $M$, 
$i(M, S)$ denotes the Lie algebra of all $X \in i(M)$ 
stabilizing $S$ infinitesimally. 
\smallskip \\ \
\textsc{Lemma 7.0}. 
(A) 
{\it 
If $S$ is a parallel subbundle of $T^{r}_{s}(M)$, 
then the identity connected component $I^{o}(M)$ 
of the group $I(M)$ of all isometries on $M$ 
stabilizes $S$. 
Especially, 
$i(M) = i(M, S)$. 
}
\\ \
(B) 
{\it 
Let $M$ be a positive quaternionic-K\"ahler manifold with respect to 
the quaternion-K\"ahler structure $S$ of $M$ 
and the complex-contact twistor space $(Z, E)$. 
Then $i(M) = i(M, S)$, and that 
$a(Z, E)$ is isomorphic to $i(M) \otimes_{\mathbb{R}} \mathbb{C}$ 
as complex Lie algebras.
}
\smallskip \\ \
{\it Proof}. 
(A) (cf. \cite[p.248, Thm.VI.4.6 (1)]{KN1963}). 
Take an arbitary $X \in i(M)$, 
which is complete since $M$ is compact 
\cite[p.14, Prop.I.1.6]{KN1963}. 
Let $\alpha_{X}(t)$ 
be the 1-parameter group of isometries on $M$ 
generated by $X$. 
Let $x$ be an arbitrary point of $M$. 
Let $\tau$ be the orbit $x_{t} = \alpha_{X}(t) (x)$ 
of $x$.  
Let $\tau_{t}^{s}$ be the parallel displacement along the curve 
$\tau$ from $x_{s}$ to $x_{t}$. 
For each $t \in \mathbb{R}$, 
we consider a linear transformation 
$C_{t} := \tau_{0}^{t} \circ \alpha_{X}(t)_{*}$ of $T M |_{x}$. 
Then 
$C_{t}$ is a 1-parameter group 
of linear transformations of $T M |_{x}$ 
\cite[p.245, Prop.4.1]{KN1963}, 
which is contained in the linear holonomy group $\Psi(x)$ ($x \in M$) 
since $M$ is compact by B.~Kostant \cite[p.247, Thm.VI.4.5]{KN1963}. 
Let $C_{t}$ be extended to a 1-parameter group of automorphisms 
of $T^{r}_{s}(M) |_{x}$, which stabilizes $S |_{x}$ since $S$ is parallel. 
Then 
$\alpha_{X}(t)_{*}(S |_{x}) 
= \tau_{t}^{0} ( C_{t}(S |_{x})) 
= \tau_{t}^{0} ( S |_{x}) 
= S_{x_{t}} 
= S_{\alpha_{X}(t)(x)}$ 
for $t \in \mathbb{R}$. 
Hence, $i(M) = i(M, S)$. 
Then 
$I^{o}(M)$ stabilizes $S$ since $I^{o}(M)$ is generated as a group 
by a subset $\{ \alpha_{X}(t); t \in \mathbb{R}, X \in i(M) \}$. 
\\ \
(B) 
$S$ is a parallel subbundle of $End(T M) := T^{1}_{1}(M)$ 
(cf. \cite{Th1986}, \cite[p.146, \S 2]{NT1987}). 
From the claim (A), it follows that 
$i(M) = i(M, S)$ since $M$ is compact \cite[p.103]{Ss1989}. 
Because of \cite[Theorems 3.1, 3.2; p.159, Cor.2]{NT1987}, 
$a(Z, E) \cong i(M, S) \otimes_{\mathbb{R}} \mathbb{C} 
= i(M) \otimes_{\mathbb{R}} \mathbb{C}$. 
~\qed 
\smallskip \\ \
{\it Proof of Corollary 1.3}. 
By assumption, $L_{E}$ is immersive for $(Z, E)$. 
Put ${\cal G} := a(P_{E}, \theta_{E})$. 
By Lemma 6.0, 
$T_{\eta}' P_{E} = \{ X |_{\eta}; X \in {\cal G} \}$ 
at each $\eta \in P_{E}$. 
By Proposition 6.3 (B, D), 
${\cal G}$ is a simple complex Lie algebra of ${\rm n}_{\cal G} > 1$ 
and that $(Z, E) \cong (X_{\cal G}, E_{\theta_{\cal G}})$. 
Then 
${\rm n}_{X_{\cal G}} = {\rm n}_{Z} \geq 3$. 
By Lemma 5.2 (B), 
${\rm rank}({\cal G}) > 1$. 
By Lemmas 4.1 (B, E (e)) and 7.0 (B), 
${\cal G} \cong a(Z, E) \cong i(M) \otimes_{\mathbb{R}} \mathbb{C} = \mathfrak{g}$, 
so that ${\rm rank} (\mathfrak{g}) = {\rm rank}({\cal G}) > 1$. 
By Lemma 5.3, 
$(X_{\cal G}, E_{\theta_{\cal G}}) 
\cong (X_{\mathfrak{g}}, E_{\theta_{\mathfrak{g}}})$. 
Hence, 
$(Z, E) \cong (X_{\mathfrak{g}}, E_{\theta_{\mathfrak{g}}})$, 
so that 
${\rm Aut}(Z, E) \cong {\rm Aut}(X_{\mathfrak{g}}, E_{\theta_{\mathfrak{g}}}) 
\cong {\rm Aut}(\mathfrak{g})$ by Theorem 1.2 (1). 
~\qed 
\smallskip \\ \
{\it Acknowledgements}. 
The author is indebted to 
the reviewers of the last manuscripts 
for the critical and helpful comments, 
inspired by which this article is completed. 
He is also indebted to the following Professors: 
Masataka Tomari; Claude LeBrun, 
the late Masaru Takeuchi; Yasuyuki Nagatomo; 
Ryoichi Kobayashi, anonymous referees of JMSJ, 
Jaroslaw Buczy\'nski and Bogdan Alexandrov 
for pointing out crucial mistakes on the previous ones. 
He is cordially grateful to the following Professors: 
Nobuharu Onoda; Seiki Nishikawa; 
the late Tadashi Nagano; JMSJ and arXiv.org 
for their hospitality of encounter with these comments. 
He is also grateful to the following Professors: 
Hajime Urakawa, the late Tsunero Takahashi; 
Yusuke Sakane, Shigeru Mukai; 
Tomonori Noda; 
Hiroshi Asano, Hajime Kaji 
and Kiwamu Watanabe 
for their invaluable comments 
and encouradgements. 
During these periods, 
he was supported in part 
by the Grand-in-Aid for Encouragement of Young Scientists, 
Nos. 03740054, 05740052 and 06740055, 
Ministry of Education, Science and Culture; 
and by the Grand-in Aid for Scientific Research on Priority Areas 
(C)/(2), Nos. 10640046 and 15540066, 
Japan Society of the Promotion of Science. 
 
\end{document}